\documentclass[12pt,a4paper]{article}
\usepackage[utf8]{inputenc}
\usepackage[T1]{fontenc}
\usepackage{lmodern,amsmath,amsthm,amsfonts,amssymb,graphicx,float,microtype,thmtools,mathtools,verbatim}
\usepackage{tocloft}
\graphicspath{{figs/}}
\usepackage{paralist} 
\usepackage[shortlabels]{enumitem}
\setlist[itemize]{topsep=0ex,itemsep=0ex,parsep=0.3ex}
\setlist[enumerate]{topsep=0ex,itemsep=0ex,parsep=0.3ex}
\usepackage[table,usenames,dvipsnames,svgnames,table]{xcolor}
\usepackage{imakeidx}
\newcommand{\defn}[1]{\textcolor{Maroon}{\emph{#1}}\index{#1}}
\newcommand{\hdefn}[2]{\textcolor{Maroon}{\emph{#1-#2}}\index{#2@#1-#2}}
\newcommand{\mathdefn}[1]{\textcolor{Maroon}{#1}\index{$#1$}}
\newcommand{\DefNoIndex}[1]{\textcolor{Maroon}{#1}}

\makeindex
\usepackage[unicode=true]{hyperref}
\hypersetup{
colorlinks,
breaklinks=true,
linkcolor={blue!60!black},
citecolor={black},
urlcolor={blue!60!black},
pdftitle={Proof of the Clustered Hadwiger Conjecture}}
\usepackage[capitalise, compress, nameinlink, noabbrev]{cleveref}
\crefname{lem}{Lemma}{Lemmas}
\crefname{thm}{Theorem}{Theorems}
\crefname{cor}{Corollary}{Corollaries}
\crefname{prop}{Proposition}{Propositions}
\crefname{good}{}{}
\creflabelformat{good}{#2(g#1)#3}
\crefname{extra}{}{}
\creflabelformat{extra}{#2(x#1)#3}
\crefname{triple}{}{}
\creflabelformat{triple}{#2(t#1)#3}
\crefname{paren}{}{}
\creflabelformat{paren}{#2(#1)#3}
\crefname{tparen}{}{}
\creflabelformat{tparen}{#2(c#1)#3}
\crefformat{equation}{(#2#1#3)}
\Crefformat{equation}{Equation #2(#1)#3}

\setlist[enumerate,2]{label=(\roman*),ref=(\roman{enumi}.\roman*)}

\newcommand{\authorreveal}[1]{}

\newenvironment{clmproof}{\noindent\textit{Proof.}}{\hfill$\diamondsuit$}

\usepackage[longnamesfirst,numbers,sort&compress]{natbib}
\makeatletter
\def\NAT@spacechar{~}
\makeatother
\usepackage[margin=30mm]{geometry}
\renewcommand{\baselinestretch}{1.05}
\allowdisplaybreaks

\DeclarePairedDelimiter{\floor}{\lfloor}{\rfloor}
\DeclarePairedDelimiter{\ceil}{\lceil}{\rceil}

\renewcommand{\epsilon}{\varepsilon}
\renewcommand{\emptyset}{\varnothing}
\renewcommand{\ge}{\geqslant}
\renewcommand{\le}{\leqslant}
\renewcommand{\geq}{\geqslant}
\renewcommand{\leq}{\leqslant}

\DeclareMathOperator{\dist}{dist}
\DeclareMathOperator{\tw}{tw}

\DeclareMathOperator{\interior}{Int}
\newcommand{\RR}{\mathbb{R}}
\newcommand{\JJ}{\mathcal{J}}
\newcommand{\PP}{\mathcal{P}}

\newcommand{\TT}{\mathcal{T}}
\newcommand{\GG}{\mathcal{G}}

\newcommand{\LL}{\mathcal{L}}
\newcommand{\NN}{\mathbb{N}}

\newcommand{\torso}[2]{{#1}\langle{#2}\rangle}
\newcommand{\ttorso}[3]{\torso{\torso{#1}{#2}}{#3}}
\newcommand{\ltorso}[2]{{#1}\{{#2}\}}

\DeclareMathOperator{\rtop}{top}
\newcommand{\ttop}[2]{\rtop(\torso{#1}{#2})}

\theoremstyle{plain}
\newtheorem{thm}{Theorem}
\newtheorem{lem}[thm]{Lemma}
\newtheorem{cor}[thm]{Corollary}

\newtheorem{prop}[thm]{Proposition}
\newtheorem{obs}[thm]{Observation}
\newtheorem{clm}{Claim}[thm]
\crefname{obs}{Observation}{Observations}

\newtheorem*{lem*}{lem}
\theoremstyle{definition}

\newtheorem*{conj*}{Conjecture}

\begin{document}
\title{\bf\boldmath\fontsize{18pt}{18pt}\selectfont
Proof of the Clustered Hadwiger Conjecture\footnote{ Mathematics Subject Classification: 05C15, 05C83 \\
An extended abstract of this paper appeared in
\href{https://doi.org/10.1109/FOCS57990.2023.00116}{\emph{Proc. 64th IEEE Symposium on Foundations of Computer Science}} (FOCS '23), pp.~1921--1930.}}

\author{%
Vida Dujmovi{\'c}\,\footnotemark[3]\quad
Louis Esperet\,\footnotemark[4] \quad
Pat Morin\,\footnotemark[5] \quad
David~R.~Wood\,\footnotemark[2]
}

\maketitle

\begin{abstract}
  Hadwiger's Conjecture asserts that every $K_h$-minor-free graph is properly $(h-1)$-colourable. We prove the following improper analogue of Hadwiger's Conjecture: for fixed $h$, every $K_h$-minor-free graph is $(h-1)$-colourable with monochromatic components of bounded size. The number of colours is best possible regardless of the size of the monochromatic components. This solves an open problem of Edwards, Kang, Kim, Oum and Seymour [\emph{SIAM J. Discrete Math.} 2015], and concludes a line of research initiated in 2007. Similarly, for fixed $t\geq s$, we show that every $K_{s,t}$-minor-free graph is $(s+1)$-colourable with monochromatic components of bounded size. The number of colours is best possible, solving an open problem of van den Heuvel and Wood [\emph{J.~London Math.\ Soc.} 2018].  We actually prove a single theorem from which both of the above results are immediate corollaries. For an excluded apex minor, we strengthen the result as follows: for fixed $t\geq s\geq 3$ and for any fixed apex graph $X$, every $K_{s,t}$-subgraph-free $X$-minor-free graph is $(s+1)$-colourable with monochromatic components of bounded size. The number of colours is again best possible.
\end{abstract}

\renewcommand{\thefootnote}{\fnsymbol{footnote}}

\footnotetext[3]{School of Computer Science and Electrical Engineering, University of Ottawa, Ottawa, Canada (\texttt{vida.dujmovic@uottawa.ca}). Research supported by NSERC.}

\footnotetext[4]{Laboratoire G-SCOP  (CNRS, Univ.\ Grenoble Alpes), Grenoble, France
(\texttt{louis.esperet@grenoble-inp.fr}). Partially supported by ANR Projects GrR
(ANR-18-CE40-0032) and Twin-width (ANR-21-CE48-0014), and by LabEx PERSYVAL-lab (ANR-11-LABX-0025).}

\footnotetext[5]{School of Computer Science, Carleton University, Ottawa, Canada (\texttt{morin@scs.carleton.ca}). Research  supported by NSERC.}

\footnotetext[2]{School of Mathematics, Monash University, Melbourne, Australia (\texttt{david.wood@monash.edu}). Research supported by the Australian Research Council and by NSERC.}

\renewcommand{\thefootnote}{\arabic{footnote}}

\newpage
\tableofcontents
\newpage

\section{Introduction}
\label{Introduction}

\subsection{Hadwiger's Conjecture}

Our starting point is Hadwiger's Conjecture~\citep{Hadwiger43}, which suggests a deep relationship between graph colourings and graph minors\footnote{We consider simple, finite, undirected graphs~$G$ with vertex-set~${V(G)}$ and edge-set~${E(G)}$. See later sections and \citep{Diestel5} for graph-theoretic definitions not given here.}. A \defn{colouring} of a graph $G$ is a function that assigns one colour to each vertex of $G$. For an integer $k\geq 1$, a \hdefn{$k$}{colouring} is a colouring using at most $k$ colours. A colouring of a graph is \defn{proper} if each pair of adjacent vertices receives distinct colours. The \defn{chromatic number} \defn{$\chi(G)$} of a graph $G$ is the minimum integer $k$ such that $G$ has a proper $k$-colouring. A graph $H$ is a \defn{minor} of a graph $G$ if $H$ is isomorphic to a graph that can be obtained from a subgraph of $G$ by contracting edges. A graph~$G$ is \hdefn{$H$}{minor-free} if~$H$ is not a minor of~$G$. Let $K_h$ be the complete graph on $h$ vertices. \citet{Hadwiger43} famously conjectured that every $K_h$-minor-free graph is properly $(h-1)$-colourable. This is widely considered to be one of the most important open problems in graph theory. By Wagner's characterization of $K_5$-minor-free graphs~\citep{Wagner37}, the case $h=5$ is equivalent to the 4-Colour Theorem~\citep{AH89,RSST97}. The conjecture is true for $h\leq 6$ \citep{RST-Comb93} and is open for $h\geq 7$. The best known upper bound on the chromatic number of $K_h$-minor-free graphs remained of order $O(h\sqrt{\log h})$~\citep{Kostochka82,Kostochka84,Thomason84,Thomason01} from the 1980s until a breakthrough by \citet{NPS23}, who proved that $K_h$-minor-free graphs are $O(h(\log h)^{1/4+\epsilon})$-colourable (for any fixed $\epsilon>0$). In another breakthrough, this bound was subsequently improved to $O(h\log\log h)$ by \citet{DP25}. It is open whether every $K_h$-minor-free graph is $O(h)$-colourable. See Seymour's survey~\citep{SeymourHC} for more on Hadwiger's Conjecture.

\subsection{The Clustered Hadwiger Conjecture}\label{sec:hdw}

As mentioned above, one of the main ways to approach Hadwiger's Conjecture has been to try to minimise the number of colours in a proper colouring of a $K_h$-minor-free graph. A second natural approach is to fix a number of colours close to Hadwiger's bound (at $h-1$, for instance), and try to obtain a colouring that is close to being proper. This leads to the notion of improper colouring of $K_h$-minor-free graphs; see~\citep{OOW19,Liu24,KM07,Norin15,vdHW18,LW1,LW2,LW3,LW4,Kawa08,Wood10,EKKOS15,LO18,KO19,DN17,NSSW19,NSW22} and the references therein. A \defn{monochromatic component} with respect to a colouring of a graph $G$ is a connected component of the subgraph of $G$ induced by all the vertices assigned a single colour. A colouring has \defn{clustering} $c$ if every monochromatic component has at most $c$ vertices. Note that a colouring with clustering $1$ is precisely a proper colouring. The \defn{clustered chromatic number} \defn{$\chi_\star(\GG)$} of a graph class $\mathcal{G}$ is the minimum integer $k$ for which there exists an integer $c$ such that every graph in  $\GG$ is $k$-colourable with clustering $c$ (or $\infty$ if no such $k$ exists). See \citep{WoodSurvey} for an extensive survey on this topic, and see \citep{BBEGLPS24,EW23} for connections between clustered colouring and asymptotic dimension in geometric group theory and site percolation in probability theory. One of the earliest works on clustered colouring, by \citet{KMRV97}, used it as a tool to design algorithms for evolving databases.

Consider the clustered chromatic number of the class of $K_h$-minor-free graphs. So-called standard examples provide a lower bound of $h-1$ (regardless of the clustering value); see \cref{Standard} below. A natural weakening of Hadwiger's Conjecture, sometimes called the Clustered Hadwiger Conjecture, asserts that every $K_h$-minor-free graph is $(h-1)$-colourable with clustering at most $f(h)$ for some function $f$. This conjecture was first asked as an open problem by  \citet{EKKOS15}. This line of research was initiated in 2007 by \citet{KM07}, who proved the first $O(h)$ upper bound on the number of colours. Their bound was $\ceil{\frac{31}{2}h}$, which was successively improved to $\ceil{\frac{7h-3}{2}}$ by \citet{Wood10}\footnote{The result of \citet{Wood10} depended on a result announced by Norin and Thomas~\citep{Thomas09}, which has not yet been fully written.}, to $4h-4$ by \citet{EKKOS15}, to $3h-3$ by \citet{LO18}, to $2h-2$ independently by \citet{Norin15}, \citet{vdHW18} and \citet{DN17}, and most recently to $h$ by \citet{LW3}. Note that \citet{DN17} showed that $h-1$ colours suffice for $h\leq 9$, and that \citet{EKKOS15} proved that every $K_h$-minor-free graph has an $(h-1)$-colouring in which each monochromatic component has bounded maximum degree (which is significantly weaker than having bounded size).

Our first contribution is to prove the Clustered Hadwiger Conjecture,\footnote{\citet{DN17} announced in 2017 that a forthcoming paper, which has not yet been fully written, will also prove the Clustered Hadwiger Conjecture.} thereby solving the above-mentioned open problem of \citet{EKKOS15}.

\begin{thm}
\label{Kh}
There is a function $f$ such that every $K_h$-minor-free graph is $(h-1)$-colourable with clustering at most $f(h)$.
\end{thm}

\subsection{Excluding a Complete Bipartite Minor}

Consider the clustered chromatic number of the class of $K_{s,t}$-minor-free graphs, where $K_{s,t}$ is the complete bipartite graph with parts of size $s\geq 1$ and $t\geq s$.
This question is of particular interest since the answer turns out not to depend on $t$. Van den Heuvel and Wood~\citep{vdHW18} proved a lower bound of $s+1$ (for $t\geq\max\{s,3\}$; see \cref{Standard} below), and observed that results of \citet{EKKOS15} and \citet{OOW19} imply
an upper bound of $3s$, which was improved to $2s+2$ by \citet{DN17},
and further improved to $s+2$ by \citet{LW2}. Our second main result resolves the question.

\begin{thm}
\label{Kst}
There is a function $f$ such that every $K_{s,t}$-minor-free graph is $(s+1)$-colourable with clustering at most $f(s,t)$.
\end{thm}

Combined with the above-mentioned lower bound, \cref{Kst} shows that the clustered chromatic number of the class of $K_{s,t}$-minor-free graphs equals $s+1$ for $t\geq\max\{s,3\}$. This resolves an open problem posed by \citet{vdHW18}.

\subsection{Colin de Verdi\`ere Parameter}
\label{ColindeVerdiere}

The Colin de Verdi\`ere parameter $\mu(G)$ is an important graph invariant introduced by \citet{CdV90,CdV93}; see~\citep{HLS,Schrijver97} for surveys. It is known that $\mu(G)\leq 1$ if and only if $G$ is a disjoint union of paths, $\mu(G)\leq 2$ if and only if $G$ is outerplanar, $\mu(G)\leq 3$ if and only if $G$ is planar, and $\mu(G)\leq 4$ if and only if $G$ is linklessly embeddable. A famous conjecture of \citet{CdV90} states that $\chi(G)\leq \mu(G)+1$, which implies the 4-Colour Theorem and is implied by Hadwiger's Conjecture. It was conjectured by \citet{WoodSurvey} that the clustered chromatic number of the class of graphs with Colin de Verdi\`ere parameter at most $\mu$ equals $\mu+1$. Since any graph with Colin de Verdi\`ere parameter $\mu$ excludes the complete graph $K_{\mu+2}$ as a minor, all partial results on the Clustered Hadwiger Conjecture mentioned in \cref{sec:hdw} naturally give partial results on this conjecture. In particular, the result of \citet{LW3} proving that graphs with no $K_h$-minor are $h$-colourable with clustering $f(h)$ for some function $f$ implies a weaker version of the conjecture with $\mu+2$ colours instead of $\mu+1$. The same result is also implied by an earlier theorem of \citet{LW2} on graphs excluding a minor and a complete bipartite subgraph (see \cref{sec:Kstintro} for more details about this result).

We resolve Wood's conjecture \citep{WoodSurvey}.

\begin{thm}
\label{ClusteredColin}
The clustered chromatic number of the class of graphs with Colin de Verdi\`ere parameter $\mu$ equals $\mu+1$.
\end{thm}

The lower bound in \cref{ClusteredColin} is proved in \citep{WoodSurvey}. The upper bound follows immediately from either \cref{Kh} or \cref{Kst}, since graphs with Colin de Verdi\`ere parameter $\mu$ are $K_{\mu+2}$-minor-free \citep{CdV90,CdV93} and are $K_{\mu,\max\{\mu,3\}}$-minor-free~\citep{HLS}.

\subsection{A Common Generalization}
\label{Generalization}

In this work, we prove a common generalization of \cref{Kh,Kst} (and thus \cref{ClusteredColin}) using the following definition. Let \defn{$\JJ_{s,t}$} be the set of all graphs $K_s\oplus T$ where $T$ is a $t$-vertex tree. Here the \defn{complete join} $G_1\oplus G_2$ of $G_1$ and $G_2$ is the graph obtained from the disjoint union of graphs $G_1$ and $G_2$ by adding all edges between $G_1$ and $G_2$. A graph is \defn{$\JJ_{s,t}$-minor-free} if it contains no graph in $\JJ_{s,t}$ as a minor\footnote{Let \defn{$P_{s,t}$} be the graph obtained from $K_{s,t}$ by adding a path on the $t$-vertex side. Every tree on $t^2$ vertices contains a $K_{1,t}$ minor or a path on $t$ vertices. It follows that every graph in $\JJ_{s,t^2}$ contains $K_{s+1,t}$ or $P_{s,t}$ as a minor. Conversely, $P_{s,t+s}$ and $K_{s+1,t+s}$ contain a graph in $\JJ_{s,t}$ as a minor (obtained by contracting a suitable $s$-edge matching). This says that (ignoring dependence on $t$) being $\JJ_{s,t}$-minor-free is equivalent to being $K_{s+1,t}$-minor-free and $P_{s,t}$-minor-free. Therefore \cref{Jst} could be stated for graphs that are $K_{s+1,t}$-minor-free and $P_{s,t}$-minor-free. We choose to work with $\JJ_{s,t}$ for convenience.}.

\begin{thm}
\label{Jst}
There is a function $f$ such that every $\JJ_{s,t}$-minor-free graph is $(s+1)$-colourable with clustering at most $f(s,t)$.
\end{thm}

Since $K_{s,t}$ is a subgraph of each graph in $\JJ_{s,t}$, \cref{Jst} implies \cref{Kst}. Since $\JJ_{h-2,2}=\{K_h\}$, \cref{Jst} with $s=h-2$ implies \cref{Kh}.

\subsection{Excluding a Complete Bipartite Subgraph}\label{sec:Kstintro}

As mentioned above, \citet{LW2} proved that $K_{s,t}$-minor-free graphs are $(s+2)$-colourable with bounded clustering. They actually proved a stronger result, where the number of colours is determined by an excluded complete bipartite subgraph, as expressed in the following results. For a graph $X$, a graph $G$ is \hdefn{$X$}{subgraph-free} if no subgraph of $G$ is isomorphic to $X$.

\begin{thm}[{\protect\citep[Theorem~5]{LW2}}]
\label{LiuWoodHplanar}
There is a function $f$ such that for any integers $s, t, w \in \NN$, every $K_{s,t}$-subgraph-free graph of treewidth at most $w$ is $(s+1)$-colourable with clustering at most $f(s,t,w)$.
\end{thm}

We use an extension of \cref{LiuWoodHplanar} as one of the tools in the current paper (see \cref{lem:LW}).

\begin{thm}[{\protect\citep[Theorem~2]{LW2}}]
\label{LiuWoodH}
There is a function $f$ such that for any integers $t\geq s\geq 1$ and for any graph $X$, every $K_{s,t}$-subgraph-free $X$-minor-free graph is $(s+2)$-colourable with clustering at most $f(s,t,X)$.
\end{thm}

In both these theorems, the number of colours is best possible~\citep{LW2}. It is worth contrasting our \cref{Kst} with \cref{LiuWoodH}.  \cref{Kst} uses only $s+1$ colours, but applies only to a graph $G$ that contains no $K_{s,t}$ minor.  \cref{LiuWoodH} uses $s+2$ colours, but applies in the more general setting in which $G$ contains no $K_{s,t}$-subgraph and excludes any other (fixed) graph $X$ as a minor.

\cref{LiuWoodH} implies that $K_h$-minor-free graphs are $(h+1)$-colourable with bounded clustering (since $K_h$ is a minor of $K_{h-1,h-1}$). \citet{LW3} pushed this proof further to reduce the number of colours to $h$, as mentioned above. These results are presented over a series of three articles \citep{LW1,LW2,LW3}. The main tool introduced in the first article of the series \citep{LW1} shows (via a technical 70-page proof) that $K_{s,t}$-subgraph-free graphs of bounded layered treewidth are $(s+2)$-colourable with bounded clustering. Three colours is best possible in the $s=1$ case, but for $s\geq 2$, it is open whether $s+1$ colours (which would be best possible) suffice in this setting~\citep{LW1,LW2}. Bounded layered treewidth is not a minor-closed property. However, for minor-closed classes, having bounded layered treewidth is equivalent to excluding an apex graph as a minor \citep{DMW17}. Here a graph is \defn{apex} if it can be made planar by deleting at most one vertex. Our next result shows that $s+1$ colours suffice for apex-minor-free $K_{s,t}$-subgraph-free graphs.

\begin{thm}
\label{apex1}
There is a function $f$ such that for any integers $t\geq s\geq 3$ and for any apex graph $X$, every $K_{s,t}$-subgraph-free $X$-minor-free graph is $(s+1)$-colourable with clustering at most $f(s,t,X)$.
\end{thm}

Several notes on \cref{apex1} are in order:
\begin{itemize}
\item \cref{LiuWoodHplanar} is equivalent to saying that $K_{s,t}$-subgraph-free graphs excluding a fixed planar graph as a minor are $(s+1)$-colourable with bounded clustering. \cref{apex1} strengthens this result to the setting of apex-minor-free graphs (for $s\geq 3$).
\item The bound on the number of colours in \cref{apex1} is tight, simply because $s+1$ colours is tight for bounded treewidth graphs~\citep{LW1,LW2}.
\item Some non-trivial lower bound on $s$ is needed in \cref{apex1}, since the hexagonal grid graph is $K_5$-minor-free and $K_{1,7}$-subgraph-free, but every $2$-colouring has unbounded clustering by the Hex Lemma~\citep{HT19}. \cref{apex1} with $s=2$ is open, although recently solved in the case of graphs of bounded Euler genus by \citet{Dvorak25}.
\item Our proof of \cref{apex1} in the case $s\ge 4$ is reasonably short and simple, and is presented in \cref{apex1proof}.  The case $s=3$ is more difficult, and requires tools for dealing with $K_{2,t}$-subgraphs in surfaces that are also required by the proofs of \cref{Kh,Kst,Jst}.
\end{itemize}

For the sake of clarity, we now summarise how the present paper uses results of \citet{LW1,LW2,LW3}. \cref{lem:LW0} is a result of \citet{LW2} for colouring $K_{s,t}$-subgraph-free graphs of bounded treewidth. We use this to provide a simple proof of $(h+4)$-colourability for $K_h$-minor-free graphs as a way to introduce some of the key ideas used in our main proof. \cref{lem:LW} extends this lemma, and we provide a full proof that uses one lemma by \citet{LW1} (our \cref{heavy_neighbour_bound}), which has a simple one-paragraph proof. \cref{lem:LW2} is another result for colouring bounded treewidth graphs that extends a result of \citet{LW3}. Again, we provide a full proof. The present paper does not use layered treewidth or the 70-page proof of \citet{LW1} mentioned above.

\subsection{Clustered 4-Colouring Theorems}

It is well-known that the clustered chromatic number of the class of planar graphs equals $4$, where the upper bound follows from the 4-Colour Theorem or a weaker result of \citet{CCW86}, and the lower bound follows from the $s=3$ case of \cref{Standard} below. In fact, much more general results are known for $4$-colouring with bounded clustering.

The \defn{Euler genus} of a surface with~$h$ handles and~$c$ cross-caps is~${2h+c}$. The \defn{Euler genus} of a graph~$G$ is the minimum integer $g\geq 0$ such that there is an embedding of~$G$ in a surface of Euler genus~$g$; see \citep{MoharThom} for more about graph embeddings in surfaces. \citet{DN17} proved the following elegant generalization of the 4-Colour Theorem.

\begin{thm}[\citep{DN17}]\label{Genusg}
There is a function $f$ such that every graph of Euler genus $g$ has a $4$-colouring with clustering at most $f(g)$.
\end{thm}

It follows from Euler's formula that graphs of Euler genus $g$ are $K_{3,2g+3}$-minor-free. Thus the $s=3$ case of \cref{Kst} generalises \cref{Genusg} to the setting of $K_{3,t}$-minor-free graphs. This setting is substantially more general, since disjoint unions of $K_5$ are $K_{3,3}$-minor-free, but have unbounded Euler genus. This highlights the utility of considering excluded complete bipartite minors.

\cref{apex1} with $s=3$ provides a further generalization where the $K_{3,t}$-minor-free assumption is relaxed to apex-minor-free and $K_{3,t}$-subgraph-free.
Again, this is a substantial generalization since, for example, the $1$-subdivision of $K_{3,n}$ is $K_5$-minor-free and $K_{3,3}$-subgraph-free, but contains a $K_{3,n}$ minor.

Dvo{\v{r}}{\'a}k and Norin's proof of \cref{Genusg} uses the so-called `island' method. Our proof uses and builds on this approach; see \cref{three_islands} below.

\subsection{Lower Bounds}

The number of colours in \cref{Kh,Kst,apex1} is best possible because of the following well-known `standard' example~\citep{LMST08,EKKOS15,WoodSurvey,LW2,vdHW18}, which we include for completeness.

\begin{prop}
\label{Standard}
For each integer $s\geq 1$, there is a graph class $\GG_s$ such that every graph in $\GG_s$ has treewidth $s$ and is $K_{s,s+2}$-subgraph-free, and for any integer $c\geq 1$, there exists $G\in\GG_s$ such that every $s$-colouring of $G$ has a monochromatic component on more than $c$ vertices.
\end{prop}

\begin{proof}
We proceed by induction on $s$, where $\GG_1$ is the class of all paths on at least two vertices, which trivially satisfies the desired properties. Now assume the claim for some $s\geq 1$. As illustrated in
\cref{LowerBoundGraphs}, let $\GG_{s+1}$ be the class of all graphs obtained by taking arbitrarily many disjoint copies of graphs in $\GG_s$ and adding one dominant vertex.

\begin{figure}[!ht]
  \begin{center}
    \begin{tabular}{ccc}
      \includegraphics[page=1,trim={1cm 0 .8cm 0},clip,scale=1.1]{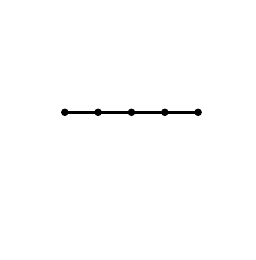} &
      \includegraphics[page=2,scale=1.1]{g123} &
      \includegraphics[page=3,scale=1.1]{g123}
    \end{tabular}
  \end{center}
  \caption{Members of the families $\GG_1$, $\GG_2$, and $\GG_3$.}
  \label{LowerBoundGraphs}
\end{figure}

Every graph in $\GG_{s+1}$ has treewidth $s+1$ and is $K_{s+1,s+3}$-subgraph-free. For any $c$, there is a graph $H$ in $\GG_s$ such that any $s$-colouring of $H$ has a monochromatic component on more than $c$ vertices. Let $G$ be obtained from $c$ disjoint copies of $H$ by adding a dominant vertex $v$. So $G\in\GG_{s+1}$. In any colouring of $G$ with clustering $c$, some copy of $H$ avoids the colour assigned to $v$, as otherwise the monochromatic component containing $v$ would have at least $c+1$ vertices. By assumption, $H$ has at least $s+1$ colours, implying $G$ has at least $s+2$ colours. This proves the claim.
\end{proof}

\cref{Standard} shows that $h-1$ colours in \cref{Kh} is best possible, since every graph with treewidth $s$ is $K_{s+2}$-minor-free. Similarly, \cref{Standard} shows that $s+1$ colours in \cref{Kst} is best possible, since the graph in \cref{Standard} is $K_{s,\max\{s,3\}}$-minor-free~\citep{vdHW18}. To see that the number of colours in \cref{apex1} is best possible, note that every graph in $\GG_s$ contains no $(s+1)\times(s+1)$ grid as a minor (since the $(s+1)\times(s+1)$ grid has treewidth $s+1$). Thus $\GG_s$ excludes a planar (and thus apex) graph as a minor. If, for any function $f$, every graph in $\GG_s$ is $k$-colourable with clustering at most $f(s)$, then $k\geq s+1$ by \cref{Standard}.

\subsection{Proof Outline}

This section gives a high-level sketch of the proof of \cref{Jst}. We highlight the main challenges that arise when adapting existing techniques to prove our results, and give a rough idea of how we overcome these challenges.
The key to our proofs is the novel use of `graph product structure theory' in partnership with the Graph Minor Structure Theorem of \citet{RS-XVI}. Graph product structure theory is a recently developed field that describes graphs in a complicated graph class as subgraphs of a product of simpler graphs (along with some other operations).

The starting point for these recent developments is the Planar Graph Product Structure Theorem of \citet{DJMMUW20}, which says that every planar graph $G$ is isomorphic to a subgraph of the strong product\footnote{The \defn{strong product} of graphs $A$ and $B$ is the graph $\mathdefn{A\boxtimes B}$ with vertex set $V(A)\times V(B)$, where distinct vertices $(x,y)$ and $(x',y')$ are adjacent if and only if $x=x'$ or $xx'\in E(A)$, and $y=y'$ or $yy'\in E(B)$.} of a graph $H$ of treewidth $8$ and a path $P$.  This often allows a question about planar graphs ($G$) to be reduced to an analogous question about graphs of bounded treewidth ($H$), which is usually easier to solve. Extending the solution for $H$ to the product $H\boxtimes P$ is often straightforward.
This approach has been used to solve several decades-old problems in mathematics and theoretical computer science, including bounding the queue-number of planar graphs~\citep{DJMMUW20}, bounding the nonrepetitive chromatic number of planar graphs~\citep{DEJWW20},
the construction of nearly optimal adjacency labelling schemes for planar graphs~\citep{DEGJMM21,EJM23}, amongst other results~\citep{BKW26,BDJM}.

This approach is, of course, limited to graph families that admit this type of product structure. \citet{DJMMUW20} proved that apex-minor-free graphs are the most general minor-closed classes contained in the product of a bounded treewidth graph and a path, so this reduction cannot be applied to the graph classes considered by \cref{Kh,Kst,Jst} (except for a few small values of $h$ or $s$).

Since not all minor-closed classes admit the type of product structure described above, we start by applying the Graph Minor Structure Theorem of \citet{RS-XVI}, which describes graphs excluding a fixed minor in terms of graphs called `torsos' that are formed from a surface-embedded subgraph, by adding vortices in the surface, and apex vertices with unrestricted neighbourhoods. The torsos are then pasted together in a tree-like way (described by a tree-decomposition) using clique-sums.

A natural strategy for colouring a $\JJ_{s,t}$-minor-free graph is to colour each torso one-by-one beginning with the root torso. With this strategy, when colouring a subsequent non-root torso, the vertices it shares with its parent torso are already precoloured with colours assigned in previous steps. To avoid creating monochromatic components that span multiple torsos, it is desirable that the colouring of the new torso \defn{properly extends} the precolouring, so that each vertex receives a colour distinct from those of all its precoloured neighbours.  In our setting, the precoloured vertices of a torso are contained in one clique-sum whose size is bounded by a function of $s$ and $t$. This bound may be much larger than $s$, which is a major obstacle since \cref{Jst} promises a colouring with at most $s+1$ colours, and it is impossible to properly extend a precolouring if some vertex is adjacent to $s+1$ precoloured vertices that use all $s+1$ colours.

For the purposes of this high-level description, call a clique-sum \defn{small} if it has at most $s$ vertices, and \defn{large} otherwise (although still of size bounded by a function of $s$ and $t$). Using the approach described above, the precoloured vertices of torsos are part of either a small or a large clique-sum. Since the natural strategy is not applicable to large clique-sums, we perform those clique-sums before any colouring. In particular, we call a maximal collection of torsos pasted together by only large clique-sums a \defn{curtain}. (This choice of name will become clear later.) Now one can think of the Robertson--Seymour decomposition as a collection of curtains glued together using small clique-sums, which we call a \defn{tree of curtains}. To colour the entire graph, it suffices to colour each curtain one-by-one beginning with the curtain that contains the root torso.  The critical advantage of proceeding this way is that each curtain comes with at most $s$ precoloured vertices.
 
To colour each individual curtain, we use a variant of the Graph Minor Structure Theorem by \citet{DvoTho}, which classifies apex vertices as either major or non-major. Edges incident to major apex vertices are unrestricted, but non-major apex vertices are not adjacent to vertices in the part of the surface-embedded subgraph avoiding the vortices. A graph is \hdefn{$k$}{apex} if it has $k$ vertices whose removal leaves a planar graph. \citet{DvoTho} showed that for $k$-apex-minor-free graphs, the number of major apex vertices is at most $k-1$. In our case, $\JJ_{s,t}$ includes an $(s-2)$-apex graph. So the number of major apex vertices is at most $s-3$. From now on, we make no distinction between non-major apex vertices and vortex vertices.

An important consequence of the \citet{DvoTho} result is that (after some manipulation) any large clique-sum only involves (major and non-major) apex vertices and vertices from vortices of the two torsos being summed. Therefore, none of the clique-sums used to make a curtain touch the surface-embedded parts of the torsos, except at the boundaries of vortices.  In our language, each of these clique-sums connects the \defn{top} of one torso to the \defn{near-top} of the other.

The material discussed above is presented in \cref{Structure}. To illustrate the utility of these preliminary ideas, \cref{Structure} also describes a short proof that $K_h$-minor-free graphs have clustered colourings using $h+4$ colours, and we prove \cref{apex1} for $s\ge 4$.

Our task now is to colour each curtain, given a set of at most $s$ precoloured vertices that are at the top of the curtain.  Recall that each torso in the curtain is described by a set of at most $s-3$ major apex vertices, a surface-embedded subgraph, and a collection of vortices and non-major apex vertices.

The next step employs graph product structure theory. A result of \citet{DJMMUW20} says that if we remove the major apex vertices from a torso, the resulting graph is a subgraph of the strong product of a bounded treewidth graph and a path (generalising the result for planar graphs mentioned above). This product structure can be described in terms of partitions and layerings. For each torso (after removing its major apex vertices), we obtain a partition of the vertex-set with connected parts, and a layering such that the intersection of each part and each layer has bounded size. Loosely speaking, each part in the partition is `skinny' with respect to the layering.  In addition, the first layer of the layering contains exactly the vortices and the non-major apex vertices of the torso, while the second layer includes the surface-embedded neighbours of the vortex boundaries. This means that each large clique-sum used to paste together a pair of torsos in the curtain $G$ is restricted to the major apex vertices, the first layer of one torso, and the first two layers of the other torso.  Thus, each large clique-sum takes place within the first two layers of the layerings of the respective torsos, connecting the top of one to the near-top of the other.

\citet{DJMMUW20} proved that the minor obtained by contracting each part of the partition results in a quotient graph that has bounded treewidth.  However, we do not perform the contractions that would create this bounded treewidth quotient graph, since doing so could interfere with a clique-sum and could introduce complete bipartite subgraphs that we wish to avoid for reasons that will become clear shortly. Instead, we reduce to a bounded treewidth graph by `raising the curtain' as follows. Within each torso of a curtain $G$, we contract the parts of the corresponding partition while avoiding the first five layers.  This produces a minor $G_{\uparrow}$ of $G$ on at most six layers and having bounded treewidth. Since these contractions avoid the first five layers of each torso, they avoid the vertices used in the clique-sums between torsos.  Therefore, these clique-sums remain in $G_\uparrow$, and  $G_\uparrow$ is also a curtain with a surface-embedded subgraph within each torso.  Unlike $G$, $G_{\uparrow}$ has bounded treewidth.

However, it is not enough that $G_{\uparrow}$ has bounded treewidth. Indeed, the extremal examples have bounded treewidth (see \cref{Standard}). It is also critical that $G_{\uparrow}$ is a minor of $G$. However, because of the possibility of edges added to each torso in order to perform a clique-sum, it is not immediate that $G_{\uparrow}$ is obtained by contracting connected subgraphs in $G$. To overcome this issue, we show that the tree of curtains we construct is `lower-minor-closed'. We expect that this property, and in particular \cref{ApexMinorFreeStructure}, may be of independent interest.

In \cref{BoundedTreewidth}, we prove two lemmas (that are technical extensions of results by \citet{LW2,LW3}) for colouring bounded treewidth graphs containing no $K_{s,t}$-subgraph. Our goal is to apply these results to $G_{\uparrow}$. However, without further work, $K_{s,t}$-subgraphs may be present in $G_{\uparrow}$.

Any such problematic subgraph contains a large $K_{2,p}$-subgraph in the deeper layers of some surface-embedded subgraph of $G_{\uparrow}$.  We eliminate these problematic $K_{2,p}$-surface-embedded subgraphs in two steps:
\begin{inparaenum}[(i)]
    \item We do further contractions of connected surface-embedded subgraphs of $G_{\uparrow}$ that are separated from the rest of $G_{\uparrow}$ by short surface-separating cycles and major apex vertices.
    \item We remove certain redundant surface degree-$2$ vertices.
\end{inparaenum}

The result of this process is a minor $G_{\uparrow\bullet}$ of $G_{\uparrow}$.  We augment $G_{\uparrow\bullet}$ with a set of `special' vertices, each of which dominates a connected component of the subgraph of $G_{\uparrow\bullet}$ induced by the sixth layer of $G_{\uparrow\bullet}$.  In particular, each vertex $v\in V(G_{\uparrow\bullet})\setminus V(G)$ is adjacent to exactly one such vertex $\alpha$.  The addition of these special vertices only increases the treewidth by a small amount. On this augmented graph, we apply an enhanced version of the divide-and-conquer strategy for $K_{s,t}$-subgraph-free bounded treewidth graphs of \citet{LW2,LW3} to find an $(s+1)$-colouring with bounded clustering where, for each special vertex $\alpha$, the subgraph dominated by $\alpha$ completely avoids the colour used by $\alpha$.

We now lower the curtain that was previously raised, and assign each vertex $w$ of $G$ the colour given to the vertex of $G_{\uparrow\bullet}$ that $w$ was contracted into.  The resulting colouring of $G$ does not have bounded clustering. We deal with this in two stages.

In the first stage, we use the island method of \citet{EO16}.  A \defn{$(d,c)$-island} $I$ in a graph $G$ is a set of vertices of size at most $c$ such that each vertex in $I$ has at most $d$ neighbours not in $I$.
Recall that each vertex $v$ of $G_\uparrow$ is obtained by contracting a `skinny' subgraph $X_v$ of $G$; $X_v$ contains only a constant number of vertices from each layer of $G$.  However, each vertex $v$ of $G_{\uparrow\bullet}$ is obtained by contracting an arbitrarily large connected subgraph $X'_v$ of $G_{\uparrow}$, so the corresponding subgraph $X_v$ of $G$ is not necessarily skinny.  However, the process used to obtain $G_{\uparrow\bullet}$ from $G_\uparrow$ ensures that any such vertex $v\in V(G_{\uparrow\bullet})\setminus V(G_\uparrow)$ has at most four neighbours in the embedded part of $G_{\uparrow\bullet}$ that contains $v$ and each of these is a vertex in $V(G_{\uparrow})$. From this, it follows that the neighbourhood of $X_v$ in $G$ is skinny.  As a consequence (see \cref{skinny_components}), $X_v$ is either skinny or $X_v$ contains an $(s,c)$-island of $G$.  By exhaustively removing a sequence $I_1,\ldots,I_m$ of such $(s,c)$-islands (letting $I:=\bigcup_{i=1}^m I_i$), we arrive at a situation in which $X_v-I$ is skinny for each $v\in V(G_{\uparrow\bullet})$.

At this point, we have a colouring of $G-I$ in which each monochromatic component is skinny; that is, each monochromatic component contains a bounded number of vertices from each layer.  This is still not enough because a monochromatic component may be `long'; it may intersect an unbounded number of layers.  However, for each component $C$ of the graph obtained by removing the first five layers of $G$, there is a special vertex $\alpha$ such that the colour of $\alpha$ is not used by any vertex in $C$.  This reserved colour of $\alpha$ is then used (see \cref{component_breaking}) to break any long monochromatic components of $G-I$ that intersect $C$ into `short' pieces, so that each monochromatic component intersects only a bounded number of layers. Thus, each monochromatic component of $G-I$ is short and skinny and therefore has bounded size.

The final step in colouring the curtain $G$ is to colour the $(s,c)$-islands $I_m,\ldots,I_1$ in this order.  Since $s+1$ colours are available and $I_i$ is an $(s,c)$-island in $G-(\bigcup_{j=1}^{i-1} I_j)$, the vertices in $I_i$ can be assigned colours without increasing the size of any monochromatic component in $G-\bigcup_{j=1}^{i} I_j$ and without creating any new monochromatic component of size greater than $|I_i|\le c$.

In this way, we obtain an $(s+1)$-colouring of the curtain $G$ with bounded clustering. Then we use the top-down colouring strategy mentioned earlier to extend this colouring to an $(s+1)$-colouring of the tree of curtains, which is the original $\JJ_{s,t}$-minor-free graph. This completes the high-level description of the proof of \cref{Jst}.

The full proof of \cref{Jst} is completed in \cref{JstMinorFree}. In \cref{apex1proof3}, we complete the proof of \cref{apex1} by explaining the changes needed for the case $s=3$.  \cref{algorithms} concludes the paper by outlining polynomial-time algorithms for computing the colourings in \cref{Kh,Kst,Jst,apex1}.

\Cref{Kh,Kst} give the first optimal bounds on the clustered chromatic number of $K_{h}$-minor-free and $K_{s,t}$-minor-free graphs.  At least as important as these results are the definitions and tools that we develop, including curtain decompositions, raised curtains, and skinniness-preserving contractions to control $K_{2,t}$-subgraphs.  Curtain decompositions have already found applications to other decomposition and colouring problems \citep{DHHJLMMRW25}, and we expect that they will soon find more.  Since the extremal examples for many problems include large $K_{2,t}$-subgraphs, we expect the tools for controlling $K_{2,t}$-subgraphs will also find further applications in different contexts. Indeed, the Planar Graph Product Structure Theorem mentioned above was initially developed to bound the queue-number of planar graphs~\citep{DJMMUW20}, but has since been used to resolve a number of longstanding open problems on planar graphs and other graph classes~\citep{DEJWW20,DEGJMM21,BDJM,BKW26,DFMS21,DEMWW22,EJM23}.

\subsection{Notation}

We use the following notation for a graph $G$. For $v\in V(G)$, let $\mathdefn{N_G(v)}:=\{w\in V(G) \colon vw\in E(G)\}$ and $\mathdefn{N_G[v]}:=N_G(v)\cup\{v\}$. For $S\subseteq V(G)$, let
$\mathdefn{N_G[S]}:=\bigcup_{v\in S}N_G[v]$ and $\mathdefn{N_G(S)}:=N_G[S]\setminus S$.

For two vertices $v$ and $w$ in the same component of a graph $G$, let \defn{$\dist_G(v,w)$} be the number of edges in a shortest path from $v$ to $w$ in $G$. If $v$ and $w$ are in distinct components of $G$, then $\DefNoIndex{\dist_G(v,w)}:=\infty$. For a subset $S\subseteq V(G)$, we write $\mathdefn{\dist_G(v,S)}:=\min\{\dist_G(v,w) \colon w\in S\}$. For $d\in\mathbb{N}$, let $\mathdefn{N^d_G[S]}:=\{v\in V(G) \colon \dist_G(v,S)\leq d\}$.

\section{\boldmath Structure of \texorpdfstring{$k$}{k}-Apex-Minor-Free Graphs}
\label{Structure}

\subsection{Graph Minor Structure Theorem}

The Graph Minor Structure Theorem of \citet{RS-XVI} describes the structure of graphs excluding a fixed minor using four ingredients: tree-decompositions, graphs on surfaces, vortices, and apex vertices. To describe this formally, we need the following definitions.

A \defn{tree-decomposition} of a graph $G$ is a collection $(B_x :x\in V(T))$ of subsets of $V(G)$ (called \defn{bags}) indexed by the vertices of a tree $T$, such that (a) for every edge $uv\in E(G)$, some bag $B_x$ contains both $u$ and $v$, and (b) for every vertex $v\in V(G)$, the set $\{x\in V(T):v\in B_x\}$ induces a non-empty (connected) subtree of $T$. For each edge $xy\in E(T)$ the set $B_x\cap B_y$ is called an \defn{adhesion set}. The \defn{adhesion} of $(B_x:x\in V(T))$ is $\max\{|B_x\cap B_y| \colon xy\in E(T)\}$. The \defn{width} of $(B_x:x\in V(T))$ is $\max\{|B_x| \colon x\in V(T)\}-1$. A \defn{path-decomposition} is a tree-decomposition in which the underlying tree is a path, simply denoted by the corresponding sequence of bags $(B_1,\dots,B_n)$. The \defn{treewidth} of a graph $G$, denoted by \defn{$\tw(G)$}, is the minimum width of a tree-decomposition of $G$.

If $(B_x:x\in V(T))$ is a tree-decomposition of a graph $G$, then the \defn{torso} \defn{$\torso{G}{B_x}$} of a bag $B_x$ is the graph obtained from the induced subgraph $G[B_x]$ by adding edges so that $B_x\cap B_y$ is a clique for each edge $xy\in E(T)$.

A \defn{rooted tree} consists of a tree $T$ and a distinguished vertex of $T$ called the \defn{root}. The \defn{depth} of a vertex $v$ in a tree $T$ rooted at $r\in V(T)$ is $\dist_T(r,v)$. A tree-decomposition $(B_x:x\in V(T))$ of a graph $G$ is \defn{rooted} if $T$ is rooted. If $T$ is rooted at $r\in V(T)$, then $\torso{G}{B_r}$ is called the \defn{root torso}.

A \defn{surface} is a compact connected 2-dimensional manifold without boundary. For any terms related to graphs on surfaces that are not defined here, the reader is referred to the monograph by \citet{MoharThom}.

Let $G_0$ be a graph embedded in a surface $\Sigma$. A closed topological disc $D$ in $\Sigma$ is \hdefn{$G_0$}{clean} if the interior of $D$ is disjoint from $G_0$, and the boundary of $D$ only intersects $G_0$ in vertices of $V(G_0)$. Let $x_1,\dots,x_n$ be the vertices of $G_0$ on the boundary of $D$ in the order they appear on the boundary of $D$. A \hdefn{$D$}{vortex} consists of a graph $H$ and a path-decomposition $(B_1,\dots,B_n)$ of $H$ such that $x_i\in B_i$ for each $i\in\{1,2,\dots,n\}$, and $V(G_0\cap H)=\{x_1,\dots,x_n\}$. The vertices $x_1,\ldots,x_n$ are called the \defn{boundary vertices} of the $D$-vortex.

\begin{figure}
  \begin{center}
    \begin{tabular}{c@{\hspace{1cm}}c}
      \includegraphics[page=5]{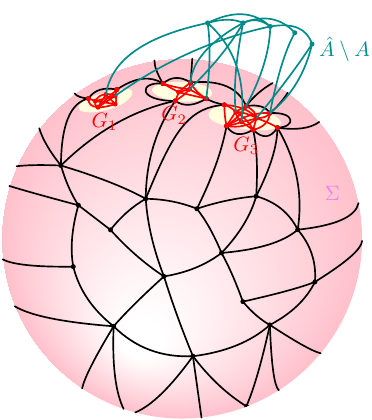} &
      \includegraphics[page=6]{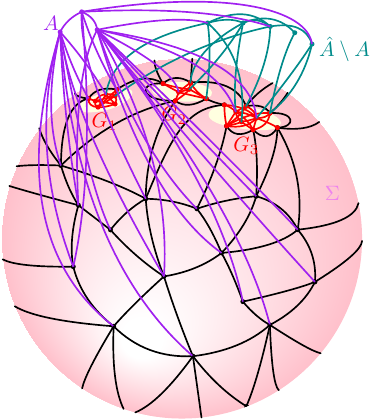}
    \end{tabular}
  \end{center}
  \definecolor{darkcyan}{rgb}{0 0.545 0.545}
  \definecolor{purple}{rgb}{0.627 0.125 0.941}
  \caption{A $(0,\ell)$-almost-embedded graph and a $(3,\ell)$-almost-embedded graph.  The $G_0$-clean topological discs used to define vortices are yellow.  Edges and vertices that participate in vortices are \textcolor{red}{red}.  Non-major apex vertices and their incident edges are \textcolor{darkcyan}{dark cyan}. Major apex vertices and their incident edges are \textcolor{purple}{purple}.  Non-top edges and vertices are black.}
  \label{almost_embedding}
\end{figure}

As illustrated in \cref{almost_embedding}, for integers $a,\hat{a},g,r\geq 0$ and $w\geq1$,
an \hdefn{$(a,\hat{a},g,r,w)$}{almost-embedding} of a graph $G$ is a tuple $\mathcal{E}:=(A,\hat{A},G_0,G_1,\dots,G_r)$ such that:
\begin{itemize}
\item  $A\subseteq \hat{A}\subseteq V(G)$ with $|A|\le a$ and $|\hat{A}|\leq \hat{a}$;
\item $G_0,G_1,\dots,G_r$ are subgraphs of $G$ such that $G-\hat{A} = G_{0} \cup G_{1} \cup \cdots \cup G_r$;
\item $N_G(v)\subseteq \hat{A}\cup V(G_1)\cup\cdots\cup V(G_r)$ for each $v\in \hat{A}\setminus A$;
\item $G_{1}, \dots, G_r$ are pairwise vertex-disjoint;
\item $G_{0}$ is embedded in a surface $\Sigma$ of Euler genus at most $g$;
\item there are pairwise disjoint  $G_0$-clean  topological discs $D_1,\dots,D_r$ in $\Sigma$; and
\item $G_i$ is a $D_i$-vortex with a path-decomposition $(B_1,\dots,B_{n_i})$ of width at most $w$, for each $i\in\{1,2,\dots,r\}$.
\end{itemize}

The \defn{top} of $G$ (with respect to $\mathcal{E}$) is $\hat{A}\cup V(G_1)\cup \cdots\cup V(G_r)$. The \defn{near-top} of $G$ is the union of the top of $G$ and $N_{G_0}[\bigcup_{i=1}^r V(G_i)\cap V(G_0)]$.

The graph $G_0$ is called the \defn{embedded part} of $\mathcal{E}$. The vertices in $\hat{A}$ are called \defn{apex vertices}; those in $A$ are called \defn{major} apex vertices and those in $\hat{A}\setminus A$ are called \defn{non-major} apex vertices.  Major apex vertices can be adjacent to any vertex of $G$ but the neighbourhood of non-major apex vertices is restricted to the top of $G$.  The precise number, $a$, of major apex vertices in an almost-embedding is critical in this work, so we say that a graph is \hdefn{$(a,\ell)$}{almost-embeddable} if it has an  $(a,\hat{a},g,r,w)$-almost-embedding where each of $\hat{a},g,r,w$ is at most $\ell$. A graph $G$ equipped with an  $(a,\ell)$-almost-embedding of $G$ is said to be \hdefn{$(a,\ell)$}{almost-embedded}.

We need the following result, which is probably well known to experts in the field, but for which we have not been able to find any reference. A closed curve in a surface $\Sigma$ is \defn{contractible} if it bounds a topological disc of $\Sigma$.

\begin{lem}
\label{lem:planardisc}
Let $G$ be a 2-connected graph embedded in a surface $\Sigma$ so that every cycle in $G$ is contractible. Then there is a face $f$ of the embedding such that $\Sigma-f$ is a topological disc of $\Sigma$ bounded by a cycle of $G$. In particular, $G$ is planar.
 \end{lem}

\begin{proof}
  Since $G$ is 2-connected, it has an open ear-decomposition \cite{Whitney32a}; that is, there is a sequence of subgraphs $G_1,\ldots,G_t=G$ of $G$ such that $G_1$ is a cycle and for each $i\ge 1$, $G_{i+1}$ is obtained from $G_{i}$ by adding a path $P_{i+1}$ that has distinct endpoints $u_i,v_i\in V(G_{i})$ and is internally disjoint from $V(G_{i})$.

  We now prove by induction that for each $i\ge 1$, there is a face $f_i$ of the induced embedding of $G_i$ in $\Sigma$ such that $\Sigma-f_i$ is a topological disc of $\Sigma$ bounded by a cycle $C_i$ of $G_i$. This is clear for $i=1$, by taking $C_1$ to be the unique cycle of $G_1$, so we can assume that $i\ge 2$. Assume that the result holds for $G_{i-1}$. If $P_i\subseteq \Sigma-f_{i-1}$, then the result also holds for $G_i$ by setting $f_i=f_{i-1}$ and $C_i=C_{i-1}$. Since $G_i$ is embedded in $\Sigma$, we may assume that $P_i-\{u_{i-1},v_{i-1}\}$ lies in $f_{i-1}$, and thus $u_{i-1},v_{i-1}$ both lie on the boundary cycle $C_{i-1}$. Let $H$ be the subgraph of $G$ consisting of the union of $C_{i-1}$ and $P_i$. Since the three cycles of $H$ are contractible, it follows from the proof of \cite[Proposition 4.3.1]{MoharThom} that $P_i$ divides $f_{i-1}$ into two faces, one of which (call it $f_i$) has the property that $\Sigma-f_i$ is a disc that contains $G_i$, and which is bounded by one of the two cycles of $H$ distinct from $C_{i-1}$.
\end{proof}

A graph $X$ is \hdefn{$k$}{apex} if $X-A$ is planar for some $A\subseteq V(X)$ of size at most $k$. We use the following version of the Graph Minor Structure Theorem, due to \citet{DvoTho}.  (A version of \cref{DvoThoOriginal} also appears in \citet[Theorem~16]{DK23}.)

\begin{thm}[{\protect\citep[Theorem~12]{DvoTho}}]
\label{DvoThoOriginal}
    For every integer $k\ge 1$ and every $k$-apex graph $X$ there exists an integer $\ell_0$ such that every $X$-minor-free graph $G$ has a  tree-decomposition $(B_x\colon x\in V(T))$ such that for each node $x\in V(T)$:
    \begin{enumerate}[\rm(i)]
        \item the torso $\torso{G}{B_x}$ is a $(k-1,\ell_0)$-almost-embedded graph; and
        \item every 3-cycle in the embedded part of $\torso{G}{B_x}$ is the boundary of a 2-cell face (that is, a face that is an open topological disc).
    \end{enumerate}
\end{thm}

For graphs $H$ and $G$, an \hdefn{$H$}{model} in $G$ is a set $\{G_x:x\in V(H)\}$ of pairwise vertex-disjoint connected subgraphs of $G$, such that for each $xy\in E(H)$ there is an edge of $G$ between $G_x$ and $G_y$. Each subgraph $G_x$ is called a \defn{branch set} of the model.
Observe that $H$ is a minor of $G$ if and only if there is an $H$-model in $G$. An $H$-model $\{G_x:x\in V(H)\}$ in $G$ is \defn{faithful} if $x\in V(G_x)$, for  each $x\in V(H)$. That is, in a faithful model, each branch set is indexed by some vertex of $G$ that is in the branch set it indexes.
If there is a faithful $H$-model in $G$, then $H$ is a \defn{faithful} minor of $G$. This terminology is due to \citet{Grohe17}.

Let $\TT:=(B_x:x\in V(T))$ be a tree-decomposition of a graph $G$ in which each torso $\torso{G}{B_x}$ is $(a,\ell)$-almost-embedded.  For each $x\in V(T)$, the \defn{lower torso} \defn{$\ltorso{G}{B_x}$} is the supergraph of $G[B_x]$ obtained by adding each edge $uv\in E(\torso{G}{B_x})\setminus E(G[B_x])$ if neither $u$ nor $v$ is in the top of $\torso{G}{B_x}$.  Equivalently, $\ltorso{G}{B_x}$ is the subgraph of $\torso{G}{B_x}$ in which an edge $uv\in E(\torso{G}{B_x})$ is removed if and only if $uv\not\in E(G[B_x])$ and at least one of $u$ or $v$ is in the top of $\torso{G}{B_x}$.
We say that $\TT$ is \defn{lower-minor-closed} if, for each $x\in V(T)$, $G$ contains a faithful $\ltorso{G}{B_x}$-model $\{G_v:v\in B_x\}$ such that, for each vertex $v$ in the top of $\torso{G}{B_x}$, $G_v$ consists only of the vertex $v$.  Note that this implies that $\ltorso{G}{B_x}$ is a minor of $G$.

The definitions of torso and lower-minor-closed are motivated by the fact that we will eventually contract subgraphs of $G[B_x]$ that are connected in $\torso{G}{B_x}$ (but not necessarily connected in $G[B_x]$) in order to obtain a bounded treewidth graph $G[B_x]_{\uparrow}$.  Each connected subgraph of $\torso{G}{B_x}$ that we contract will avoid all vertices in the top of $\torso{G}{B_x}$, so it is also a connected subgraph of $\ltorso{G}{B_x}$.  This ensures that, even though we contracted sets of vertices that are not connected in $G[B_x]$, the contracted graph $G[B_x]_{\uparrow}$ is a minor of $\ltorso{G}{B_x}$, and thus is also a minor of $G$, assuming that $\TT$ is lower-minor-closed.  This ensures that $G[B_x]_{\uparrow}$ is $X$-minor-free if $G$ is $X$-minor-free.

The tree-decomposition of \cref{DvoThoOriginal} can be modified to obtain the decomposition described in \cref{DvoThoCorollaryWithMinors} below. The proof of \cref{DvoThoCorollaryWithMinors}, which involves some careful restructuring of the $(k-1,\ell_0)$-almost-embedded torsos in the tree-decomposition of \cref{DvoThoOriginal}, is presented in \cref{minors_in_proof}.  It is also possible to prove \cref{Jst,apex1} with a version of this lemma that does not guarantee Property~(3) (see \cref{DvoThoCorollary} in \cref{minors_in_proof}), by working alternately with torsos $\torso{G}{B_x}$ and induced subgraphs $G[B_x]$.  We prefer to work with \cref{DvoThoCorollaryWithMinors} because it makes our proof conceptually simpler and is useful, even critical, in other applications~\cite{DHHJLMMRW25}.

\begin{restatable}{lem}{DvoThoCorollaryWithMinors}
\label{DvoThoCorollaryWithMinors}
  For every integer $k\ge 1$ and every $k$-apex graph $X$ there exists an integer $\ell$ such that every $X$-minor-free graph $G$ has a rooted tree-decomposition $\TT:=(B_x\colon x\in V(T))$ such that:
  \begin{enumerate}[\rm(1)]
    \item for each $x\in V(T)$, the torso $\torso{G}{B_x}$ is a $(k-1,\ell)$-almost-embedded graph;
    \item for each edge $xy$ of $T$ where $y$ is the parent of $x$;
    \begin{enumerate}[\rm(a)]
      \item $B_x\cap B_y$ is contained in the top of $\torso{G}{B_x}$,
      \item $B_x\cap B_y$ is contained in the near-top of $\torso{G}{B_y}$ or
   $|B_x\cap B_y|\le k+2$, and
      \item $B_x\cap B_y$ contains at most three vertices not in the top of $\torso{G}{B_y}$; and
    \end{enumerate}
    \item $\TT$ is lower-minor-closed.
  \end{enumerate}
\end{restatable}

Point (3) implies the existence of a faithful $\ltorso{G}{B_x}$-model $\mathcal{M}_x:=\{G_v:v\in B_x\}$ in $G$ for each $x\in V(T)$.   We may assume that, for each $x\in V(T)$ and each $v\in B_x$, $G_v$ is vertex-minimal and edge-minimal.
Under this assumption, for each such $\ltorso{G}{B_x}$-model and each $v\in B_x$, we have $V(G_v)\cap B_x= \{v\}$ and $V(G_v)\cap B_z\subseteq \{v\}$ for any node $z$ that is not a descendant of $x$.  The first of these facts follows from the faithfulness of $\mathcal{M}_x$, so that $v'$ is a vertex of $G_{v'}$ and therefore not a vertex of $G_{v}$ for any $v'\in B_x\setminus\{v\}$. The second follows from the fact that, if $y$ is the parent of $x$, then each vertex in $B_x\cap B_y$ is in the top of $\torso{G}{B_x}$.  If $v\in B_x\cap B_y$, then $G_v=(\{v\},\emptyset)$. Otherwise, the vertices in $B_x\cap B_y$ (none of which are in $G_v$) separate $v$ from every vertex in $B_z\setminus (B_x\cap B_y)$. Thus, any vertex of $G_v$ in $B_z\setminus B_x\cap B_y$ could be removed from $G_v$, which is not possible by the minimality assumption.

\subsection{Pre-Curtains}\label{sec:precurtain}

\begin{figure}[htbp]
  \begin{center}
    \includegraphics{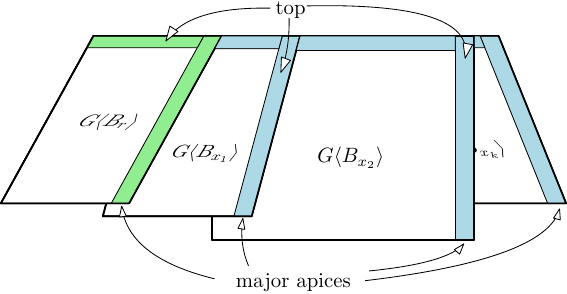}
  \end{center}
  \caption{A $(k,\ell)$-pre-curtain}
  \label{pre_curtain}
\end{figure}

As illustrated in \cref{pre_curtain}, a graph $G$ is a \hdefn{$(k,\ell)$}{pre-curtain} if it has a rooted tree-decomposition $\mathcal{T}:=(B_x:x\in V(T))$ in which each torso $\torso{G}{B_x}$ is a $(k,\ell)$-almost-embedded graph, and for each edge $xy$ of $T$ where $y$ is the parent of $x$, $B_x\cap B_y$ is contained in the top of $\torso{G}{B_x}$ and in the near-top of $\torso{G}{B_y}$.  We say that $G$ is a $(k,\ell)$-pre-curtain \defn{described by} $\TT$.  The \defn{top} of $G$ is the union of the tops of the torsos of $\TT$ and the \defn{near-top} of $G$ is the union of the near-tops of the torsos of $\TT$. Note that each vertex of $G$ that is not in the top of $G$ appears in exactly one bag of $\TT$.  Recall that the \defn{root torso} of $G$ is the torso associated with the root bag of $\TT$.

A graph $J$ is a \defn{tree of $(k,\ell)$-pre-curtains} if it has a rooted tree-decomposition $\TT:=(B_x:x\in V(T))$ such that:
\begin{compactitem}
  \item for each $x\in V(T)$, $\torso{J}{B_x}$ is a $(k,\ell)$-pre-curtain, described by some rooted tree-decomposition $\TT_x$,
  \item for each edge $xy$ of $T$ where $y$ is the parent of $x$, $B_x\cap B_y$ has size at most $k+3$ and is contained in the top of the root torso of $\torso{J}{B_x}$ in $\TT_x$.
\end{compactitem}
Again, we say that $J$ is a tree of $(k,\ell)$-pre-curtains \defn{described by} $\TT$.

Let $J$ be a tree of pre-curtains described by the tree-decomposition $\TT:=(B_x:x\in V(T))$. For each $x\in V(T)$, the \defn{top} of the pre-curtain $\torso{J}{B_x}$ is the union of the tops of the (almost-embedded) torsos of $\TT_x$.  For each $x\in V(T)$, the \defn{lower torso $\ltorso{J}{B_x}$} is the supergraph of $J[B_x]$ obtained by adding each edge $uv\in E(\torso{J}{B_x})\setminus E(J[B_x])$ if neither $u$ nor $v$ is in the top of $\torso{J}{B_x}$.  The tree of pre-curtains $J$ is \defn{lower-minor-closed} if, for each $x\in V(T)$, $J$ contains a faithful $\ltorso{J}{B_x}$-model $\{J_v:v\in B_x\}$ such that, for each vertex $v$ in the top of $\torso{J}{B_x}$, $J_v$ consists only of the vertex $v$.  The wording of these definitions is deliberately identical to that of lower torso and lower-minor-closed in a tree-decomposition whose torsos are almost-embedded graphs; the difference is that each torso $\torso{J}{B_x}$ in a tree of pre-curtains is a pre-curtain described by a tree-decomposition $\TT_x$.

\begin{lem}\label{tree_of_raw_curtains}
    For every integer $k\ge 1$ and every $k$-apex graph $X$ there exists an integer $\ell$ such that every $X$-minor-free graph $J$ is a lower-minor-closed tree of $(k-1,\ell)$-pre-curtains.
\end{lem}

\begin{proof}
  Let $\TT_0:=(B_x:x\in V(T_0))$ be the rooted tree-decomposition of $J$ guaranteed by \cref{DvoThoCorollaryWithMinors}.
  For an edge $xy$ of $T_0$ where $y$ is the parent of $x$, we say that $xy$ is \hdefn{$k$}{heavy} if $B_x\cap B_y$ is contained in the near-top of $\torso{J}{B_y}$, and call it \hdefn{$k$}{light} otherwise. Observe that, by property (2b) of \cref{DvoThoCorollaryWithMinors}, $|B_x\cap B_y|\le k+2$ for each $k$-light edge $xy$ of $T_0$.

  By removing all $k$-light edges of $T_0$ we obtain a forest with component trees $T_{1},\ldots,T_{p}$ rooted at $r_1,\ldots,r_p$, respectively.  For each $i\in\{1,\ldots,p\}$, let $B_i:=\bigcup_{x\in V(T_i)}B_x$. Let $T$ be the tree with vertex set $V(T):=\{1,\ldots,p\}$ in which $ij\in E(T)$ if and only if some $k$-light edge of $T_0$ has an endpoint in $V(T_i)$ and an endpoint in $V(T_j)$.  Root $T$ at the index $r\in V(T)$ such that $T_r$ contains the root of $T_0$.

  Observe that $\TT:=(B_i:i\in V(T))$ is a rooted tree-decomposition of $J$ and, for each $i\in V(T)$, $\TT_i:=(B_x:x\in V(T_i))$ is a rooted tree-decomposition of $J_i:=\torso{J}{B_i}$.  We claim that $J$ is a lower-minor-closed tree of $(k-1,\ell)$-pre-curtains described by $\TT$, where the curtain $\torso{J}{B_i}$ is described by $\TT_i$.

  First we establish that $\torso{J}{B_i}$ is a $(k-1,\ell)$-pre-curtain described by $\TT_i$, for each $i\in V(T)$.  Let $i\in V(T)$ and let $xy$ be an edge of $T_i$ where $y$ is the parent of $x$.  Since $xy$ is $k$-heavy, $B_x\cap B_y$ is contained in the near-top of $\torso{J}{B_y}$.  By property (2a) of \cref{DvoThoCorollaryWithMinors}, $B_x\cap B_y$ is contained in the top of $\torso{J}{B_x}$.  Therefore, for each $i\in\{1,\ldots,p\}$, $\torso{J}{B_i}$ is a $(k-1,\ell)$-pre-curtain described by $\TT_i$.

  Next we establish that $J$ is a tree of $(k-1,\ell)$-pre-curtains described by $\TT$.  Each edge $ij$ of $T$ with $j$ the parent of $i$ corresponds to a $k$-light edge $xy$ of $T_0$ with $y\in V(T_j)$ and $x\in V(T_i)$.  Since $T$ is rooted at $r$, $y$ is the parent of $x$ in $T_0$.  By property~(2a) of \cref{DvoThoCorollaryWithMinors}, $B_x\cap B_y$ is contained in the top of $\torso{J}{B_x}$, which is the root torso of the pre-curtain $\torso{J}{B_i}$. Thus $\TT$ satisfies the first requirement for a tree of $(k-1,\ell)$-pre-curtains.  Since $xy$ is $k$-light, $|B_i\cap B_j|=|B_x\cap B_y|\le k+2$.  Thus $\TT$ satisfies the second requirement for a tree of $(k-1,\ell)$-pre-curtains.

  It remains to show that the tree of $(k-1,\ell)$-pre-curtains $J$ described by $\TT$ is lower-minor-closed.  That is, we must show that for each $i\in\{1,\ldots,p\}$, there is a faithful $\ltorso{J}{B_i}$-model $\mathcal{M}_i:=\{J_v:v\in B_i\}$ in $J$ in which, for each vertex $v$ in the top of $\torso{J}{B_i}$, $J_v$ is the $1$-vertex graph that contains $v$. The model $\mathcal{M}_i$ is defined as follows.  By property (3) of \cref{DvoThoCorollaryWithMinors}, for each $x\in V(T_0)$, there exists a faithful $\ltorso{J}{B_x}$-model $\mathcal{M}_x:=\{J^x_v:v\in B_x\}$ in $J$ in which, for each vertex $v$ in the top of $\torso{J}{B_x}$, $J^x_v$ is the $1$-vertex graph that contains $v$.  For each $v\in B_i$, if $v$ is in the top of $\torso{J}{B_i}$, then define $J_v$ to be the $1$-vertex graph that contains $v$.  If $v$ is not in the top of $\torso{J}{B_i}$ then there is a unique node $x$ of $T_i$ such that $v\in B_x$.  Then define $J_v$ to be the component of $J_v^x-(B_i\setminus\{v\})$ that contains $v$.
  Now we verify that $\mathcal{M}_i$ is a faithful $\ltorso{J}{B_i}$-model in $J$.  Obviously $\mathcal{M}_i$ is faithful because $\mathcal{M}_x$ is faithful for each $x\in V(T_i)$.

  Next we verify that the branch-sets of $\mathcal{M}_i$ are pairwise vertex-disjoint.  Let $v$ and $v'$ be distinct vertices in $B_i$, and let $x$ and $x'$ be the minimum-depth nodes of $T_i$ such that $v\in B_x$ and $v'\in B_{x'}$. Then $J_v\subseteq J^{x}_v- (B_i\setminus\{v\})$ and $J_{v'}\subseteq J^{x'}_{v'}-(B_i\setminus\{v'\})$.  If $x=x'$ then $J_v$ and $J_{v'}$ are vertex-disjoint since they are each subgraphs of the branch sets $J^x_v$ and $J^x_{v'}$ in the $\ltorso{J}{B_x}$-model $\mathcal{M}_x$.  If $x\neq x'$ then we may assume without loss of generality that $x$ is not an ancestor of $x'$.  Let $y$ be the parent of $x$. Since $T_i$ is a (connected) tree that contains $x$ and $x'$, $y$ is in $T_i$. By the definition of $x$, $v\not\in B_y$, so $v\not\in B_x\cap B_y$.  Therefore the component of $J-(B_x\cap B_y)$ that contains $v$ does not contain $v'$.  If $v'\not\in B_x\cap B_y$ then this implies that $J_{v'}$ and $J_v$ are vertex-disjoint since they are contained in different components of $J-(B_x\cap B_y)$.  If $v'\in B_x\cap B_y$ then $v'$ is in the top of (the almost-embedded graph) $\torso{J}{B_x}$, so $v'$ is in the top of (the pre-curtain) $\torso{J}{B_i}$, so $J_{v'}$, which only contains the vertex $v'$, is also vertex-disjoint from $J_v$.

  Finally, we verify that for each edge $vw\in E(\ltorso{J}{B_i})$, $J_v$ and $J_w$ are adjacent in $J$.  If $vw\in E(J[B_i])$, then this follows immediately from the definitions of the branch sets and the faithfulness of $\mathcal{M}_x$ for each $x\in V(T_i)$.  Otherwise, $vw\in E(\ltorso{J}{B_i})\setminus E(J[B_i])$, which implies that $T$ contains an edge $ij$ with $B_i\cap B_j\supseteq\{v,w\}$.  The existence of $ij\in E(T)$ implies that there exists a $k$-light edge  $xy\in E(T_0)$ where $x\in V(T_i)$, $y\in V(T_j)$, and $\{v,w\}\subseteq B_x\cap B_y$.  Since $vw\in E(\ltorso{J}{B_i})\setminus E(J[B_i])$, neither $v$ nor $w$ are in the top of the pre-curtain $\torso{J}{B_i}$.  Therefore $y$ cannot be the parent of $x$ in $T_0$ (as that would imply $\{v,w\} \subseteq B_x\cap B_y \subseteq \ttop{J}{B_x} \subseteq \ttop{J}{B_i}$).  Thus $x$ is the parent of $y$ in $T_0$.  Thus, $v$ and $w$ are both in $B_x$ and neither is in the top of $\torso{J}{B_x}$.
  Therefore, $B_x$ is the only bag of $\TT_i$ that contains $v$ and $B_x$ is the only bag of $\TT_i$ that contains $w$. Therefore $J_v=J^x_v -(B_i\setminus\{v\})=J^x_v$ and $J_w=J^x_w -(B_i\setminus\{w\})=J^x_w$.  Therefore $J_v$ and $J_w$ are adjacent in $J$ since $vw\in E(\ltorso{J}{B_x})$ and $\mathcal{M}_x$ is a faithful $\ltorso{J}{B_x}$-model in $J$.
\end{proof}

The notions of lower torso and lower-minor-closed will not be immediately needed in this section, but they will be crucial later when specializing pre-curtains to curtains.

Say that a colouring of a graph $G$ \defn{properly extends} a precolouring of a set $Y\subseteq V(G)$ if no monochromatic component of $G$ has a vertex in $Y$ and a vertex in $V(G)\setminus Y$; that is, each monochromatic component is a subgraph of $G[Y]$ or of $G-Y$.  In particular, this happens if and only if each vertex in $N_G(Y)$ is coloured differently from all its neighbours in $Y$. \cref{tree_of_raw_curtains} allows the problem of colouring an $X$-minor-free graph to be reduced to the problem of colouring a tree of pre-curtains. We now show that for this task, it is sufficient to be able to properly extend precolourings in pre-curtains.

\begin{lem}\label{Lungs}
  Let $J$ be a graph that is a tree of $(k,\ell)$-pre-curtains described by $\TT:=(B_x:x\in V(T))$ and such that, for any $(k,\ell)$-pre-curtain $\torso{J}{B_x}$, any precolouring with at most $m$ colours of a set $S$ of at most $k+3$ vertices in the top of the root torso of $\torso{J}{B_x}$  can be properly extended to an $m$-colouring of $J[B_x]$ with clustering at most $c$. Then $J$ has an $m$-colouring with clustering at most $c$.
\end{lem}

\begin{proof}
  Let  $(B_x:x\in V(T))$ be a tree-decomposition of $J$ rooted at a node $r\in V(T)$ as in the definition of a tree of $(k,\ell)$-pre-curtains.  Let $x_0,\ldots,x_p$ be the nodes of $T$, ordered so that $x_0$ is the root of $T$ and, for each $i\ge 1$, the parent of $x_i$ is among $x_0,\ldots,x_{i-1}$.  By the definition of a tree of $(k,\ell)$-pre-curtains, $\torso{J}{B_{x_0}}$ is a $(k,\ell)$-pre-curtain. By assumption, $J[B_{x_0}]$, as a spanning subgraph of $\torso{J}{B_{x_0}}$, has an $m$-colouring with clustering at most $c$.  Suppose now that $J[B_{x_0}\cup\cdots\cup B_{x_{i-1}}]$ has an $m$-colouring with clustering at most $c$, and let $x_j$ be the parent of $x_i$.  By the definition of $(k,\ell)$-pre-curtain, $S_i:=B_{x_j}\cap B_{x_i}$ has size at most $k+3$ and is contained in the top of the root torso of the $(k,\ell)$-pre-curtain $\torso{J}{B_{x_i}}$.  Treat $S_i$ as a precoloured set in the top of the root torso of $\torso{J}{B_{x_i}}$, which by assumption can be properly extended to an $m$-colouring of $J[B_{x_i}]$ with clustering at most $c$.  In this way, we properly extend the colouring of $J[B_{x_0}\cup\cdots\cup B_{x_{i-1}}]$ to an $m$-colouring of $J[B_{x_0}\cup\cdots\cup B_{x_{i}}]$.  Since this is a proper extension, any monochromatic component in this colouring is contained in $J[B_{x_0}\cup\cdots\cup B_{x_{i-1}}]$ or is contained in $J[B_{x_i}]$.  In either case, each monochromatic component has size at most $c$.  Doing this for $i=1,\ldots,p$ in order gives an $m$-colouring of $J$ with clustering at most $c$.
\end{proof}

\subsection{A Digression}
\label{Digression}

For expository purposes, we now pause to show how the notion of $(k,\ell)$-pre-curtains can be used, with some existing results, to quickly establish that the clustered chromatic number of $K_{h}$-minor-free graphs is at most $h+4$.  This is weaker than \cref{Kh}, but stronger than all previous results on clustered colouring of $K_h$-minor-free graphs except for those of \citet{LW2,LW3}, which rely on a technical lemma of \citep{LW1} with a 70-page proof. The purpose of this section is to introduce some of the techniques used to prove \cref{apex1}, which are then expanded upon to prove \cref{Jst}. We need the following result, which is a consequence of \citet[Corollary~18]{LW2}. \Cref{lem:LW}, proven in \cref{BoundedTreewidthKstFree}, is an extension of this result.

\begin{lem}\label{lem:LW0}
There is a function $f$ such that for every $K_{s,t}$-subgraph-free graph $G$ of treewidth at most $k$ and every $P\subseteq V(G)$ such that each vertex in $V(G)\setminus P$ has at most $s$ neighbours in $P$, any precolouring of $P$ with $s+1$ colours can be properly extended to an $(s+1)$-colouring of $G$ with clustering at most $f(k,s,t,|P|)$.
\end{lem}

We note that the assumption that each vertex of $V(G)\setminus P$ has at most $s$ neighbours in $P$ is necessary: if a vertex $v\in V(G)\setminus P$ has $s+1$ neighbours in $P$, each coloured with a different colour, then the precolouring cannot be properly extended to $v$.

\begin{prop}\label{h_plus_four}
  There is a function $f$ such that for any integer $h\ge 5$, every $K_h$-minor-free graph has an $(h+4)$-colouring with clustering at most $f(h)$.
\end{prop}

\begin{proof}
  Let $J$ be a $K_{h}$-minor-free graph. Since $K_4$ is planar, $K_h$ is $(h-4)$-apex.  Therefore, by \cref{tree_of_raw_curtains}, for some integer $\ell=\ell(h)$, $J$ is a tree of $(h-5,\ell)$-pre-curtains described by a tree-decomposition $\TT:=(B_x:x\in V(T))$.  Fix some $x\in V(T)$, let  $\TT_x:=(C_y:y\in V(T_x))$ be the tree-decomposition that describes the $(h-5,\ell)$-pre-curtain $\torso{J}{B_x}$, and let $G:=J[B_x]$. Let $S$ be a set of at most $h-2$ precoloured vertices contained in the top of the root torso of $\torso{J}{B_x}$. By \cref{Lungs}, it suffices to show that the precolouring of $S$ can be properly extended to an $(h+4)$-colouring of $G$ with clustering at most $f(h)$. Since $J$ is $K_h$-minor-free and $G$ is a subgraph of $J$, $G$ is also $K_h$-minor-free. Contracting a matching of size $h-2$ in $K_{h-1,h-1}$ produces $K_h$. So $G$ is $K_{h-1,h-1}$-minor-free, implying that $G$ is $K_{h-1,h-1}$-subgraph-free.

  Let $G'$ be the subgraph of $G$ induced by the top of the pre-curtain $\torso{J}{B_x}$.
  The tree-decomposition $\TT_x$ induces a tree-decomposition $\TT'_x:=(C_y':y\in V(T_x))$ of $G'$, where $C'_y:=C_y\cap V(G')$ for each $y\in V(T_x)$. We now show that each torso $\torso{G'}{C'_y}$ in this decomposition has treewidth at most some $k:=k(h)$. For each $y\in V(T_x)$, let $(A_y, \hat{A}_y, G^y_0, G^y_1, \ldots, G^y_r)$ be the $(h-5,\ell)$-almost-embedding of the torso $\torso{\torso{J}{B_x}}{C_y}$. The graph $\torso{G'}{C'_y}-\hat{A}_y$ can be written as the union of a graph $H_0$ embedded on some surface $\Sigma$ of Euler genus at most $\ell$, and vortices $H_1,\ldots,H_{r'}$ as in the definition of an almost-embedding (with $r'\le \ell$ and such that each $H_j$ has a path-decomposition of width at most $\ell$), with the additional property that each vertex of $H_0$ lies on the boundary of a vortex (by the definition of a top). Let $H_0'$ be the supergraph of $H_0$ obtained as follows:  For each vortex $H_j$ associated with the $H_0$-clean topological disc $D_j$, add a vertex $v_j$ in the interior of $D_j$ that is adjacent to each vertex on the boundary of $D_j$. Next, add edges so that each pair of consecutive vertices encountered while traversing the boundary of $D_j$ are adjacent.

  Then $H_0'$ can be embedded on $\Sigma$ and thus has Euler genus at most $\ell$. Moreover, each connected component of $H_0'$ has diameter at most $3r-1\le 3\ell$ (any two vertices on the boundary of a vortex are at distance at most 2 in $H_0'$, and the surface edges of $H_0'$ are connecting vortices). \citet{Eppstein-Algo00} proved that graphs of bounded Euler genus and bounded diameter have bounded treewidth, thus $H_0'$ has treewidth at most $k':=k'(\ell)$. Thus $H_0$, which is a subgraph of $H_0'$, also has treewidth at most $k'$.  \citet[Lemma 10]{DvoTho} proved that in this case $\bigcup_{i=0}^r H_i=\torso{G'}{C'_y}$ has treewidth at most $k'':=\ell(k'+1)-1$, where $\ell$ accounts for the fact that the vortices $H_1,\ldots,H_r$ all have a path-decomposition of width at most $\ell$.  Therefore, $\torso{G'}{C'_y}-\hat{A}_y$ has a tree-decomposition  of width at most $k''$, so $\torso{G'}{C'_y}$ has a tree-decomposition $\TT_y$ of width at most $k''+|\hat{A}_y|\le k''+\ell=:k$.  For each $yz\in E(T_x)$, some bag of $\TT_y$ contains $C'_y\cap C'_z$ and some bag of $\TT_z$ contains $C'_y\cap C'_z$. Joining these two bags by an edge, for each $yz\in E(T_x)$ gives a tree-decomposition of $G'$ having width at most $k$.  Therefore the treewidth of $G'$ is at most $k$.

  Since $S$ is a subset of the top of the root torso of $\torso{J}{B_x}$, $S$ is a subset of the top of $\torso{J}{B_x}$.  Therefore $S\subseteq V(G')$. Since $|S|\le h-2$, every vertex not in $S$ has at most $h-2\le h-1$ neighbours in $S$. By \cref{lem:LW0}, the precolouring of $S$ can be properly extended to an $h$-colouring of $G'$ with clustering $f(k,h-1,h-1,h-2)$.  Without loss of generality, we may assume this colouring of $G'$ uses colours $\{1,\ldots,h\}$.  Now, consider the graph $G-V(G')$.  Each component of $G-V(G')$ has Euler genus at most $\ell$. By \cref{Genusg}, $G-V(G')$ has a $4$-colouring using colours $\{h+1,\ldots,h+4\}$ with clustering $f'(\ell)$, for some function $f'$.  The resulting colouring of $G$ is clearly a proper extension of the precolouring of $S$ and has clustering $\max\{f(k,h-1,h-1,h-2),f'(\ell)\}$.
\end{proof}

We remark that the proof of \cref{h_plus_four} is easily adapted to prove the following result (whose statement is more technical, but which is stronger and immediately implies \cref{h_plus_four}):

\begin{prop}\label{EasykapexKstFree}
  There is a function $f$ such that for every $k$-apex graph $X$, 
  for all integers $t\geq s\geq 1$, every $X$-minor-free $K_{s,t}$-subgraph-free graph is $\max\{k+7,s+5\}$-colourable with clustering at most $f(k,s,t)$.
\end{prop}

\begin{proof}[Proof Sketch]
  Let $J$ be an $X$-minor-free and $K_{s,t}$-subgraph-free graph. By \cref{tree_of_raw_curtains}, $J$ is a tree of $(k-1,\ell)$-pre-curtains. Let $m=\max\{k+7,s+5\}$. Let $G$ be a subgraph of $J$ induced by the vertex set of some $(k-1,\ell)$-pre-curtain in the decomposition of $J$. In order to show the desired result, by \cref{Lungs}, it suffices to show that the precolouring of any set $S$ of at most $(k-1)+3=k+2$ vertices in the top of the root torso of the pre-curtain can be properly extended to an $m$-colouring of $G$. As above, we observe that the subgraph $G'$ of $G$ induced by the top of the pre-curtain has bounded treewidth, and since $G'$ (as a subgraph of $G$ and $J$) is $K_{s,t}$-subgraph-free, it is $K_{m-5,t}$-subgraph-free. Since $|S|\le k+2\le m-5$, it follows from \cref{lem:LW0} that the precolouring of $S$ can be properly extended to an $(m-4)$-colouring of $G'$ with bounded clustering. As above, we combine this with a 4-colouring of $G-V(G')$ with bounded clustering, to obtain an $m$-colouring of $G$ properly extending the precolouring of $S$, as desired.
\end{proof}

\subsection{Layered Partitions and Curtains}

We now need the following definitions from the literature. A \defn{partition} of a graph $G$ is a collection $\PP$ of non-empty sets of vertices in $G$ such that each vertex of $G$ is in exactly one element of $\PP$. Each element of $\PP$ is called a \defn{part}. The \defn{quotient} of $G$ with respect to $\PP$ is the graph, denoted by \defn{$G/\PP$}, with vertex set $\PP$ where distinct parts $A,B\in \PP$ are adjacent in $G/\PP$ if and only if some vertex in $A$ is adjacent in $G$ to some vertex in $B$.  A partition of $G$ is \defn{connected} if the subgraph induced by each part is connected. Observe that in this case, the quotient is exactly the minor of $G$ obtained by contracting each part $P\in \PP$ into a single vertex $P$.  It is often convenient to omit singleton sets when defining partitions and quotient graphs, and only require that $\PP$ is a set of disjoint subsets of $V(G)$. In this case, $G/\PP:=G/\PP'$ where $\PP':=\PP\cup\{\{v\}:v\in V(G)\setminus \bigcup_{P\in \PP}P\}$.

Partitions, products and clustered colouring are intimately related, since it follows from the definitions that the following are equivalent for any graph $G$:
\begin{itemize}
\item $G$ is $k$-colourable with clustering $c$,
\item $G$ is isomorphic to a subgraph of $H\boxtimes K_c$ for some properly $k$-colourable graph $H$,
\item $G$ has a partition $\PP$ such that every part of $\PP$ has at most $c$ vertices and the quotient $G/\PP$ is properly $k$-colourable.
\end{itemize}

A \defn{layering} of a graph $G$ is an ordered partition $\LL=(L_1,L_2,\dots)$ of $V(G)$ such that for every edge $vw\in E(G)$, if $v\in L_i$ and $w\in L_j$, then $|i-j| \leq 1$.   For any integer $i$, we use the shorthands $\mathdefn{L_{\le i}}:=\bigcup_{j\le i} L_j$ and $\mathdefn{L_{\ge i}} :=\bigcup_{j\ge i} L_j$.  We say that $\LL$ is \defn{upward-connected} if for every $i\ge 2$ every vertex in $L_i$ has a neighbour in $L_{i-1}$. A typical example of a layering is a Breadth-First-Search (BFS) layering: set a root vertex $v_j$ in each connected component and for each $i\ge 1$, let  $L_i$ be the set of vertices at distance exactly $i-1$ from one of the vertices $v_j$. Every BFS layering is upward-connected. We use the following useful property of upward-connectivity

\begin{obs}
\label{UpwardConnected}
If $\LL=(L_1,L_2,\dots)$ is an upward-connected layering of a graph $G$, then, for any integer $i\ge 1$, the number of components of $G[L_{\le i}]$ is at most the number of components of $G[L_1]$.
\end{obs}

If $\LL$ is a layering of a graph $G$, then a set $X$ of vertices in $G$ is \hdefn{$\ell$}{skinny} with respect to $\LL$ if $|X\cap L|\le\ell$ for each $L\in\LL$.  We use the following lemma to transform a colouring in which each monochromatic component is $\ell$-skinny to a colouring in which each monochromatic component has bounded size. A monochromatic component $C$ in some colouring $\varphi$ of a graph $G$ is called a \hdefn{$\varphi$}{monochromatic component} in $G$ (this notation is helpful to avoid confusion when we consider several colourings of a given graph, as in the proof of the next result).

\begin{lem}\label{component_breaking}
  Let $G$ be a graph, let $\LL:=(L_1,L_2,\ldots)$ be a layering of $G$, and let $t \ge 1$ be an integer.  Let $\varphi:V(G)\to \{1,\ldots,c\}$ be a $c$-colouring of $G$ such that every $\varphi$-monochromatic component of $G$ is $\ell$-skinny with respect to $\LL$, and for each component $C$ of $G[L_{\ge t}]$, there exists a colour $\phi_C\in \{1,\ldots,c\}\setminus\{\varphi(v):v\in V(C)\}$ that is not used in $C$.  Then $G$ has a $c$-colouring $\varphi'$ with clustering at most $(2c+t-1)\ell$ such that
  \begin{enumerate}[(a)]
    \item\label[paren]{breaker_untouched_top} for each $v\in L_{\le t}$, $\varphi'(v)=\varphi(v)$;
    \item\label[paren]{breaker_untouched_even} for each integer $i\ge 0$, each $v\in L_{t+2i}$ has $\varphi'(v)=\varphi(v)$;
    \item\label[paren]{breaker_component_colours} for each component $C$ of $G[L_{\ge t}]$ and each $v\in V(C)$, $\varphi'(v)\in\{\varphi(v),\phi_C\}$;
    \item\label[paren]{breaker_only_recolour_2vf} for each $v \in L_{\ge t}$, if $\varphi'(v) \neq \varphi(v)$, then $v \in L_i$ for some integer $i \ge t+1$ with $i\equiv t+1+2(\varphi(v)-1)\pmod{2c}$.
  \end{enumerate}
\end{lem}

\begin{proof}
  As illustrated in \cref{breaking_figure}, define a colouring $\varphi'$ of $G$ as follows:  For each component $C$ of $G[L_{\ge t}]$, each colour $\phi\in\{1,\ldots,c\}\setminus\{\phi_C\}$, each integer $i\ge t+1$ with $i\equiv t+1+2(\phi-1) \pmod{2c}$, and each $v\in \{ w\in V(C)\cap L_i:\varphi(w)=\phi\}$, set $\varphi'(v):=\phi_{C}$.  For any $v\in V(G)$ not defined by this rule, set $\varphi'(v):=\varphi(v)$. Clearly, $\varphi'$ is a  $c$-colouring of $G$ and it is straightforward to verify that $\varphi'$ satisfies conditions (a) to (d).  We now show that $\varphi'$ has clustering at most $(2c+t-1)\ell$.

  \begin{figure}
      \centering
      \begin{tabular}{c@{\hspace{2em}}c}
      \includegraphics[page=1]{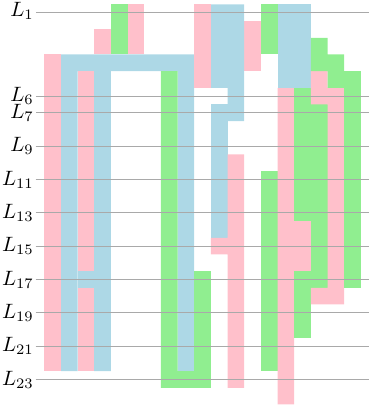} &
      \includegraphics[page=2]{breakers} \\
      $\varphi$ & $\varphi'$
      \end{tabular}
      \caption{The proof of \cref{component_breaking} with $t=6$ and $c=3$}
      \label{breaking_figure}
  \end{figure}

  \emph{Claim.} For each edge $vw$ of $G$, if $\varphi(v)\neq\varphi(w)$ then $\varphi'(v)\neq \varphi'(w)$.

  \begin{clmproof}
    Suppose for the sake of contradiction that $\varphi(v)\neq\varphi(w)$ and $\varphi'(v)=\varphi'(w)$. Without loss of generality, $\varphi(v)\neq\varphi'(v)$. By construction, $v\in V(C) \cap L_i$ where $C$ is a component of $G[L_{\geq t}]$, and $i\equiv t+1+2(\varphi(v)-1)\pmod{2c}$, and $\varphi'(v)=\phi_C$. Say $w\in L_{i'}$. Since $\LL$ is a layering, $i'\geq i-1 \geq t$. Hence $w$ is also in $C$. If $\varphi'(w)= \varphi(w)$, then $\phi_C=\varphi'(v)=\varphi'(w)=\varphi(w)$, which contradicts the assumption that $\phi_C$ is not used in $C$. Thus $\varphi'(w)\neq \varphi(w)$. Since $w\in V(C)$, $\varphi'(w)=\phi_C$ by property~\ref{breaker_component_colours}. By construction, $i'\equiv t+1+2(\varphi(w)-1)\pmod{2c}$. Since $vw$ is an edge and $\LL$ is a layering, $|i-i'|\leq 1$. Furthermore, since $2c$ is even and both $i$ and $i'$ have the same parity as $t+1$ modulo $2c$, $i$ and $i'$ must have the same parity. The only way they can have the same parity and $|i-i'| \le 1$ is if $i = i'$. But if $i=i'$, then $t+1+2(\varphi(v)-1) \equiv t+1+2(\varphi(w)-1) \pmod{2c}$, which implies $2(\varphi(v)-\varphi(w)) \equiv 0 \pmod{2c}$, and thus $\varphi(v)-\varphi(w) \equiv 0 \pmod c$. Since $\varphi(v),\varphi(w)\in\{1,\dots,c\}$, we must have $\varphi(v)=\varphi(w)$. This contradiction shows that $\varphi'(v)\neq\varphi'(w)$.
  \end{clmproof}

  It follows from the claim that each $\varphi'$-monochromatic component $A'$ is contained in some $\varphi$-monochromatic component $A$. Consider the following two cases:
  \begin{compactenum}
    \item $A'$ and $A$ have different colours $\phi'$ and $\phi$, respectively.  That is, $\phi=\varphi(v)\neq\varphi'(v)=\phi'$ for each $v\in V(A')$.   In this case $A'$ is contained in some component $C$ of $G[L_{\ge t}]$, $\phi'=\phi_C$, and $A'=A[L_i]$ for some $i \ge t+1$ with $i\equiv t+1+2(\phi-1)\pmod{2c}$ (since only vertices in layers $L_k$ with $k\equiv t+1 \pmod 2$ change their colour).  Therefore, $|V(A')|= |V(A[L_i])|\le \ell$, since $A$ is $\ell$-skinny with respect to $\LL$.

    \item $A'$ and $A$ have the same colour $\phi$.  In other words, $\varphi'(v)=\varphi(v)=\phi$ for each $v\in V(A')$.  Consider some component $C$ of $G[L_{\ge t}]$.  If $\phi_C=\phi$, then $A'$ does not intersect $C[L_{\ge t}]$ by the definition of $\phi_C$.  Let $C_1,\ldots,C_r$ be the components of $G[L_{\ge t}]$ for which $\phi_{C_i}\neq \phi$. Then, for each $j\in\{1,\ldots,r\}$, $C_j[L_i]$ has no vertex of colour $\phi$ for any $i\ge t+1$ with $i\equiv t+1+2(\phi-1) \pmod{2c}$.  Therefore the vertices of $A'[L_{\ge t+1}]$ are contained in $2c-1$ consecutive layers of $\LL$.  Since $A'$ is contained in $A$,  and $A$ is $\ell$-skinny with respect to $\LL$, this implies that $|V(A'[L_{\ge t+1}])|\le (2c-1)\ell$.  Furthermore, $|V(A'[L_{\le t}])|=|V(A[L_{\le t}])|\le t\ell$.  Therefore, $|V(A')|\le (2c+t-1)\ell$.
  \end{compactenum}
  Therefore, $\varphi'$ is a $c$-colouring of $G$ with clustering at most $(2c+t-1)\ell$.
\end{proof}

A \defn{layered partition} $(\PP, \LL)$ of a graph $G$ consists of a partition $\PP$ and a layering $\LL$ of $G$. Adjectives used to modify $\PP$ or $\LL$ have the obvious meaning when applied to $(\PP,\LL)$.  For example, $(\PP,\LL)$ is \defn{connected} if $\PP$ is connected and $(\PP,\LL)$ is \defn{upward-connected} if $\LL$ is upward-connected.  We say that $(\PP,\LL)$ is a \hdefn{$(w,\ell)$}{layered partition} of $G$ if each $P\in\PP$ is $\ell$-skinny with respect to $\LL$, and the quotient graph $G/\PP$ has treewidth at most $w$. We write \hdefn{$\ell$}{layered partition} instead of $(\ell,\ell)$-layered partition  for simplicity.

Layered partitions are intimately related to graph product structure theory since it follows from the definitions that a graph $G$ has a $(w,\ell)$-layered partition if and only if $G$ is isomorphic to a subgraph of $H\boxtimes P\boxtimes K_\ell$ for some graph $H$ with treewidth at most $w$ and for some path $P$.

If $\mathcal{E}:=(A,\hat{A},G_0,G_1,\dots,G_r)$ is an almost-embedding of a graph $G$, then a layering $\LL$ of $G-A$ is \defn{neat} with respect to $\mathcal{E}$ if the first layer of $\LL$ consists of exactly $(\hat{A}\setminus A)\cup V(G_1\cup\dots\cup G_r)$. In other words, the first layer of $\LL$ consists of the top of $G$ minus the major apex vertices of $G$. This implies that the first two layers of a neat layered partition $\LL$ contain the near-top of $G$ minus the major apex vertices of $G$ (the second layer may also contain some additional vertices).

The next lemma is a mild variant of a result by \citet{DJMMUW20}, and it is one of the key tools that distinguishes our proof from those of \citet{LW1,LW2,LW3}. In particular, we use layered partitions, whereas \citet{LW1,LW2,LW3} use layered treewidth.

\begin{lem}
\label{StronglyAlmostEmbeddableLayeredPartition}
  For every connected $(a,\ell)$-almost-embedded graph $G$ with major apex vertex set $A$, $G-A$ has a $(13\ell,6\ell)$-layered partition $(\PP,\LL)$ that is connected, upward-connected, and neat with respect to the almost-embedding.
\end{lem}

\begin{proof}[Proof Sketch]
  Since \cref{StronglyAlmostEmbeddableLayeredPartition} is not stated explicitly in  \citep{DJMMUW20}, we now explain the small modifications to the proof of \citep[Lemma~28]{DJMMUW20} that are needed to deduce it. The proof of  Lemma 28 in \citep{DJMMUW20} starts by considering the embedded part $G_0$ of the almost-embedding of $G$.  Choose some vertex $v$ of $G_0$ on the boundary of some vortex and make $v$ adjacent to all boundary vertices of all vortices. Since the number, $r$, of vortices is bounded, the resulting graph has bounded Euler genus. The authors of \citep{DJMMUW20} then consider a BFS layering starting at $v$ in the resulting graph, and apply a specific result for graphs on surfaces.  This produces an upward-connected layered partition of $G_0$ in which all the vertices on the boundary of a vortex appear on the first two layers of the resulting layering, and the proof continues by adding all vortex vertices and apex vertices to the first layer.

  To obtain the result stated in \cref{StronglyAlmostEmbeddableLayeredPartition} we start this process slightly differently. Introduce a new vertex $x^*$ adjacent to every vertex of $G_0$ that is on the boundary of a vortex (the resulting graph still has bounded Euler genus). Consider a BFS layering of the resulting graph starting at $x^*$, and then apply the result for graphs on surfaces, as before.  This gives a layered partition in which $x^*$ is in the first layer and the second layer contains precisely the vertices lying on the boundary of a vortex.  After adding the non-major apex vertices and vortex vertices to the second layer and discarding the first layer (that is, $\{x^*\}$), we obtain a neat layering of $G-A$. Moreover, since we started with a BFS layering, each vertex outside the first layer has a neighbour in the previous layer, so the resulting layering is upward-connected.
\end{proof}

A $(k,\ell)$-almost-embedded graph $D$ with major apex vertex set $A$ that is equipped with a neat $\ell$-layered partition $(\PP,\LL)$ of $D-A$ is called a \hdefn{$(k,\ell)$}{drape},
and we say that $D$ is \defn{described by} $(\PP,\LL)$.\footnote{In reality it is rather $D-A$ that is described by $(\PP,\LL)$, but writing $D$ instead of $D-A$ is sometimes convenient, since it allows us to avoid naming the major apex set systematically.} Note that  a $(k,\ell_1)$-almost-embedded graph with a neat $\ell_2$-layered partition is also a $(k,\ell)$-almost-embedded graph with a neat $\ell$-layered partition for $\ell=\max(\ell_1,\ell_2)$, so it is not restrictive to ask for the same $\ell$ in the almost-embedded graph and the layered-partition.

A $(k,\ell)$-pre-curtain $G$ described by a tree-decomposition $\TT:=(B_x:x\in V(T))$ is a \hdefn{$(k,\ell)$}{curtain} if for each $x\in V(T)$, the torso $\torso{G}{B_x}$ is a $(k,\ell)$-drape.
Let $G$ be  a $(k,\ell)$-curtain described by $\TT:=(B_x:x\in V(T))$.  The graph $\mathdefn{\langle G\rangle}:=\bigcup_{x\in V(T)}\torso{G}{B_x}$ is called a \hdefn{$(k,\ell)$}{valance}.\footnote{A \defn{valance} refers to a style of curtain with a fancy decorative top used to hide the curtain rod.}  Note that, under these definitions, $G$ is a $(k,\ell)$-curtain described by $\mathcal{T}$ if and only if $\langle G\rangle$ is a $(k,\ell)$-curtain described by $\TT$.

We say that $G$ is \defn{upward-connected} if, for each $x\in V(T)$, the $\ell$-layered partition $(\PP_x,\LL_x)$ that describes the drape $\torso{G}{B_x}$ is upward-connected. We say that $G$ is \defn{$\PP$-connected} if, for each $x\in V(T)$, the vertex partition $\PP_x$ in the $\ell$-layered partition $(\PP_x,\LL_x)$ that describes $\torso{G}{B_x}$ is a connected partition of $\torso{G}{B_x}$ (so $\torso{G}{B_x}[P]$ is connected, for each $P\in\PP_x$).

A tree of curtains is defined similarly as a tree of pre-curtains, except that each pre-curtain is now replaced by a curtain. Explicitly, a graph $J$ is a \defn{tree of $(k,\ell)$-curtains} if it has a rooted tree-decomposition $\TT:=(B_x:x\in V(T))$ such that:
\begin{compactitem}
  \item for each $x\in V(T)$, $\torso{J}{B_x}$ is a $(k,\ell)$-curtain,
  \item for each edge $xy$ of $T$ where $y$ is the parent of $x$, $B_x\cap B_y$ has size at most $k+3$ and is contained in the top of the root torso of $\torso{J}{B_x}$.
\end{compactitem}
Again, we say that $J$ is a tree of $(k,\ell)$-curtains \defn{described by} $\mathcal{T}$. If every torso in a tree of $(k,\ell)$-curtains is upward-connected, then we say that it is a \defn{tree of upward-connected $(k,\ell)$-curtains}.  If every torso in a tree of $(k,\ell)$-curtains is $\PP$-connected, then we say that it is a \defn{tree of $\PP$-connected $(k,\ell)$-curtains}.

The following is a direct consequence of \cref{Lungs}, since (trees of) curtains are (trees of) pre-curtains.

\begin{cor}\label{Lungs-curtain}
  Let $J$ be a graph that is a tree of $(k,\ell)$-curtains described by $\TT:=(B_x:x\in V(T))$ and such that, for any $(k,\ell)$-curtain $\torso{J}{B_x}$ of $J$, any precolouring of a set $S$ of at most $k+3$ vertices in the top of the root torso of $\torso{J}{B_x}$  can be properly extended to an $m$-colouring of $J[B_x]$ with clustering at most $c$. Then $J$ has an $m$-colouring with clustering at most $c$.
\end{cor}

The following theorem, which is an immediate consequence of \cref{tree_of_raw_curtains,StronglyAlmostEmbeddableLayeredPartition},
provides a structural description of graphs excluding a fixed minor that we believe is of independent interest and might have further applications. Recall that the notion of a \emph{lower-minor-closed} tree-decomposition was defined in Section~\ref{Structure}, shortly before the statement of \cref{tree_of_raw_curtains}.

\begin{thm}
\label{ApexMinorFreeStructure}
For any integer $k\ge 1$ and for any $k$-apex graph $X$, there is an integer $\ell\geq 2$ such that every $X$-minor-free graph is a lower-minor-closed tree of $\PP$-connected upward-connected $(k-1,\ell)$-curtains.
\end{thm}

\begin{figure}[htb]
 \centering
 \includegraphics[scale=1.2]{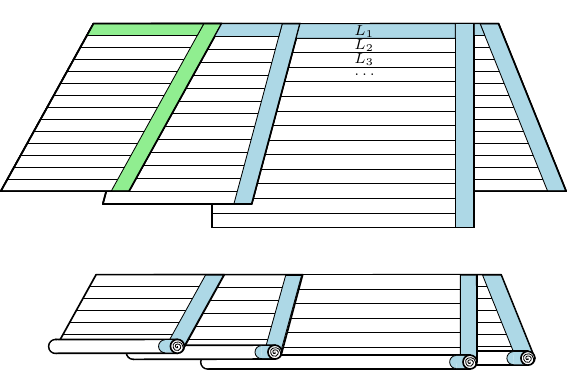}
 \caption{A curtain (above) and a raised curtain (below). We think of the layers as drawn from top to bottom (with the top of the root torso of the curtain being depicted in green). Each torso $\torso{G}{B_x}$ is glued to its parent $\torso{G}{B_y}$ along a subset of the top $A_x \cup L^x_1$ of $\torso{G}{B_x}$ and a subset of the near-top $A_y \cup L^y_1 \cup L^y_2$ of $\torso{G}{B_y}$ (major apex vertices are depicted in blue---or green, for the major apex vertices $A_r$ of the root torso---and are not part of the layerings).}
 \label{fig:curtains}
\end{figure}

Let $G$ be a $\mathcal{P}$-connected upward-connected $(k,\ell)$-curtain described by a rooted tree-decomposition $\TT:=(B_x:x\in V(T))$.  We obtain the raised curtain $G_{\uparrow}$ as follows (see \cref{fig:curtains}): For each $x\in V(T)$, let $(\PP_x,\LL_x)$ be the $\ell$-layered partition that describes the $(k,\ell)$-drape $\torso{G}{B_x}$ and write $\LL_x=(L_1^x,L_2^x,\ldots)$.
For each $x\in V(T)$, let $\mathdefn{G[B_x]_{\uparrow}}$ be the graph obtained from the induced subgraph $G[B_x]$ by contracting each component $C$ of $G[P\cap L^x_{\ge 6}]$
into a single vertex $v_C$, for each $P\in\PP_x$.  Let $L^x_{\uparrow}:=V(G[B_x]_{\uparrow})\setminus(\bigcup_{i=1}^5 L^x_i)$ be the set of vertices obtained from these contractions.\footnote{The value $6$ may seem arbitrary here and, indeed, the definition makes sense for any integer $b\ge 3$ since the sets $L^x_{\ge 3}$, for $x\in V(T)$ are pairwise disjoint.  We fix the value $6$ for the sake of simplicity, because it is the smallest integer that is large enough for all the applications of curtains in the present paper.}

\label{raised_definitions} 
For $x\in V(T)$, let $A_x$ denote the major apex set of $\torso{G}{B_x}$. For each edge $xy\in E(T)$ where $y$ is the parent of $x$, the adhesion set $B_x\cap B_y$ is contained in $(A_x \cup L^x_1)\cap (A_y \cup L^y_1 \cup L^y_2)$. So for distinct $x,y\in V(T)$, the connected subgraphs contracted to create $G[B_{x}]_{\uparrow}$ and the connected subgraphs contracted to create $G[B_{y}]_{\uparrow}$ have no vertices in common. The \defn{raised curtain} is the graph $\mathdefn{G_\uparrow}:= \bigcup_{x\in V(T)} G[B_{x}]_{\uparrow}$ obtained from $G$ by applying all the above edge contractions (for each $x\in V(T)$).

For each $x\in V(T)$, $G[B_x][L^x_{\ge 6}]=\torso{G}{B_x}[L^x_{\ge 6}]$.  In particular, all of the edges contracted in $G[B_x]$ to create $G[B_x]_\uparrow$ are also present in $\torso{G}{B_x}$. Define \defn{$\torso{G}{B_x}_\uparrow$} to be the graph obtained from $\torso{G}{B_x}$ by performing the same edge contractions used to create $G[B_x]_\uparrow$ and define the \defn{raised valance} $\mathdefn{\langle G\rangle_\uparrow}:=\bigcup_{x\in V(T)}\torso{G}{B_x}_\uparrow$.  Clearly, $G[B_x]_\uparrow\subseteq \torso{G}{B_x}_\uparrow$ for each $x\in V(T)$, so $G_\uparrow\subseteq\langle G\rangle_\uparrow$.

\begin{lem}\label{raised_treewidth2}\label{raised_treewidth}
  Let $G$ be a $\PP$-connected $(k,\ell)$-curtain and let $G_{\uparrow}$ be the corresponding raised curtain.
         Then $G_{\uparrow}$ is a minor of $G$
         with treewidth at most $k+6\ell(\ell+1)-1$,
         and $G_{\uparrow}$ is a $(k,\ell)$-curtain
       with at most six layers in each drape.
\end{lem}

\begin{proof}
  Let $\mathcal{T}:=(B_x:x\in V(T))$ be the tree-decomposition that describes the $(k,\ell)$-curtain $G$.  Let $x\in V(T)$.  Then $\torso{G}{B_x}$ is a $(k,\ell)$-drape.  Let $\hat{A}_x$ be the set of (at most $\ell$) apex vertices of $\torso{G}{B_x}$, let $A_x\subseteq \hat{A}_x$ be the set of (at most $k$) major apex vertices of $\torso{G}{B_x}$, and let $(\PP_x,\LL_x)$ be the connected layered partition of
  $\torso{G}{B_x}-A_x$ that describes the $(k,\ell)$-drape $\torso{G}{B_x}$.

  For each $x\in V(T)$, let $(L^x_1,L^x_2,\ldots):=\LL_x$ and let $G^x_{\ge 6}:=\torso{G}{B_x}[L^x_{\ge 6}]$.  Observe that $G^x_{\ge 6}=G[B_x][L^x_{\ge 6}]$ for each $x\in V(T)$ and that the graphs in $\{G^x_{\ge 6}:x\in V(T)\}$ are pairwise vertex disjoint.  It is immediate from the definition that $G_\uparrow$ is a minor of $G$ (because $G_{\uparrow}$ is obtained from $G$ by contracting connected subgraphs of $G_{\ge 6}:=\bigcup_{x\in V(T)} G^x_{\ge 6}$, each of which is a connected subgraph of $G$).

  We now upper bound the treewidth of $G_\uparrow$ by upper bounding the treewidth of the raised valance $\langle G\rangle_\uparrow$, which is a supergraph of $G_\uparrow$.   For each $x\in V(T)$, let $L^x_\uparrow:=V(\torso{G}{B_x}_\uparrow)\setminus (A_x\cup L^x_{\le 5})$, so $\mathcal{L}^x_\uparrow:=(L^x_1,\ldots,L^x_5,L^x_\uparrow)$ is a layering of $\torso{G}{B_x}_\uparrow-A_x$. For each $x\in V(T)$, let $C_x:= V(\torso{G}{B_x}_\uparrow)=A_x\cup L^x_{\le 5}\cup L^x_\uparrow$. For each edge $xy$ of $T$, $C_x\cap C_y=B_x\cap B_y$, so $\mathcal{C}:=(C_x:x\in V(T))$ is a tree-decomposition of $\langle G\rangle_\uparrow$.
  For each $x\in V(T)$, let $H_x:=(\torso{G}{B_x}-A_x)/\PP_x$.
  By the definition of $(k,\ell)$-drape, $\tw(H_x)\le\ell$, so $H_x$ has a tree-decomposition
  $\mathcal{C}_x:=(C^x_z:z\in V(T_x))$ of width at most $\ell$, for each $x\in V(T)$.  We want to use $\mathcal{C}_x$ to construct a tree-decomposition of $\torso{G}{B_x}_\uparrow$.  However, the bags of $\mathcal{C}_x$ contain vertices of $H_x$ (each of which corresponds to a part in $\PP_x$), so we first need to define, for each part $P\in\PP_x$, the set $P_\uparrow$ of vertices that $P$ contributes to $\torso{G}{B_x}_\uparrow$.

  For each $x\in V(T)$ and each $P\in\PP_x$, contracting each component $C$ of $G[P\cap L^x_{\ge 6}]$ produces a vertex $v_C$ in $\torso{G}{B_x}_\uparrow$. For each $x\in V(T)$ and each $P\in\PP_x$, let
  \[
    P_\uparrow:=(P\cap L^x_{\le 5})\cup \{v_C:\text{$C$ is a component of $G[P\cap L^x_{\ge 6}]$}\} \enspace .
  \]
  Observe that, for each $x\in V(T)$, $\mathcal{P}_x:=\{P_\uparrow:P\in\PP_x\}$ is a vertex-partition of $\torso{G}{B_x}_\uparrow-A_x$ and that $H_x = (\torso{G}{B_x}_\uparrow-A_x)/\mathcal{P}_x$.

  We claim that, for each $x\in V(T)$ and each $P\in\PP_x$, $|P_\uparrow|\le 6\ell$. Indeed, $|P\cap L^x_i|\le\ell$ for each $i\in\{1,\ldots,5\}$ so $|P\cap L^x_{\le 5}|\le 5\ell$.  The only remaining issue is the number of components of $G[P\cap L^x_{\ge 6}]$, which we now show is at most $\ell$.  Let $P\in\PP_x$ for some $x\in V(T)$. Since $G$ is a $\mathcal{P}$-connected curtain, $\torso{G}{B_x}[P]$ is connected.
  \begin{itemize}
    \item If $P\subseteq L^x_{\le 5}$ then $G[P\cap L^x_{\ge 6}]$ has zero components.
    \item If $P\subseteq L^x_{\ge 6}$ then $G[P\cap L^x_{\ge 6}]$ has one component.
    \item Otherwise, $P\cap L^x_{\le 5}\neq\emptyset$ and $P\cap L^x_{\ge 6}\neq\emptyset$. Since $\torso{G}{B_x}[P]$ is connected, for each component $C$ of $G[P\cap L^x_{\ge 6}]$, $G$ contains an edge $vw$ with $v\in L^x_{\le 5}$ and $w\in V(C)$.  Since $\LL_x$ is a layering of $\torso{G}{B_x}-A_x$, $v\in L^x_5$ and $w\in L^x_6\cap V(C)$. Therefore, each component of $G[P\cap L^x_{\ge 6}]$ contains a vertex in $L^x_6$. Therefore, the number of components of $G[P\cap L^x_{\ge 6}]$ is at most $|P\cap L^x_6|\le \ell$.
  \end{itemize}

  For each $x\in V(T)$ and each $z\in V(T_x)$, let $D^x_z:=A_x\cup \bigcup_{P\in C^x_z}P_\uparrow$. Then $\bigcup_{z\in V(T_x)}D^x_z=V(\torso{G}{B_x}_\uparrow)=C_x$  and $\mathcal{D}_x:=(D^x_z:z\in V(T_x))$ is a tree-decomposition of $\torso{G}{B_x}_\uparrow$ in which each bag $D^x_z$ has size at most $|A_x|+6\ell |C^x_z|\le k+ 6\ell(\ell+1)$.
  Let $xy$ be an edge of $T$. As observed above, $B_x\cap B_y=C_x\cap C_y$. By the definition of torso, $C_x\cap C_y$ induces a clique $K$ in $\torso{\langle G\rangle_\uparrow}{C_x}=\torso{G}{B_x}_\uparrow$. Since $\mathcal{D}_x$ is a tree-decomposition of $\torso{G}{B_x}_\uparrow$, $T_x$ contains a node $x'$ with $B_x\cap B_y\subseteq D^x_{x'}$.  Similarly, $T_y$ contains a node $y'$ with $B_x\cap B_y\subseteq D^y_{y'}$. By joining $T_x$ and $T_y$ with the edge $x'y'$ we obtain a tree-decomposition of $\torso{G}{B_x}_\uparrow \cup \torso{G}{B_y}_\uparrow$. Doing this for each edge $xy$ of $T$ yields a tree-decomposition $\mathcal{D}$ of $\langle G\rangle_\uparrow$ of width at most $k+ 6\ell(\ell+1)-1$. Therefore $\tw(G_\uparrow)\le\tw(\langle G\rangle_\uparrow)\le k+ 6\ell(\ell+1)-1$.

  It remains to argue that $G_\uparrow$ is a $(k,\ell)$-curtain with six layers in each drape.  We claim that the tree-decomposition $\mathcal{C}$ describes the $(k,\ell)$-curtain $G_{\uparrow}$.  To see this, observe that for each $x\in V(T)$:
  \begin{enumerate}[label=(\alph*)]
      \item\label{still_almost_embedded} $\torso{G_\uparrow}{C_x}=\torso{G}{B_x}_\uparrow$ is a $(k,\ell)$-almost-embedded graph with major apex set $A_x$, since it is obtained from the $(k,\ell)$-almost-embedded $\torso{G}{B_x}$ by contracting edges in its embedded part.  Furthermore $\torso{G}{B_x}_\uparrow$ and $\torso{G}{B_x}$ have the same top and near-top.
      \item\label{still_layered} $\LL^x_\uparrow:=(L^x_1,\ldots,L^x_5,L^x_\uparrow)$ is a layering of $\torso{G_\uparrow}{C_x}-A_x$.
      \item\label{still_partitioned} $\PP^x_\uparrow:=\{P_\uparrow:P\in\PP_x\}$ is a vertex-partition of $\torso{G_\uparrow}{C_x}-A_x$, the quotient $(\torso{G}{B_x}_\uparrow-A_x)/\PP^x_\uparrow=H_x$ has treewidth at most $\ell$, and $|\PP^x_\uparrow\cap L|\le \ell$ for each $L\in\LL^x_\uparrow$.
  \end{enumerate}
  For each edge $xy$ of $T$, $B_x\cap B_y=C_x\cap C_y$.  This fact, combined with observation \ref{still_almost_embedded} ensures that $G_\uparrow$ is a pre-curtain described by $\mathcal{C}$.  Observations \ref{still_layered} and \ref{still_partitioned} show that $\torso{G_\uparrow}{C_x}$ is a $(k,\ell)$-drape with at most six layers that is described by $(\PP^x_\uparrow,\LL^x_\uparrow)$, for each $x\in V(T)$. Therefore $G_\uparrow$ is a $(k,\ell)$-curtain with at most six layers in each drape.
\end{proof}

Let $J$ be a lower-minor-closed tree of $\PP$-connected $(k,\ell)$-curtains described by a tree-decomposition $\TT:=(B_x:x\in V(T))$. Then, for $x\in V(T)$, the torso $\torso{J}{B_x}$ is not necessarily a minor of $J$, which prevents us from directly applying \cref{raised_treewidth} to claim that $\torso{J}{B_x}_\uparrow$ is a minor of $J$. This is the main reason for the introduction of the notion of lower torsos.

\begin{lem}\label{raised_minor}
   Let $J$ be a lower-minor-closed tree of $\PP$-connected $(k,\ell)$-curtains described by a tree-decomposition $\TT:=(B_x:x\in V(T))$. Then, for each $x\in V(T)$, the graph $\ltorso{J}{B_x}_\uparrow$ obtained from $\ltorso{J}{B_x}$ by performing the same contractions used to raise the curtain $\torso{J}{B_x}$ is a minor of $J$.
\end{lem}

\begin{proof}
    Let $x$ be a node of $\TT$.  Then there exists a tree-decomposition $\TT_x:=(C_y:y\in V(T_x))$ that describes  the $\PP$-connected $(k,\ell)$-curtain $\torso{J}{B_x}$.  For each $y\in V(T_x)$, let $A_y$ be the set of major apex vertices in the $(k,\ell)$-almost-embedded graph $\torso{\torso{J}{B_x}}{C_y}$ and let $(\PP_y,\LL_y:=(L^y_1,L^y_2,\ldots))$ be the connected layered partition of $\torso{\torso{J}{B_x}}{C_y}-A_y$ that describes the $(k,\ell)$-drape $\torso{\torso{J}{B_x}}{C_y}$.  Let $C$ be a component of $\torso{\torso{J}{B_x}}{C_y}[P\cap L^y_{\ge 6}]$
    for
    some $P\in\PP_y$.  Then
    \[
       C\subseteq \torso{\torso{J}{B_x}}{C_y}[P\cap L^y_{\ge 6}]
       = \torso{J}{B_x}[P\cap L^y_{\ge 6}]
       = \ltorso{J}{B_x}[P\cap L^y_{\ge 6}]  \enspace .
    \]
    Therefore, each set $V(C)$ induces a connected subgraph of $\ltorso{J}{B_x}$. In other words, $\ltorso{J}{B_x}_\uparrow$ is obtained from $\ltorso{J}{B_x}$ by contracting connected subgraphs, so $\ltorso{J}{B_x}_\uparrow$ is a minor of $\ltorso{J}{B_x}$. Since $J$ is lower-minor-closed, $\ltorso{J}{B_x}$ is a minor of $J$. By transitivity of the minor relation, $\ltorso{J}{B_x}_\uparrow$ is a minor of $J$.
\end{proof}

\subsection{\boldmath \texorpdfstring{\cref{apex1}}{Theorem~7}: Proof for \texorpdfstring{$s\ge 4$}{s~≥~4}}
\label{apex1proof}

Recall that \cref{apex1} states that there is a function $f$ such that for any apex graph $X$ and any integers $t\ge s\ge 3$, any $X$-minor-free $K_{s,t}$-subgraph-free graph is $(s+1)$-colourable with clustering at most $f(X,s,t)$.  Here we prove \cref{apex1} for $t\ge s\ge 4$.  The case $s=3$ requires more sophisticated tools that are also used to establish \cref{Kh,Kst,Jst}, so will be proven later, in \cref{apex1proof3}.

\begin{proof}[Proof of \cref{apex1} for $s\ge 4$]
   Let $X$ be an apex graph and let $J$ be an $X$-minor-free $K_{s,t}$-subgraph-free graph. By \cref{ApexMinorFreeStructure}, $J$ is a lower-minor-closed tree of $\PP$-connected $(0,\ell)$-curtains for some $\ell:=\ell(X)$.  Let $\TT_0:=(B_x:x\in V(T_0))$ be the tree-decomposition that describes $J$.

   Before giving details, we first outline the structure of the proof.  For each $x\in V(T_0)$, instead of colouring the torso $\torso{J}{B_x}$ directly, we colour the lower torso $\ltorso{J}{B_x}$. By definition, $J[B_x]\subseteq \ltorso{J}{B_x}\subseteq \torso{J}{B_x}$, and $\ltorso{J}{B_x}$ contains all the edges that are contracted when raising the curtain $\torso{J}{B_x}$.  When performing the same contractions used to raise $\torso{J}{B_x}$ in $\ltorso{J}{B_x}$, we obtain a subgraph $\ltorso{J}{B_x}_\uparrow$ of the raised curtain $\torso{J}{B_x}_\uparrow$.  Like $\torso{J}{B_x}_\uparrow$,  $\ltorso{J}{B_x}_\uparrow$ has bounded treewidth, by \cref{raised_treewidth}.  Unlike $\torso{J}{B_x}_\uparrow$, $\ltorso{J}{B_x}_\uparrow$ does not contain $K_{s,t'}$-subgraphs for arbitrarily large $t'$ because it is a minor of $J$. This allows us to use \cref{lem:LW0} to properly extend an (appropriately restricted) precolouring of $s$ vertices in $\ltorso{J}{B_x}_\uparrow$ to a full colouring of $\ltorso{J}{B_x}_{\uparrow}$. We finish by using this colouring of $\ltorso{J}{B_x}_\uparrow$ with \cref{component_breaking} to colour $\ltorso{J}{B_x}\supseteq J[B_x]$. The resulting procedure for colouring individual lower torsos $\ltorso{J}{B_x}$ is then used with \cref{Lungs-curtain} to colour the entire tree of curtains $J$.

   Let $x\in V(T_0)$, let $G:=\torso{J}{B_x}$, and let $(C_y:y\in V(T_x))$ be a tree-decomposition that describes the $(0,\ell)$-curtain $G$. For each $y\in V(T_x)$, let $(\PP_y,\LL_y:=(L^y_1,L^y_2,\ldots))$ be the connected layered partition that describes the $(0,\ell)$-drape $\torso{G}{C_y}$. Since $\torso{G}{C_y}$ is a $(0,\ell)$-drape, it does not have any major apices, and thus $C_y=\bigcup_{i\ge 1} L_i^y$.
   Let $H:=\ltorso{J}{B_x}$ be the lower torso. Then $J[B_x]\subseteq H\subseteq G$.

   We will prove that there is a function $c$ such that, for each $x\in V(T_0)$, any $(s+1)$-colouring of a set $S$ of at most $s$ vertices in the top of the root torso of $G=\torso{J}{B_x}$ can be properly extended to an $(s+1)$-colouring of $H:=\ltorso{J}{B_x}$ with clustering at most $c(s,t)$.  Since $J[B_x]\subseteq H$ for each $x\in V(T_0)$,  \cref{Lungs-curtain} then implies that $J$ is $(s+1)$-colourable with clustering $c(s,t)$.

   Perform contractions in the $(0,\ell)$-curtain $G$ to obtain the raised curtain $G_\uparrow$.  Perform the same contractions in $H$ to obtain a graph $H_\uparrow$.  Then $H_\uparrow$ is a minor of $H$ since it is obtained by contracting connected subgraphs of $H$.  (Each contracted subgraph is a component of $\torso{G}{C_y}[P\cap L^y_{\ge 6}]=H[P\cap L^y_{\ge 6}]$, for some $y\in V(T_x)$ and some $P\in\PP_y$.) Furthermore, since $J$ is lower-minor-closed, $H$ is a minor of $J$.

   By \cref{raised_treewidth}, $G_\uparrow$ is a $(0,\ell)$-curtain with treewidth at most $6\ell(\ell+1)-1$ in which each drape has at most six layers.   Let $(D_y:y\in V(T_x))$ be a tree-decomposition that describes $G_\uparrow$. For each $y\in V(T_x)$, let $(\PP_y,\LL_y:=(L^y_1,\ldots,L^y_5,L^y_\uparrow))$ be the connected layered partition that describes the $(0,\ell)$-drape $\torso{G_\uparrow}{D_y}$.

   Let $L_\uparrow:=\bigcup_{y\in V(T_x)}L^y_\uparrow$ and,
   for each $i\ge 1$, let $L_i:=\bigcup_{y\in V(T_x)}L^y_i$.
   ($L_1$ and $L_2$ are not necessarily disjoint.)
   Let $H^{+}_\uparrow$ be the graph obtained from $H_\uparrow$ by adding a single vertex $\alpha$ adjacent to every vertex in $L_\uparrow$, as illustrated in \cref{apex_minor_free_proof}.

   \begin{figure}
    \begin{center}
    \includegraphics{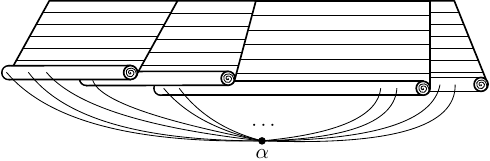}
    \end{center}
    \caption{Adding the vertex $\alpha$ adjacent to the bottom layer of $H_{\uparrow}$.}
    \label{apex_minor_free_proof}
  \end{figure}

  We will apply \cref{lem:LW0} to $H^{+}_{\uparrow}$ with the precoloured set $P:=S\cup\{\alpha\}$ (recall that $S$ is a set of at most $s$ precoloured vertices in the top of the root torso of $G$).  Without loss of generality, we may assume that the vertices in $S$ are coloured with colours in $\{1,\ldots,s\}$ and we choose the colour of $\alpha$ to be $s+1$.  In order to apply \cref{lem:LW0} we must satisfy three conditions:
  \begin{enumerate}
      \item Each vertex of $H^+_\uparrow$ must have at most $s$ precoloured neighbours.  Each vertex of $L_{\uparrow}$ has only one neighbour, $\alpha$, in $P$.  Each vertex in $\bigcup_{i=1}^5 L_i$ has at most $s$ neighbours (all contained in $S$) in $P$. Therefore each vertex of $H^{+}_\uparrow$ has at most $s$ neighbours in $P$.
      \item $H^+_\uparrow$ must have bounded treewidth.  The treewidth of $H^{+}_\uparrow$ is at most $1+\tw(H_\uparrow)\le 1+\tw(G_\uparrow)\le 6\ell(\ell+1)$.

      \item $H^+_\uparrow$ must be $K_{s,t'}$-subgraph-free for some bounded $t'$. Suppose by way of contradiction that $H^+_\uparrow$ contains a $K_{s,t'}$-subgraph $Y$ with parts $U$ of size $s$ and $V$
      of size $t'$, for some large value of $t'$.  We first rule out the possibility that $Y$ is a subgraph of $H$. In this case, $Y$ is also a subgraph of $G$, which has a tree-decomposition $\mathcal{T}_x=(C_y:y\in V(T_x))$.

      If every vertex of $U$ is in the top of $G$ then every edge of $Y$ is incident to a vertex in the top of $G$.  Since $Y\subseteq\ltorso{J}{B_x}$ this implies that $Y\subseteq J[B_x]$, by the definition of lower torso.  For $t'\ge t$, this is not possible, since $J$ is $K_{s,t}$-subgraph-free.  Thus, at least one vertex $u\in U$ is not in the top of $G$.  By symmetry at least one vertex $v$ of $V$ is not in the top of $G$.

      By the definition of (pre-)curtain, every vertex of $G$ that participates in an adhesion of $\mathcal{T}_x$ is in the top of $G$. Therefore neither $u$ nor $v$ participates in any adhesion of $\mathcal{T}_x$, so there is a unique $y\in V(T_x)$ such that $u\in C_y$ and a unique $y'\in V(T_x)$ such that $v\in C_{y'}$.  Since $uv$ is an edge of $\torso{J}{B_x}$, $y=y'$.  Since every vertex in $V$ is adjacent to $u$, $V\subset C_y$. Since every vertex of $U$ is adjacent to $v$, $U\subseteq C_y$.
      To summarize, for $t'\ge t$, we have that $U\cup V\subseteq C_y$, at least one vertex $u\in U$ is not in the top of $G$, and at least one vertex $v$ in $V$ is not in the top of $G$.

      The embedded part of $\torso{G}{C_y}$ is embedded on a surface of Euler genus $g\le \ell$. It follows from Euler's Formula that the embedded part of $\torso{G}{C_y}$ is $K_{4,g+3}$-subgraph-free.  We claim that at least $t'':=t'-g-2\ge 1$ vertices of $V$ are in the top of $\torso{G}{C_y}$. Otherwise, let $w$ be any vertex of $V$ not in the top of $\torso{G}{C_y}$. Then $N_{H}(w)\supseteq U$ is contained in the embedded part of $\torso{G}{C_y}$.  Thus, all vertices of $U$ are in the embedded part of $\torso{G}{C_y}$ and at least $g+3$ vertices of $V$ are not in the top of $\torso{G}{C_y}$.
      Therefore the embedded part of $\torso{G}{C_y}$ contains the vertices and edges of a $K_{s,g+3}$-subgraph $Y'$ of $Y$.  Since $s\ge 4$, this yields a contradiction.

      Thus, there is a subset $V'\subseteq V$ of size at least $t''$ that lies in the top of $\torso{G}{C_y}$.  Let $Y':=Y[U\cup V']$. Then every edge of $Y'$ is incident to a vertex in the top of $G=\torso{G}{B_x}$. As above, this leads to the conclusion that $Y'\subseteq J[B_x]$.
      This is not possible for $t''\ge t$ because $J$ is $K_{s,t}$-subgraph-free.  We conclude that $Y$ is not a subgraph of $H$ (for $t'\ge t+ g+2$).

      Since $Y$ is not a subgraph of $H$, but $Y$ is a subgraph of $H^+_\uparrow$, $Y$ contains a vertex in $V(H^+_\uparrow)\setminus V(H)=\{\alpha\}\cup L_\uparrow$. It follows that $Y$ is a subgraph of  $H^{+}_\uparrow[L_4\cup L_5\cup L_{\uparrow}\cup \{\alpha\}]$.  Each component of $H^{+}_\uparrow[L_4\cup L_5\cup L_{\uparrow}]$ is embedded in a surface of Euler genus at most $g\le \ell$.  It follows from Euler's Formula that  $H^{+}_\uparrow[L_4\cup L_5\cup L_{\uparrow}]$ is $K_{3,2g+3}$-subgraph-free and $K_{4,g+3}$-subgraph-free.

      If $Y$ contains $\alpha$, then $Y-\alpha\subseteq H^{+}_\uparrow[L_5\cup L_{\uparrow}]$ contains a $K_{s,t'-1}$ or a $K_{s-1,t'}$ subgraph.  Since $s\ge 4$, this implies that $t'<\max\{2g+3,g+4\}$.  Otherwise, if $Y$ contains a vertex in $L_\uparrow$ then $Y\subseteq H^{+}_\uparrow[L_4\cup L_5\cup L_{\uparrow}]$, so $t'< g+3$.
      Thus, $H^{+}_\uparrow$ is $K_{s,t'}$-subgraph-free for $t'=\max\{g+4,2g+3,t+g+2\}$.
  \end{enumerate}
  Therefore, $H^{+}_{\uparrow}$ satisfies all the requirements of \cref{lem:LW0}, so $H_{\uparrow}^{+}$ has an $(s+1)$-colouring that properly extends the precolouring of $P$ and has clustering at most $c:=f(6\ell(\ell+1),s,\max\{g+4,2g+3,t+g+2\},s+1)$.

  Extend the colouring of $H^{+}_{\uparrow}$ to a colouring of $H$ in the obvious way: For each vertex $v\in L_\uparrow$, there is a connected subgraph $X_v$ of $G[L_{\ge 6}]$ that was contracted to create $v$.  Assign each vertex in $X_v$ the colour of $v$. The resulting colouring of $H$ does not necessarily have bounded clustering because $X_v$ can be arbitrarily large. However, by the definition of  $(\cdot,\ell)$-curtain, $V(X_v)$ is $\ell$-skinny with respect to the layering $\LL:=(L_1\cup L_2,L_3,L_4,\ldots)$.  This implies that each monochromatic component of $H$ is $c\ell$-skinny with respect to $\LL$. Furthermore, since the colouring of $H^{+}_\uparrow$ properly extends the precolouring of $P=S\cup\{\alpha\}$, each vertex in $L_\uparrow$, and therefore each vertex in $L_{\ge 6}$ avoids the colour $s+1$ assigned to $\alpha$.  Therefore, by \cref{component_breaking} (applied with $t=5$), $H$ has an $(s+1)$-colouring with clustering at most $(2s+6)c\ell$ properly extending the precolouring of $S$, as desired.
\end{proof}

The proof of \cref{apex1} uses many of the elements that are ultimately used in our proof of \cref{Jst}. Specifically, the proof of \cref{Jst} uses the same strategy of decomposing our graph $J$ into a tree of curtains,  and then colouring each curtain $\ltorso{J}{B_x}$ by first
introducing extra vertices to the raised curtain $\ltorso{J}{B_x}_{\uparrow}$ whose colours are  used to break long skinny monochromatic components into short skinny monochromatic components (\cref{component_breaking}).  The main difference is that, in the $\JJ_{s,t}$-minor-free setting, $\torso{J}{B_x}$ is an $(s-3,\ell)$-curtain, so each drape of $\torso{J}{B_x}$ has up to $s-3$ major apex vertices.

The presence of major apex vertices causes considerable complications.  Even a single major apex vertex $a$ in a torso ruins the component breaking strategy used in the proof of \cref{component_breaking}.  To avoid this, we will colour the raised curtain $\ltorso{J}{B_x}_{\uparrow}$ so that, for each component $C$ of $(\ltorso{J}{B_x}_{\uparrow})[L_5\cup L_{\uparrow}]$, the colour assigned to any vertex $v$ in $V(C)\cap L_{\uparrow}$ is distinct from the colour assigned to each neighbour of $v$ in the major apex set $A_y$ of the drape $\torso{\ltorso{J}{B_x}}{C_y}$ that contains $v$, as well as from the colour assigned to an extra vertex $\alpha_C$, whose colour is chosen to be different from the colours of all vertices in $A_y$. This way, the colour of $\alpha_C$ can later be used as a breaker. This motivates the introduction of \defn{admissible sets} that appear in the next section.

However, this highlights another difficulty we will eventually encounter.  This difficulty is most obvious in the context of $K_h$-minor-free graphs, which \cref{Kh} promises to colour with $h-1$ colours.   In this setting, the major apex set $A_y$ of the drape $\torso{\ltorso{J}{B_x}_\uparrow}{C_y}$ that contains $C$ can have size up to $h-5$.  After excluding the $h-4$ colours used by vertices in $A_y\cup\{\alpha_C\}$ this may leave only $(h-1)-(h-4)=3$ colours available for the vertices in $C$.  In general, this would not be possible since the graph $C[L_{\uparrow}]$ may not have a $3$-colouring with bounded clustering. For example, all graphs in the graph class $\mathcal{G}_3$ described in \cref{Standard} are planar and have treewidth $3$, so $(\torso{\ltorso{J}{B_x}_\uparrow}{C_y})[L^y_{\uparrow}]$ could contain arbitrarily large members of  $\mathcal{G}_3$.  In this case, we will be forced to produce a $3$-colouring of $(\torso{\ltorso{J}{B_x}_\uparrow}{C_y})[L^y_{\uparrow}]$ that has arbitrarily large monochromatic components.

Using some new tricks (\cref{skinny_components}) and old tricks (islands), we show that each monochromatic component in $\ltorso{J}{B_x}_{\uparrow}$ corresponds to a monochromatic component in $G$ that is skinny with respect to $\LL$, so it can be broken into bounded size monochromatic components using \cref{component_breaking}. This is the only step of the proof that makes use of the upward-connectivity condition  provided by \cref{ApexMinorFreeStructure}.  This upward-connectivity condition holds for the torso $\torso{J}{B_x}$ but not necessarily for the lower torso $\ltorso{J}{B_x}$.  This requires us to frequently work with both torsos $\torso{J}{B_x}$ and with lower torsos $\ltorso{J}{B_x}$.

\section{Bounded Treewidth Lemmas}
\label{BoundedTreewidth}

This section proves two lemmas about bounded treewidth graphs that exclude certain subgraphs. Both lemmas depend on the following sequence of definitions.

A \defn{list-assignment} of a graph $G$ is a function $L$ such that $L(v)$ is a set of `colours' for each vertex $v\in V(G)$. An \hdefn{$L$}{colouring} is a colouring of $G$ where each vertex $v\in V(G)$ is assigned a colour in $L(v)$.

Fix an integer $s\geq 1$.  For a list-assignment $L$ of a graph $G$, let
\begin{align*}
P(L) & :=\{v\in V(G): |L(v)|=1\}\text{ and }\\
Q(L) & :=\{v\in V(G)\setminus P(L):|N_G(v)\cap P(L)| \in \{1,2,\dots,s-1\}\}.
\end{align*}
The vertices in $P(L)$ are said to be \defn{precoloured} by $L$. We say $L$ is \hdefn{$(s,p)$}{good} if:
\begin{compactenum}[(g1)]
  \item $|P(L)| \leq p$,\label[good]{p_size}
  \item $|L(v)|\ge s+1-|N_G(v)\cap P(L)|$ for all $v\in Q(L)$;\label[good]{q_size}
  \item $|L(v)|\ge 2$ for all $v\in N_G(P(L))\setminus Q(L)$;\label[good]{u_size}
  \item $|L(v)|\ge s+1$ for all $v\in V(G)\setminus N_G[P(L)]$; and\label[good]{other_size}
  \item $L(v)\cap L(u)=\emptyset$ for all $v\in Q(L)$ and $u\in N_G(v)\cap P(L)$.\label[good]{q_proper}
\end{compactenum}

Note that by the definition of $P(L)$, property \cref{u_size} above could be equivalently stated as $L(v)\ne \emptyset$ for all $v\in N_G(P(L))\setminus Q(L)$, or equivalently (using the other properties) as $L(v)\ne \emptyset$ for all $v\in V(G)$.

For two list-assignments $L$ and $L'$ of $G$, we say that $L'$ is a \defn{specialization} of $L$ if $L(v)\supseteq L'(v)$ for each $v\in V(G)$, written as \DefNoIndex{$L\supseteq L'$}. In this case, any $L'$-colouring of $G$ is also an $L$-colouring of $G$.

\begin{figure}[htb]
 \centering
 \includegraphics[scale=1]{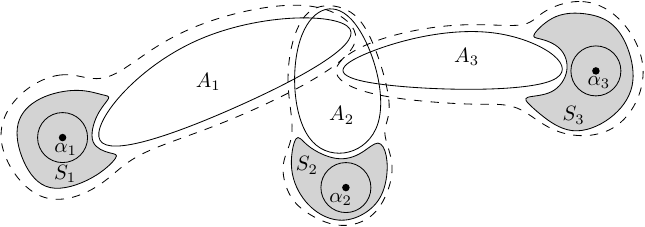}
 \caption{A $k$-admissible set.}
 \label{fig:admissible}
\end{figure}

As illustrated in \cref{fig:admissible}, for an integer $k\geq 0$, a set $\mathcal{S}:=\{(\alpha_1,A_1),\dots,(\alpha_r,A_r)\}$ is \hdefn{$k$}{admissible} in a graph $G$ if:
\begin{compactenum}[(t1)]
\item \label[triple]{a_size}
$A_j$ is a set of at most $k$ vertices in $G$,  for each $j\in\{1,\ldots,r\}$,
\item \label[triple]{no_a_in_X}
$\alpha_j$ is a vertex of $G$
and $N_G[\alpha_j]\cap A_j=\emptyset$,
for each $j\in\{1,\ldots,r\}$,
\item \label[triple]{no_a_in_s}
$(\bigcup_{j=1}^r A_j) \cap (\bigcup_{j=1}^r S_j) = \emptyset$ where \(\mathdefn{S_j}:=N^2_{G-A_j}[\alpha_j]\),
\item \label[triple]{a_in_n_s}
$A_j \subseteq N_G(S_j)$ for each $j\in\{1,\ldots,r\}$, and
\item \label[triple]{s_disjoint}
$S_1,\dots,S_r$ are pairwise disjoint.
\end{compactenum}

To avoid any possible confusion, we emphasise that the neighbourhoods in \cref{no_a_in_X} and \cref{no_a_in_s} above are defined in different graphs ($G$ in \cref{no_a_in_X} and $G-A_j$ in  \cref{no_a_in_s}).

We call $S_j$ the \hdefn{$j$}{trigger set}, since in the  colouring procedure below, colouring the first vertex in $S_j$ triggers the colouring of all uncoloured vertices in $A_j\cup\{\alpha_j\}$.

A list-assignment $L$ and a $k$-admissible set $\mathcal{S}=\{(\alpha_1,A_1),\dots,(\alpha_r,A_r)\}$ for a graph $G$ are \defn{compatible} if for each $j\in\{1,\ldots,r\}$:
\begin{compactenum}[(c1)]
    \item if $S_j \cap P(L)\neq\emptyset$, then $A_j\cup\{\alpha_j\} \subseteq P(L)$;
    \label[tparen]{t_afirst}
    \item if $S_j \cap P(L)\neq\emptyset$, then $L(a)\cap L(\alpha_j)=\emptyset$ for each $a\in A_j$; and \label[tparen]{t_apex_apex_proper}
    \item if $S_j \cap P(L)\neq\emptyset$, then $L(a)\cap L(x)=\emptyset$ for each $a\in A_j\cup\{\alpha_j\}$ and each $x\in N_G(\alpha_j)\cap N_G(a)$.
    \label[tparen]{t_apex_proper}
  \end{compactenum}

In words, whenever a vertex of a trigger set $S_j$ is precoloured, we require that $\alpha_j$ and all the vertices of $A_j$ are also precoloured (with  $\alpha_j$ having a colour distinct from the colours of the vertices from $A_j$), and moreover the colour of each vertex of $A_j\cup\{\alpha_j\}$ does not appear in the list of its neighbours from $N_G(\alpha_j)$.
In our applications, $A_j$ will be the set of major apex vertices of an almost-embedded graph, and $\alpha_j$ is a dummy vertex added to
a raised curtain arising from the embedded part of the almost-embedded graph, such that $\alpha_j$ is adjacent to no apex vertex, and $S_j$ is adjacent to no vortex and no non-major apex vertex. The colour assigned to $\alpha_j$  (which by \cref{t_apex_apex_proper,t_apex_proper} is not used by a vertex in $N_G(\alpha_j) \cup A_j$) then serves as an extra colour that can be used in \cref{component_breaking} to break up large (but skinny) monochromatic components when we lower the curtain.

We use the following result by \citet{LW1}, which has a simple proof implicit in the earlier work of \citet{OOW19}.

\begin{lem}[{\protect\citep[Lemma~6]{LW1}}]\label{heavy_neighbour_bound}
  For every integer $\ell\ge 2$, every $s\in\{1,\ldots,\ell-1\}$, every integer $t\ge s$, there exists an integer $c:=c(s,t,\ell)\geq 1$, such that for every graph $G$ with treewidth less than $\ell$ and with no $K_{s,t}$-subgraph, and for every $P\subseteq V(G)$,
  \[
     \big|\{u\in V(G)\setminus P: |N_G(u)\cap P|\ge s\}\big| \le c|P|\enspace.
  \]
\end{lem}

\subsection{\boldmath \texorpdfstring{$K_{s,t}$}{Kst}-Subgraph-Free Bounded Treewidth Graphs}
\label{BoundedTreewidthKstFree}

We now reach the first of our two lemmas for bounded treewidth graphs. Say a function $f:\RR^+\to\RR^+$ is \defn{reasonable} if $f(x)\geq x$ for all $x\in\NN$, $f$ is increasing, and $\lim_{x\to\infty} x/\log f(x)=\infty$.  The following lemma is the technical extension of \cref{lem:LW0} mentioned earlier.

\begin{lem}\label{lem:LW}
  For any integers $s,t,\ell\geq 2$, there is a reasonable function $f:\RR^+\to\RR^+$ such that for every integer $p\ge f(1)$  the following holds:
  Let $G$ be a graph with treewidth less than $\ell$ and with no $K_{s,t}$-subgraph. Let $L$ be an $(s,p)$-good list-assignment for $G$ that is compatible with an $(s-2)$-admissible set  $\mathcal{S}:=\{(\alpha_1,A_1),\dots,(\alpha_r,A_r)\}$ in $G$.
  Then $G$ has an $L$-colouring $\varphi:V(G)\to \bigcup_{v\in V(G)} L(v)$ with clustering $f(p)$, such that:
\begin{compactenum}[\rm(a)]
  \item the union of the monochromatic components that intersect $P(L)$ has at most $f(|P(L)|)$ vertices;  \label[paren]{p_mono_size}
  \item $\varphi(a)\neq\varphi(\alpha_j)$ for each $j\in\{1,\ldots,r\}$ and $a\in A_j$. \label[paren]{apex_apex_proper}
  \item $\varphi(a)\neq\varphi(x)$ for each $j\in\{1,\ldots,r\}$, $a\in A_j\cup\{\alpha_j\}$ and  $x\in N_G(\alpha_j)\cap N_G(a)$. \label[paren]{apex_proper}
\end{compactenum}
\end{lem}

\citet{LW2} prove a result similar to \cref{lem:LW} with $r=0$ and in the more general setting of tangles (rather than bounded treewidth).  Rewriting their result without reference to tangles (assuming bounded treewidth instead) results in a lemma statement equivalent to \cref{lem:LW} with $r=0$, and a proof that has roughly the same structure as the one given below. The main differences between the two are the technical modifications required to achieve  \cref{apex_apex_proper,apex_proper} when $r>0$.
\begin{proof}
Define $f$ by $f(x):=x^{\log_{4/3}2} 12\ell s (cs+1)^{12\ell s}$, where $c=c(s,t,\ell)$ is defined as in \cref{heavy_neighbour_bound}.  So $f$ is reasonable. Fix $p\geq f(1)$. We prove this result by induction on $(|V(G)|,|V(G)|-|P(L)|)$ in lexicographical order. Since the result is vacuous when $V(G)=\emptyset$, we may assume that $|V(G)|\ge 1$. Let $P:=P(L)$ and $Q:=Q(L)$. There are two easy cases to deal with before arriving at the two main cases.

\paragraph{Case A. $P=\emptyset$:}
If $A_j\neq\emptyset$ for some $j\in \{1,\dots,r\}$,  then let $v$ be any vertex in $A_j$. Otherwise, if $r\ge 1$, then let $v:=\alpha_1$.  Otherwise, let $v$ be any vertex of $G$. Set $P':=\{v\}$ and define $L'(v)$ to be a 1-element subset of $L(v)$. For each $w\in N_G(v)$, let $L'(w):=L(w)\setminus L'(v)$. Since each $L(w)$ has size at least $s+1$, $|L'(w)|\ge s$, as desired. For each vertex $x\in V(G)\setminus N_G[v]$, let $L'(x):=L(x)$. Note that $L'$ is an $(s,p)$-good list-assignment, with $P(L')=\{v\}$ and $Q(L')=N_G(v)$ (since $s\ge 2$).

  We now show that $L'$ is compatible with $\mathcal{S}$. If $r=0$ then \cref{t_afirst}--\cref{t_apex_proper} are vacuous. Otherwise, $v\in A_j\cup\{\alpha_j\}$ for some $j\in\{1,\ldots,r\}$.  If $v=\alpha_j$ then, by the choice of $v$, $A_j=\emptyset$ (and $j=1$), so $A_j\cup\{\alpha_j\}\subseteq \{v\}=P(L')$ which satisfies \cref{t_afirst}; \cref{t_apex_apex_proper} is vacuous; and $Q(L')=N_G(v)$ so \cref{q_proper} implies \cref{t_apex_proper}.  Otherwise $v\in A_j$ so, by \cref{no_a_in_s}, $P(L')\cap S_k=\{v\}\cap S_k=\emptyset$ for each $k\in\{1,\ldots,r\}$, so \cref{t_afirst}--\cref{t_apex_proper} are vacuous.

 Since $|P'|>|P|$, by induction, $G$ has an $L'$-colouring $\varphi$ that satisfies the requirements of the lemma for $(G,L',\mathcal{S})$. Since $L\supseteq L' $ and $P=\emptyset$, $\varphi$ is an $L$-colouring that satisfies the requirements of the lemma for $(G,L,\mathcal{S})$, as desired.

  \paragraph{Case B. $N_G(P)=\emptyset$:}
  Let $G':=G- P$. Since Case A does not hold, $P\ne\emptyset$ and $|V(G')|<|V(G)|$. By induction, there is a colouring $\varphi$ of $G'$ satisfying the conditions of the lemma for $(G',L,\mathcal{S})$ (where we consider the natural restrictions of $L$ and $\mathcal{S}$ to $G'$). Extend $\varphi$ to $G$ by assigning to each vertex $v\in P$ the unique colour in $L(v)$. Since $N_G(P)=\emptyset$, each monochromatic component intersecting $P$ in $G$ is contained in $P$, and thus the union of the monochromatic components intersecting $P$ has size at most $|P|\le p=f(1)\le f(|P|)$, while the size of the remaining monochromatic components is the same in $G$ and $G'$. Therefore $\varphi$ satisfies requirement \cref{p_mono_size}. Since $L$ is compatible with $\mathcal{S}$, \cref{t_apex_apex_proper} and \cref{t_apex_proper} imply that $\varphi$ satisfies \cref{apex_apex_proper} and \cref{apex_proper}, respectively.  Therefore, we obtain an $L$-colouring that satisfies the requirements of the lemma for $(G,L,\mathcal{S})$, as desired.

  \paragraph{Case C. $|P|< 12\ell s$:}
  Let $\hat{p}:=|P|$ and let $v_1,\ldots,v_{\hat{p}}$ be the  vertices of $P$. Define $P_0:=P$ and $Q_0:=Q$.  We will define a sequence of good list-assignments $L=L_0\supseteq L_1\supseteq\cdots\supseteq L_{\hat{p}}=:L'$ that are all compatible with $\mathcal{S}$.  We will ensure that $|P(L')|>|P(L)|$ so that we can apply induction to $(G,L',\mathcal{S})$, and obtain an $L'$-colouring $\varphi$ of $G$ that satisfies \cref{apex_apex_proper,apex_proper}.  Of course, $\varphi$ also satisfies \cref{p_mono_size} with respect to $L'$, but this is not sufficient to ensure that it satisfies \cref{p_mono_size} with respect to $L$ since $|P(L')|>|P(L)|$. The main difficulty, therefore, is to construct $L_1,\ldots,L_{\hat{p}}$ so that the union of the monochromatic components that intersect $P$ has size at most $f(|P|)$.

  We will construct a sequence of precoloured sets $P_1\subseteq \cdots\subseteq P_{\hat{p}}$ and three related sequences $(U^-_i)_{i\in [\hat{p}]}$, $(U_i)_{i\in [\hat{p}]}$, and $(Q_i)_{i\in [\hat{p}]}$ of vertices of $G$ and these will be used to define $L_1,\ldots,L_{\hat{p}}$. In order to ensure that $U^-_1$ is non-empty, it is treated slightly differently than the other sets $U^-_i$:  If $N_G(P)$ contains at least one vertex that is adjacent to at least $s$ vertices of $P=P_0$, then $U_1^-:=\{u\in N_G(P):|N_G(u)\cap P_0|\ge s\}$.  Otherwise, $U_1^-:=\{u\}$ where $u$ is an arbitrary vertex in $N_G(P_0)$. The set $U_1^-$ is well-defined and non-empty because $P$ is non-empty (otherwise Case~A would apply) and $N_G(P)$ is non-empty (otherwise Case B would apply). Consider the following sets (illustrated in \cref{fig:diagram} and explained below):
\begin{align*}
J_1 & :=\{j\in\{1,\dots,r\}:U_1^-\cap S_j \neq\emptyset,\,P_0\cap S_j=\emptyset\}\\
U_1 & := U_1^-\cup \bigcup_{j\in J_1} ((A_j\cup\{\alpha_j\}) \setminus P_0 ) \\
P_1 & :=P_0 \cup U_1\\
Q_1 & :=\{q\in N_G(P_1): |N_G(q)\cap P_1| \in\{1,2,\dots,s-1\}\}\enspace.
\end{align*}
 For each $i=2,3,\dots,\hat{p}$, let
   \begin{align}
 U^-_i &:=\{u\in N_G(P_{i-1}):|N_G(u)\cap P_{i-1}|\ge s\} \label{U_i_defn}\\
J_i & :=\{j\in\{1,\dots,r\}:U^-_i\cap S_j \neq\emptyset,\,P_{i-1}\cap S_j =\emptyset\} \label{J_i_defn}\\
 U_i & := U^-_i\cup \bigcup_{j\in J_i} ( (A_j \cup\{\alpha_j\}) \setminus P_{i-1} )\label{u_i_defn}\\
 P_i & :=P_{i-1}\cup U_i\label{P_i_defn}\\
Q_i & :=\{q\in N_G(P_i): |N_G(q)\cap P_i| \in\{1,2,\dots,s-1\}\}\enspace.\nonumber
  \end{align}
Note that the $i=1$ case only differs from the $i\in\{2,\dots,\hat{p}\}$ case in the definition of $U^-_i$.

\begin{figure}[htb]
 \centering
 \includegraphics[scale=1.2]{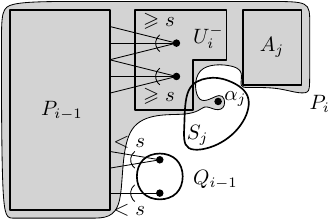}
 \caption{An illustration of the various sets defined in Case C.}
 \label{fig:diagram}
\end{figure}

These definitions warrant some explanation. Suppose we have already constructed $(s,p)$-good list assignments $L_0\supseteq L_1\supseteq\cdots\supseteq L_{i-1}$ where $P_j=P(L_j)$ and $L_j$ is compatible with $\mathcal{S}$ for each $j\in\{0,\ldots,i-1\}$.  Now consider the construction of $L_i\subseteq L_{i-1}$ and $P_i=P(L_i)=P_{i-1}\cup U_i$.  At this point, the set $P_{i-1}$ is precoloured.  The set $U^-_i$ consists of those non-precoloured vertices with at least $s$ precoloured neighbours (step \cref{U_i_defn}), with some exceptions in the $i=1$ case. We are about to precolour all the vertices in $U^-_i$, hence they are added to $P_i$ (step \cref{P_i_defn}). If the $j$-trigger set $S_j$ includes some vertex in $U^-_i$ and no vertex in the $j$-trigger set is currently precoloured, then $j$ is included in $J_i$  (step \cref{J_i_defn}), and we say that $j$ is \defn{triggered}. In this case, $A_j\cup\{\alpha_j\}\setminus P_{i-1}$ is added to $U_i$ and $P_i$ (step \cref{P_i_defn}); subsequently we will also precolour the vertices in $A_j\cup\{\alpha_j\}$ that are not already precoloured. In particular, $\alpha_j$ is added to $P_i$, implying $P_i\cap S_j\neq\emptyset$ and $j$ is triggered at most once (that is, $J_1,\dots,J_{\hat{p}}$ are pairwise disjoint).
  This algorithm ensures that the precoloured set $P_i$ satisfies \cref{t_afirst}; that is, $A_j\subseteq P_i$ whenever $P_i\cap S_j\neq\emptyset$.

  Now it is time to assign colours to the vertices of $U_i=P_i\setminus P_{i-1}$.  For each $j\in J_i$ and each $a\in A_j\cap U_i$ choose any colour in $L_{i-1}(a)$ distinct from the colour of $v_i$.  This is possible since $|L_{i-1}(a)|\ge 2$ by \cref{q_size,u_size,other_size}.  Next, for each $j\in J_i$ choose a colour for $\alpha_j$ that is different from every colour used by any of the at most $(s-2)+1=s-1$ vertices in $A_j\cup\{v_i\}$.  This is possible since $j\in J_i$ implies that $N_G(\alpha_j)\cap P_{i-1}=\emptyset$, so $|L_{i-1}(\alpha_j)|\ge s+1$ by \cref{other_size}.  This ensures that $L_i$ satisfies \cref{t_apex_apex_proper}.

  Before colouring the remaining vertices of $U_i$ we pause to consider each vertex $x\in N_G(\alpha_j)$ for each $j\in J_i$. By definition, $N_G(x)\subseteq A_j\cup\{\alpha_j\}\cup S_j$ and, since $j\in J_i$, $S_j\cap P_{i-1}=\emptyset$.  Therefore, at this moment the only neighbours of $x$ that have been assigned their final colour are in $A_j\cup\{\alpha_j\}$.  We will remove all these colours from $L_{i-1}(x)$ before continuing,  which will ensure that $L_i(x)$ satisfies \cref{t_apex_apex_proper} and still leaves at least $s+1-(s-1)=2$ colours in $L_{i-1}(x)$.  Finally, we assign colours to all the remaining vertices of $U_i$. For each such vertex $v$, \cref{q_size,u_size,other_size} imply that $|L_{i-1}(v)|\ge 2$ and this allows us to choose a colour for $v$ that is distinct from $v_i$.  Lastly, we remove colours from the lists of vertices in $Q_i$ so that they avoid their newly-precoloured neighbours in $U_i$.

  Shortly, we will give a more careful definition of the list-assignments $L_1,\ldots,L_{\hat{p}}$, and prove that each is an $(s,p)$-good list-assignment for $G$ (that is, a list-assignment satisfying \cref{p_size}--\cref{q_proper} above). Before explicitly defining the contents of these lists, we can already verify that they will satisfy \cref{p_size} by bounding the size of $P_i=P(L_i)$ for each
  $i\in\{1,\dots,\hat{p}\}$.  By \cref{heavy_neighbour_bound}, $|U^-_i| \le c|P_{i-1}|$.  By \cref{s_disjoint}, the trigger sets $S_1,\ldots,S_r$ are pairwise disjoint, so each element of $U^-_i$ is responsible for adding the contents of at most one set $A_j\cup\{\alpha_j\}$ (of size at most $s-1$) to $U_i$. Therefore $|U_i|\le s|U^-_i|\le cs|P_{i-1}|$.  (Think of each $u\in U^-_i$ as contributing itself and up to $s-1$ elements of $A_j\cup\{\alpha_j\}$ to $U_i$; for a total of $s$ elements.)  Hence $P_i$ satisfies the recurrence $|P_i| = |P_{i-1}|+|U_i| \le (cs+1)|P_{i-1}|$.  Therefore $|P_i|\le |P_0|(1+cs)^i < 12\ell s(cs+1)^{12\ell s}=f(1)\leq p$.

  We are now ready to carefully define $L_1,\ldots,L_{\hat{p}}$.  Let $L_0:=L$ and for each $i\in\{1,\ldots,\hat{p}\}$, let $c_i$ be the unique element in $L(v_i)$. While defining $L_1,\ldots,L_{\hat{p}}$ we will prove that each is $(s,p)$-good and compatible with $\mathcal{S}$. By assumption, $L_0$ is $(s,p)$-good and compatible with $\mathcal{S}$.   Therefore, in the following definition of $L_i$, we may inductively assume that $L_{i-1}$ is $(s,p)$-good and compatible with $\mathcal{S}$. We start by defining the lists of the vertices in $P_i$, and then we consider vertices in $Q_i$, vertices in  $N_G(P_i)\setminus Q_i$, and finally vertices in $V(G)\setminus N_G[P_i]$.

  \begin{enumerate}
    \item For $v\in P_{i-1}$ set $L_i(v):=L_{i-1}(v)$. Consider a vertex $v\in U_i=P_{i}\setminus P_{i-1}$. We need $L_i(v)$ to be a singleton set, so that $P(L_i)=P_i$. We distinguish three cases, which we handle in the following order:
    \begin{enumerate}[{1.}1]
      \item \label[compound]{easy_list}
      If $v\not\in \bigcup_{j=1}^{r} N_G[\alpha_j]$ then let $L_i(v)$ be a singleton set consisting of an arbitrary element of $L_{i-1}(v)\setminus \{c_i\}$. Since $L_{i-1}$ is $(s,p)$-good and $v\not\in P_{i-1}$, \cref{q_size}--\cref{other_size} imply that $|L_{i-1}(v)|\ge 2$. Therefore $L_{i-1}(v)\setminus\{c_i\}$ is non-empty and $L_i(v)$ is well-defined.

      \item \label[compound]{a_list}
      If $v=\alpha_j$ for some $j\in\{1,\ldots,r\}$ then \cref{t_afirst} implies that $U_i\cap S_j \neq\emptyset$.
      So, $A_j \subseteq P_{i-1}\cup U_i$ by the definition of $U_i$.  In particular, for each $a\in A_j$, $|L_i(a)|=1$ since $a\in P_{i-1}$ (so $|L_{i-1}(a)|=1$) or $a\in U_i$ (so $L_i(a)$ was assigned in \ref{easy_list}, above).
      If $v$ is adjacent to some vertex $u\in P_{i-1}$, then $u\in N_G(\alpha_j) \subseteq S_j$, and $v\in P_{i-1}$ by \cref{t_afirst}, which is a contradiction. Thus $v\not\in N_G(P_{i-1})$. By \cref{other_size}, $|L_{i-1}(v)|\geq s+1$.  In this case, let $L_i(v)$ be a singleton set consisting of an arbitrary colour in $L_{i-1}(v)\setminus (\bigcup_{a\in A_j} L_i(a))\setminus\{c_i\}$.  This is well-defined because $|A_j|\le s-2$ by \cref{a_size}, and $|L_i(a)|=1$ for each $a\in A_j$.  Furthermore, this choice of $L_i(v)$ ensures that $L_i$ satisfies \cref{t_apex_apex_proper}.

      \item Otherwise, $v\in N_G(\alpha_j)$ for some $j\in\{1,\ldots,r\}$ and we must take extra care to ensure that $L_i(v)$ satisfies \cref{t_apex_proper}. In particular, each time we set $L_i(v)$ for some $v\in N_G(\alpha_j)\cap U_i$ we ensure that $A_j\cup\{\alpha_j\}\subseteq P_i$ and that $L_i(v)\cap L_i(a)=\emptyset$ for each $a\in N_G(v)\cap (A_j\cup\{\alpha_j\})$.
      This ensures that \cref{t_apex_proper} is satisfied.
      There are two cases to consider:
      \begin{compactenum}[(i)]
         \item The `typical' case: Some $u\in N_G(P_{i-1})$ has at least $s$ neighbours in $P_{i-1}$.  In this case, the inclusion of $v\in U_i$ implies that $|N_G(v)\cap P_{i-1}|\ge s$. On the other hand,
         $|N_G(v)\cap P_{i-1}| \le |A_j\cup\{\alpha_j\}|+|S_j\cap P_{i-1}|\le s-1+|S_j\cap P_{i-1}|$ by \cref{a_size}.  Therefore, $|S_{j}\cap P_{i-1}|\ge 1$. By \cref{t_afirst,t_apex_proper}, $A_j \cup\{\alpha_j\}\subseteq P_{i-1}$ and $L_{i-1}(a)\cap L_{i-1}(v)=\emptyset$ for each $a\in (A_j \cup\{\alpha_j\} )\cap N_G(v)$. Choose $L_i(v)$ to be a singleton set that contains an arbitrary element of $L_{i-1}(v)\setminus\{c_i\}$.
         As in the previous case, \cref{q_size}--\cref{other_size} ensure that this is well-defined. By \cref{t_afirst,t_apex_proper}, $A_j\cup\{\alpha_j\} \subseteq P_{i-1}$ and $L_{i-1}(a)\cap L_{i-1}(v)=\emptyset$ for each $a\in (A_j \cup\{\alpha_j\} ) \cap N_G(v)$. Therefore, $L_i(v)$ satisfies \cref{t_apex_proper}.

        \item The `exceptional' case: Every $u\in N_G(P_{i-1})$ has at most $s-1$ neighbours in $P_{i-1}$.  If $i\ge 2$, this implies that  $U^-_i$, $J_i$ and $U_i$ are empty, and thus $L_i=L_{i-1}$ and there is nothing to prove. Otherwise, $i=1$, $U^-_1:=\{v\}$, $J_1=\{j\}$, and $|N_G(v)\cap P_0|\in\{1,\dots,s-1\}$. In this case, $U_1=(\{v\}\cup A_j\cup\{\alpha_j\}) \setminus P_0$. This implies that $v\in Q_0$, so \cref{q_proper} implies that $L_0(v)\cap L_0(u)=\emptyset$ for each $u\in  N_G(v)\cap P_0$.
        If $A_j\cup\{\alpha_j\}\subseteq P_0$ then $L_{0}(a)\cap L_{0}(v)=\emptyset$ for each $a\in N_G(v) \cap (A_j\cup\{\alpha_j\})$ by \cref{q_proper}, and we can again let $L_1(v)$ be a singleton set that contains an arbitrary element of $L_{0}(v)\setminus\{c_i\}$. Now assume $A_j\cup\{\alpha_j\}\not\subseteq P_0$.  Then $S_j\cap P_0=\emptyset$ by \cref{t_afirst}.
        Since $N_G(v)\subseteq A_j\cup S_j$, we have
        $N_G(v)\cap P_0=N_G(v)\cap P_0 \cap A_j$.  By \cref{q_size}, $|L_0(v)|\ge s+1-|N_G(v)\cap P_0\cap A_j|$. Since $v\in P_1$ and $v\in N_G(\alpha_j)\subseteq S_j$, we have $A_j\cup\{\alpha_j\}\subseteq U_1\subseteq P_1$. For each $a\in A_j\setminus P_0$, \cref{no_a_in_s} implies that $a\not\in \bigcup_{j'=1}^r N_G[\alpha_{j'}]$, so $L_1(a)$ is assigned in Case~(1.1) above.
        Similarly, $L_1(\alpha_j)$ is assigned in Case~(1.2) above. Thus $|L_1(a)|=1$ for each $a\in A_j\cup\{\alpha_j\}$. Therefore, $|L_0(v)\setminus \bigcup_{a\in A_j\cup\{\alpha_j\}} L_1(a)|\geq (s+1)-|A_j\cup\{\alpha_j\}|\geq 2$ by \cref{a_size}. Let $L_{1}(v)$ be any singleton subset of $(L_0(v)\setminus \bigcup_{a\in A_j\cup\{\alpha_j\}} L_1(a))\setminus\{c_i\}$.
      \end{compactenum}

\medskip We now show that \cref{t_apex_proper} holds for each $x\in N_G(\alpha_j)$ for each $j\in\{1,\dots,r\}$. If $x\not\in N_G(P_i)$ then \cref{t_apex_proper} holds since $L_{i-1}$ and $\mathcal{S}$ are compatible. Now assume that $x\in N_G(\alpha_j)\cap N_G(P_i)$. If $P_i\cap S_j= \emptyset$, then \cref{t_apex_proper} is vacuous for $L_i$, so we can assume that  $P_i\cap S_j\ne \emptyset$, which by \eqref{u_i_defn} implies that $A_j\cup\{\alpha_j\}\subseteq P_i$. We need to ensure that $L_i(x)\cap L_i(a)=\emptyset$ for each $a\in N_G(x)\cap (A_j\cup\{\alpha_j\})$.
      \begin{itemize}
          \item If $x\in Q_i$, then the steps taken in Step~2 below to ensure that $L_i$ satisfies \cref{q_proper} will ensure that $L_i(x)\cap L_i(a)=\emptyset$ for each $a\in N_G(x)\cap (A_j\cup\{\alpha_j\})$.
          \item Otherwise $x\in N_G(P_i)\setminus Q_i$.  If $A_j\cup\{\alpha_j\}\subseteq P_{i-1}$ then since $L_{i-1}$ satisfies \cref{t_apex_proper} and \cref{q_proper} we have $L_i(x)\cap L_i(a)=\emptyset$ for each $a\in N_G(x)\cap (A_j\cup\{\alpha_j\})$. Otherwise $(A_j\cup\{\alpha_j\})\cap U_i\ne \emptyset$, in which case \cref{t_apex_proper} for $L_i$ will also be ensured when defining $L_i(x)$ in Step~3 below.
      \end{itemize}
      Thus $L_i$ will satisfy \cref{t_apex_proper}.  We have already argued above that $L_i$ satisfies \cref{t_afirst,t_apex_apex_proper}, so $L_i$ will be compatible with $\mathcal{S}$.
    \end{enumerate}

At this point, we have already defined $L_i(v)$ for each $v\in U_i$ so that $|L_i(v)|=1$.  This ensures that $P_i\subseteq P(L_i)$. In the following, any further lists we define will have size at least $2$, so $P_i=P(L_i)$.

    \item For each $q\in Q_i$, define $L_i(q):=L_{i-1}(q)\setminus\bigcup_{v\in U_i\cap N_G(q)}L_{i}(v)$.  The fact that $L_{i-1}$ satisfies \cref{q_proper} and that $P_i:=P_{i-1}\cup U_i$ ensures $L_i$ satisfies \cref{q_proper}.  To check that $L_i$ satisfies \cref{q_size}, first observe that $q\not\in N_G(P_{i-1})\setminus Q_{i-1}$ since this would imply that $q\in U_{i}\subseteq P_i$.  Therefore, one of \cref{q_size} or \cref{other_size} implies that $|L_{i-1}(q)|\ge s+1 - |N_G(q)\cap P_{i-1}|$.  Since $P_{i-1}$ and $U_i$ are disjoint, $|N_G(q)\cap P_i|=|N_G(q)\cap P_{i-1}|+|N_G(q)\cap U_i|$.
    Therefore, $|L_i(q)|\ge |L_{i-1}(q)|-|N_G(q)\cap U_i| = s+1-|N_G(q)\cap P_i|$, so $L_i$ satisfies \cref{q_size}.

   \item The purpose of this step is to ensure that \cref{t_apex_proper} will eventually be satisfied for such vertices (even though we are not ready to precolour them).  For each $v\in N_G(P_{i})\setminus Q_i$, we consider two cases:
    \begin{compactenum}[{3.}1]
      \item If  $v\in N_G(\alpha_j)$ and  $N_G(v)\cap (A_j\cup\{\alpha_j\})\cap U_i\ne \emptyset$, then set \[L_i(v):=L_{i-1}(v)\setminus \bigcup_{a\in N_G(v)\cap (A_j\cup\{\alpha_j\})\cap U_i}L_i(a).\] As discussed above, this ensures that $L_i$ satisfies \cref{t_apex_proper}. We simply need to check that $L_i(v)$ still satisfies the conditions of an $(s,p)$-good list-assignment.
    All the lists $L_i(a)$, for $a\in N_G(v)\cap (A_j\cup\{\alpha_j\})\cap U_i$, have been set at Step~1.1 or Step~1.2 and therefore have size 1. It follows that $|L_i(v)|\ge |L_{i-1}(v)|-|N_G(v)\cap (A_j\cup\{\alpha_j\})\cap U_i|$.

\medskip
    If $v\not\in Q_{i-1}$ then
    $v\not\in N_G(P_{i-1})\setminus Q_{i-1}$
    since $v\not\in U_i$. Thus, $|L_{i-1}(v)|\ge s+1$ by \cref{other_size}, and thus $|L_i(v)|\ge s+1-|A_j\cup\{\alpha_j\}|\ge 2$.

    \medskip
    Now assume that $v\in Q_{i-1}$. Then $|L_{i-1}(v)|\ge s+1-|N_G(v)\cap P_{i-1}|$  by \cref{q_size}.
    If  $P_{i-1}\cap S_j\neq\emptyset$, then  $A_j\cup\{\alpha_j\}\subseteq P_{i-1}$  by \cref{t_afirst}, contradicting $(A_j\cup\{\alpha_j\})\cap U_i\ne \emptyset$. Thus $P_{i-1}\cap S_j=\emptyset$. Since $N_G(v)\subseteq S_j\cup A_j$, we have $N_G(v)\cap P_{i-1}\subseteq A_j$. Since $P_{i-1}\cap U_i=\emptyset$, it follows that
    \begin{eqnarray*}
    |L_i(v)|&\ge &|L_{i-1}(v)|-|N_G(v)\cap (A_j\cup\{\alpha_j\})\cap U_i|\\
    &\ge& s+1-|N_G(v)\cap P_{i-1}|-|N_G(v)\cap (A_j\cup\{\alpha_j\})\cap U_i|\\
    &\ge & s+1-|A_j\cup\{\alpha_j\}|\\ &\ge& 2\enspace ,
    \end{eqnarray*}
    so $L_i$ satisfies \cref{u_size}.

    \item Otherwise, set $L_i(v):=L_{i-1}(v)$.
    Since $v\not\in P_{i-1}$, \cref{q_size}--\cref{other_size} imply that $|L_{i-1}(v)| \ge 2$.  Therefore, $|L_i(v)|\ge 2$, so $L_i$ satisfies \cref{u_size}.
    \end{compactenum}

    \item For each $v\in V(G)\setminus N_G[P_{i}]$, define $L_i(v):=L_{i-1}(v)$.  Since $N_G[P_{i}] \supseteq N_G[P_{i-1}]$, \cref{other_size} implies that $|L_{i-1}(v)| \ge s+1$.  Therefore $|L_i(v)|=|L_{i-1}(v)|\ge s+1$ and $L_i$ satisfies \cref{other_size}.
  \end{enumerate}

  Let $L':=L_{\hat{p}}$ and let $P':=P(L')=P_{\hat{p}}$. We previously established that $L'$ is $(s,p)$-good and is compatible with $\mathcal{S}$.  Since $U^-_1\subseteq U_1\subseteq P_1\subseteq P'$ is non-empty, $|P'|\ge |P|+1$.  Therefore, by induction applied to $(G,L',\mathcal{S})$, $G$ has an $L'$-colouring $\varphi:V(G)\to \bigcup_{v\in V(G)} L(v)$ that satisfies the requirements of the lemma for $(G,L',\mathcal{S})$.  Since $L'$ is a specialization of $L$, $\varphi$ is also an $L$-colouring that satisfies \cref{apex_proper}.

  It remains to show that $\varphi$ satisfies \cref{p_mono_size}. To accomplish this, we now prove that any monochromatic component that intersects $P$ is contained in $P'$.  Consider the monochromatic component $C_i$ of $G$ that contains $v_i$ for an arbitrary $i\in\{1,\ldots,\hat{p}\}$, so every vertex in $C_i$ is coloured $c_i$.  We claim that $V(C_i)\subseteq P_{i-1}$.  If not, there are neighbours $w\in P_{i-1}$ and $v\in N(P_{i-1})$ such that $v,w\in V(C_i)$. Since $v\in N(P_{i-1})$, $v$ is either in $U_{i}$ or in $Q_{i-1}$.  If $v\in U_i$ then, by definition $\varphi(v)\neq c_i$. Thus, $v\in Q_{i-1}$ and by \cref{q_proper},  $L_{i-1}(v)\cap L_{i-1}(w)=\emptyset$. Since $\varphi(v)\in L_{i-1}(v)$, $\varphi(w)\in L_{i-1}(w)$, and $\varphi(v)=\varphi(w)$, this is a contradiction. This contradiction shows that $C_i$ is contained in $P_{i-1}\subseteq P'$.  By \cref{p_size}, the union of all monochromatic components intersecting $P$ has size at most $|P'|\le p = f(1)\le f(|P|)$, as required.

  Therefore $\varphi$ is an $L$-colouring of $G$ that satisfies the requirements of the lemma, as desired.

  \paragraph{Case D. $|P|\ge 12\ell s$:}
  \citet[(2.6)]{RS-II} proved that for any graph $G$ of treewidth less than $\ell$ and for any $P\subseteq V(G)$ there exists $V_1,V_2\subseteq V(G)$ such that $G=G[V_1]\cup G[V_2]$ and:
  \begin{align*}
  |V_1\cap V_2|  \le \ell, \qquad
  |P \setminus V_2|  \le 2|P\setminus V_1| \quad \text{ and } \quad
  |P\setminus V_1|  \le 2|P\setminus V_2|.
  \end{align*}
  Apply this result to $G$ and the set $P:=P(L)$ of precoloured vertices. Note that
\begin{equation}
\label{BalSep}
|P\setminus V_1|
\leq  \tfrac23 ( |P\setminus V_1| + |P\setminus V_2| )
\leq \tfrac23 |P|.
\end{equation}
Similarly,
  $|P\setminus V_2|\le\tfrac{2}{3}|P|$.
 Let $Z:= V_1\cap V_2$.

  Let $J:=\{j\in\{1,\ldots,r\}: S_j\cap Z\setminus P\neq\emptyset\}$; that is, $j\in J$ if there is a non-precoloured vertex in the $j$-trigger set and in $Z$. By \cref{s_disjoint}, $|J|\le |Z|\le\ell$. Let
  \[P':=P\cup Z\cup \bigcup_{j\in J} (A_j\cup\{\alpha_j\}).\]
    We now define a list-assignment $L'$ of $G$ that is compatible with $\mathcal{S}$, and such that $P(L')=P'$.
  First, for each $v\in P$, let $L'(v):=L(v)$. 
  By the definition of $J$ and $P'$, and since $P(L')=P'$,
  the first requirement \cref{t_afirst} for $L'$ to be compatible with $\mathcal{S}$ is satisfied.
  \begin{compactenum}
    \item For each $v\in (Z\cup \bigcup_{j\in J} A_j)\setminus (P\cup\bigcup_{j=1}^{r} N_G[\alpha_j])$, set $L'(v)$ to be an arbitrary singleton subset of $L(v)$.  This is analogous to Case C (Step 1.1), above.

    \item For each $j\in J$, let $L'(\alpha_j)$ be an arbitrary singleton subset of $L(\alpha_j)\setminus\bigcup_{a\in A_j} L'(a)$.  As in the analysis of Case~C (Step~1.2) above, this is possible because $|L(\alpha_j)|\geq s+1$, (after the completion of Step~1) $|L'(a)|=1$ for each $a\in A_j$, and $|A_j|\leq s-2$. This step ensures that $L'$ satisfies \cref{t_apex_apex_proper}.

    \item For each $j\in\{1,\ldots,r\}$ and each $x\in N_G(\alpha_j)\cap Z\setminus P$, let $L'(x)$ be an arbitrary singleton subset of $L(x)\setminus \bigcup_{a\in N_G(x)\cap (A_j\cup\{\alpha_j\})}L'(a)$.  To see that this is possible, note that
    $j\in J$ since $x\in S_j\cap Z\setminus P$, and after Steps~1 and 2, $|L'(a)|=1$ for each $a\in A_j\cup\{\alpha_j\}$. Now, as in the analysis of Case~C (Step~1.3), \cref{t_afirst}--\cref{t_apex_proper} and \cref{q_size,other_size} ensure that $|L(x)|\ge s+1-|N_G(x)\cap (A_j\cup\{\alpha_j\})\cap P|$ and \cref{q_proper} implies $L(x)\cap L(a)=\emptyset$ for each $a\in (A_j\cup\{\alpha_j\})\cap P$, so $|L(x)\setminus \bigcup_{a\in N_G(x)\cap (A_j\cup\{\alpha_j\})}L'(a)|\ge s+1-|A_j\cup\{\alpha_j\}|\ge 2$ by \cref{a_size}. Furthermore, this ensures that $L'$ satisfies \cref{t_apex_proper} so that $L'$ is compatible with $\mathcal{S}$.

    \item Let $Q':=\{v\in V(G)\setminus P':|N_G(v)\cap P'| \in\{1,\dots, s-1\}\}$.  For each $q\in Q'$, set $L'(q):=L(q)\setminus\bigcup_{v\in N_G(q)\cap P'} L'(v)$.  As in the analysis of Case~C (Step 2) above, this ensures that $L'$ satisfies \cref{q_size} and \cref{q_proper} and does not make $L'$ incompatible with $\mathcal{S}$.

    \item For each $j\in J$ and $x\in N_G(\alpha_j)\setminus (Z\cup Q'\cup P)$, set 
    \[  
      \textstyle L'(x):=L(x)\setminus(\bigcup_{a\in N_G(x)\cap (A_j\cup\{\alpha_j\})} L'(a)) \enspace .
    \]
    As in the analysis of Case~C (Step 3) above, this ensures that $L'$ satisfies \cref{u_size} and does not make $L'$ incompatible with $\mathcal{S}$.

    \item For any vertex $v\in V(G)$ not covered by any of the preceding cases, set $L'(v):=L(v)$.  This ensures that $L'$ satisfies \cref{u_size,other_size} and does not make $L'$ incompatible with $\mathcal{S}$.
  \end{compactenum}

  Thus $L'$ is a list-assignment for $G$ that satisfies \cref{q_size}--\cref{q_proper} and is compatible with $\mathcal{S}$. However, $L'$ is not necessarily $(s,p)$-good for $G$ since $P'=P(L')$ may be too large to satisfy \cref{p_size}.  Instead, we use divide and conquer. Let
  \[P_1:=(P\setminus V_2) \cup Z \cup\bigcup_{j\in J}( A_j\cup\{\alpha_j\})\quad\text{and} \quad
  P_2:=(P\setminus V_1) \cup Z \cup\bigcup_{j\in J} (A_j\cup\{\alpha_j\}).\]
  Note that $P_1\cup P_2=P'$. By \cref{a_size} and \cref{BalSep},
  \[
    |P_1|
    \le |P \setminus V_2| +|Z| + \sum_{j\in J} |A_j\cup\{\alpha_j\}|
    \le \tfrac23 |P| + \ell + \ell(s-1)
    = \tfrac23 |P| + \ell s.
  \]
  Similarly, $|P_2| \leq \tfrac23 |P| + \ell s$.
  Since $\ell s \leq\frac{|P|}{12}$, for each $i\in\{1,2\}$,
  \[ |P_i|\leq \tfrac34 |P| < |P| \leq p. \]
  For $i\in \{1,2\}$, let $G_i:=G[V_i\cup P_i]$ and let $L'_i$ be the restriction of $L'$ to $V(G_i)$. Thus $L_i'$ satisfies \cref{p_size} in $G_i$ for each $i\in\{1,2\}$.  Since $L'$ satisfies \cref{q_size}--\cref{q_proper} in $G$, $L_i'$ satisfies these properties in $G_i$, for each $i\in\{1,2\}$.  For each $i\in\{1,2\}$, let $\mathcal{S}_i:=\{(\alpha_j, A_j)\in \mathcal{S}: A_j\cup S_j\subseteq V(G_i)\}$.  Then $\mathcal{S}_i$ is $(s-2)$-admissible in $G_i$ for each $i\in\{1,2\}$.  (Note that $\mathcal{S}_1\cup\mathcal{S}_2$ is not necessarily equal to $\mathcal{S}$.) 

  Since $|P_1|<|P|$ and $P_1\supseteq P\cap V_1$, there exists at least one vertex $v\in P\cap V_2$ that is not in $Z\cup\bigcup_{j\in J}(A_j\cup\{\alpha_j\})$, and thus $v$ is not in $G_1$. For each $i\in\{1,2\}$, apply induction to $(G_i,L',\mathcal{S}_i)$ to obtain an $L_i'$-colouring, $\varphi_i$, that satisfies \cref{p_mono_size}--\cref{apex_proper} in $G_i$. Since $V(G_1)\cap V(G_2) \subseteq P'$ and $|L'(v)|=1$ for each $v\in P'$, $\varphi_1(v)=\varphi_2(v)$ for each $v\in V(G_1)\cap V(G_2)$, so these two colourings of $G$ can be combined to produce an $L'$-colouring $\varphi$ of $G$.  
  
  We now argue that $\varphi$ satisfies \cref{apex_apex_proper,apex_proper} with respect to the original admissible set $\mathcal{S}$. 
  For each $i\in\{1,2\}$ and each $j\in\{1,\ldots,r\}$ such that $(\alpha_j,A_j)\in \mathcal{S}_i$, the colours assigned to vertices in $A_j\cup N_G[\alpha_j]$ satisfy \cref{apex_apex_proper,apex_proper} in $G$ since they satisfy \cref{apex_apex_proper,apex_proper} in $G_i$.  For each $j\in\{1,\ldots,r\}$ such that $(\alpha_j,A_j)\not\in\mathcal{S}_1\cup\mathcal{S}_2$, there exists $u\in (A_j\cup S_j)\setminus V_1$ and $w\in (A_j\cup S_j)\setminus V_2$.  By definition $G[S_j]$ is connected and, by \cref{a_in_n_s}, each vertex in $A_j$ has a neighbour in $S_j$.  Thus $G$ contains a path from $u$ to $w$ whose internal vertices all belong to $S_j$ and $Z=V_1\cap V_2$ contains at least one vertex $v$ of this path.  Thus, $Z$ contains a vertex $v\in S_j$.

  If $v\in P$ then $L'(\alpha_j)\cap L'(a)=\emptyset$ for each $a\in A_j$ (by \cref{t_apex_apex_proper}).
  If $v\in Z\setminus P$ then $j\in J$, $|L'(a)|=1$ for each $a\in A_j$ (by Step~1, above) and $L'(\alpha_j)\cap (\bigcup_{a\in A_j} L'(a))=\emptyset$ (by Step~2, above). In either case $L'(\alpha_j)$ and $L'(a)$ are disjoint singletons, so $\varphi(\alpha_j)\neq\varphi(a)$ for each $a\in A_j$, meaning that the pair $(\alpha_j,A_j)$ satisfies \cref{apex_apex_proper}.
  
  If $v\in P$ then $L'(x)\cap L'(a)=\emptyset$  for each $x\in N_G(\alpha_j)$ and each $a\in N_G(x)\cap (A_j\cup\{\alpha_j\})$ (by \cref{t_apex_proper}).  
  If $v\in Z\setminus P$ then, for each $x\in N_G(\alpha_j)$, either $x\in P$ (in which case $S_j\cap P\neq\emptyset$ so $L'(x)=L(x)$ avoids $L'(a)=L(a)$ by \cref{t_apex_proper} for $L$) or $L'(x)$ is defined in one of Steps~3, 4, or 5.  In each of these steps, $L'(x)$ is defined so that $L'(x)\cap L'(a)=\emptyset$ for each $a\in N_G(x)\cap (A_j\cup\{\alpha_j\})$.
  Thus $\varphi(x)\neq\varphi(a)$ for each $a\in N_G(x)\cap (A_j\cup\{\alpha_j\})$, so that the pair $(\alpha_j,A_j)$ satisfies \cref{apex_proper}.
  
  It remains to show that $\varphi$ also satisfies \cref{p_mono_size} and that monochromatic components of $G$ have size at most $f(p)$.  We start by establishing \cref{p_mono_size}. Let $C$ denote the union of the monochromatic components intersecting $P'$ in $G$, and for $i\in \{1,2\}$, let $C_i$ be the union of the monochromatic components intersecting $P_i$ in $G_i$. Then $V(C)\subseteq V(C_1)\cup V(C_2)$, so
  \begin{align*}
    |C| \le |C_1|+|C_2|\le f(|P_1|) + f(|P_2|)     \le 2 f(\tfrac{3}{4}|P|)
    & \le 2\cdot(\tfrac{3}{4} |P|)^{\log_{4/3} 2}12\ell s(cs+1)^{12\ell s} \\
    & = |P|^{\log_{4/3} 2}12\ell s(cs+1)^{12\ell s} \\& = f(|P|) \enspace .
  \end{align*}
  Thus, $\varphi$ is an $L$-colouring of $G$ that satisfies \cref{p_mono_size}.

  Finally, we verify the bound on the clustering. Since $\varphi_1$ and $\varphi_2$ satisfy \cref{apex_apex_proper,apex_proper}, each monochromatic component of $G$ that is completely contained in $G_i$, for some $i\in\{1,2\}$, has size at most $f(p)$.  Any monochromatic component of $G$ not contained entirely in $G_1$ or in $G_2$ contains a vertex in $Z\subseteq P'$, and we have argued above that each such component has size at most $f(|P|)\le f(p)$.  Thus, $\varphi$ is an $L$-colouring of $G$ that satisfies all the requirements of the lemma, which concludes the proof.
\end{proof}

\subsection{\boldmath \texorpdfstring{$\JJ_{s,t}$}{Jst}-Minor-Free Bounded Treewidth Graphs}

Recall that $k$-admissible sets $\{(\alpha_1,A_1),\dots,(\alpha_r,A_r)\}$ were defined at the beginning of \cref{BoundedTreewidth}, and that we defined $S_j:=N^2_{G-A_j}[\alpha_j]$ for each $j\in\{1,\dots,r\}$.

The following is our second lemma for bounded treewidth graphs. In comparison with \cref{lem:LW}, this lemma makes a weaker admissibility assumption, but adds several new requirements related to $K_{2,t}$ and $K_{3,t}$ subgraphs in (or close to) $G[S_j]$.  Each of the sets $S_j$ will come from the larger numbered layers in a layering of some torso, which means that $S_j$ is contained in the embedded part of some torso. This limits the size of $K_{3,t}$-subgraphs, but not the size of $K_{2,t}$-subgraphs.  Therefore, before applying \cref{lem:LW2} we will, in \cref{surface_stuff}, develop some tools for eliminating $K_{2,t}$-subgraphs in embedded graphs.

\begin{lem}\label{lem:LW2}
For any integers $\ell\ge 1$, $s\ge 3$, $t\ge s+2$, and $k\ge s-2$ there exists an integer $p_0$ and a function $f:\mathbb{R}\to\mathbb{R}$ such that the following holds for every integer $p\ge p_0$. Let $G$ be a graph with treewidth less than $\ell$. Let $\mathcal{S}:=\{(\alpha_1,A_1),\dots,(\alpha_r,A_r)\}$ be an $(s-3)$-admissible set in $G$ such that:
  \begin{compactenum}[\rm(x1)]
    \item  $G-\{\alpha_1,\ldots,\alpha_r\}$ is $\JJ_{s,t}$-minor-free;\label[extra]{kh_minor_free}
    \item \label[extra]{sj_connected} for each $j\in\{1,\ldots,r\}$,
    $G[S_j\setminus\{\alpha_j\}]$ is connected;
    \item \label[extra]{k2_t_free} for each $j\in\{1,\ldots,r\}$,
    $G[S_j\setminus\{\alpha_j\}]$ is $K_{2,k-s+3}$-subgraph-free;
    and
    \item \label[extra]{k3_t_free} for each $j\in\{1,\ldots,r\}$,
    $G[N^4_{G-A_j}[\alpha_j] \setminus\{\alpha_j\}]$ is $K_{3,k-s+3}$-subgraph-free.
\end{compactenum}
  Let $L$ be an $(s,p)$-good list-assignment for $G$ that is compatible with $\mathcal{S}$. Then $G$ has an $L$-colouring $\varphi$ with clustering $f(p)$ such that:
  \begin{compactenum}[\rm(a)]
    \item the union of the monochromatic components that intersect $P(L)$ has at most $f(|P(L)|)$ vertices; \label[paren]{p_mono_size2}
    \item $\varphi(a)\neq\varphi(\alpha_j)$ for each $j\in\{1,\ldots,r\}$ and $a\in A_j$; and \label[paren]{apex_apex_proper2}
    \item $\varphi(a)\neq\varphi(x)$ for each $j\in\{1,\ldots,r\}$, $a\in A_j\cup\{\alpha_j\}$ and  $x\in N_G(\alpha_j)\cap N_G(a)$. \label[paren]{apex_proper2}
\end{compactenum}
\end{lem}

\begin{proof}
This proof in the $r=0$ case closely follows the proof of \citet[Lemma~18]{LW3}, which is a statement about $K_h$-minor-free graphs. Here, various modifications are needed to establish \cref{apex_apex_proper} and \cref{apex_proper}, and because this theorem is about $\JJ_{s,t}$-minor-free graphs.

Let $p_0$ be the minimum integer satisfying $p_0\geq s+t-1$ and one other inequality below. Let $f:\mathbb{R}\to\mathbb{R}$ be defined as $f(x):=t^{(x+p_0)/(t-s-1)}$. We prove the result by induction on $|V(G)|$, with this choice of $f$ fixed.
Let $P:=P(L)$, so $|P|\leq p$.
Let $t':=\max\{k+1,t(t-1)\}$.

First consider the easy case in which $G$ is $K_{s,t'}$-subgraph-free.  Since $\mathcal{S}$ is $(s-3)$-admissible, it is also $(s-2)$-admissible. Since $L$ is $(s,p)$-good and compatible with $\mathcal{S}$, we can  apply \cref{lem:LW} and obtain an $L$-colouring $\varphi$ of $G$ that satisfies \cref{apex_apex_proper,apex_proper}, has clustering at most $f'(p)$, and for which the union of monochromatic components that intersect $P$ has size at most $f'(|P|)$, where $f':\RR\to\RR$ is the reasonable function that appears as $f$ in \cref{lem:LW} for $K_{s,t'}$-subgraph-free graphs with treewidth less than $\ell$. Since $f'$ is reasonable, $\lim_{x\to\infty} x/\log f'(x) = \infty$, so $f'(x)$ is asymptotically dominated by the exponential function $t^{x/(t-s-1)}$. Thus there exists an integer $p_0:=p_0(s,t,\ell)$ sufficiently large such that $f(x)=t^{(x+p_0)/(t-s-1)}\ge f'(x)$ for all $x\ge 0$.  It follows that $f(|P|)\ge f'(|P|)$ for all $|P|$, which implies \cref{p_mono_size2}, as required.

We can now assume that $G$ contains a $K_{s,t'}$-subgraph. Consider any copy $H$ of $K_{s,y}$ with $y\ge t'$, implying $y\geq k+1$ and $y\geq t(t-1)$. Let $\{X,Y\}$ be the bipartition of $H$ with $|X|=s$ and $|Y|=y$.

  The following claim is critical in order to establish that requirements \cref{apex_apex_proper} and \cref{apex_proper} do not interfere with the structure of the original proof of \citet{LW3}. It will also be reused in the proof of \cref{Heart2}.

  \begin{clm}\label{no_surface_in_x}
    For each $j\in\{1,\ldots,r\}$,  $(X\cup Y) \cap S_j=\emptyset$.
  \end{clm}

  \begin{clmproof}
    Since $\mathcal{S}$ is $(s-3)$-admissible, $|A_j|\leq s-3$ by \cref{a_size}.

    First suppose that $A_j\cup\{\alpha_j\}\subseteq X$.  Then, since $\alpha_j\in X$, $Y\subseteq N_G(\alpha_j)$.  By \cref{no_a_in_X} and \cref{no_a_in_s},  $X\subseteq N_G(Y)\subseteq N^2_G[\alpha_j]\subseteq S_j\cup A_j$ and $|X\cap S_j|=|X\setminus A_j|\geq |X|-|A_j|\geq s-(s-3)= 3$, which implies that $G[S_j\setminus\{\alpha_j\}]$ contains a $K_{2,y}$-subgraph (since $S_j\cap A_j=\emptyset$ by \cref{no_a_in_s}), contradicting \cref{k2_t_free} (since $y\geq k-s+3$).

    Now suppose that $A_j\cup\{\alpha_j\}\not\subseteq X$ and assume, for the sake of contradiction, that $(X\cup Y) \cap S_j\neq \emptyset$.  Since $A_j\cup\{\alpha_j\}\not\subseteq X$, we have $|(A_j\cup\{\alpha_j\})\cap X| \le |A_j\cup\{\alpha_j\}|-1 \le s-3$. Let $X':=X\setminus (A_j\cup\{\alpha_j\})$ and $Y':=Y\setminus (A_j\cup\{\alpha_j\})$. Thus $|X'|=|X|-|(A_j\cup\{\alpha_j\})\cap X|\geq s-(s-3)=3$ and
    $|Y'|\geq |Y|-|A_j\cup\{\alpha_j\}|\geq y-(s-2)\geq (k+1)-(s-2) = k-s+3$.
    Thus $G[X'\cup Y']$ contains a $K_{3,k-s+3}$-subgraph.
    Since $X\cup Y$ intersects $S_j$, and $S_j\cap A_j=\emptyset$ by \cref{no_a_in_s}, we have that $X'\cup Y'$ intersects $S_j$. By construction, $\alpha_j\not\in X'\cup Y'$. Thus $G[N^2_{G-A_j}[S_j]\setminus\{\alpha_j\}]$ contains a $K_{3,k-s+3}$-subgraph, which contradicts \cref{k3_t_free}.
  \end{clmproof}

  The next claim shows that $X$ does an excellent job of separating $Y$. The proof of this simple result is the only place where the specific structure of graphs in $\JJ_{s,t}$ is used. As in the previous claim, it will be reused later in the proof of  \cref{Heart2}.

  \begin{clm}\label{less_than_t_per}
    Each component of $G-X$ contains at most  $t-1$  vertices in $Y$.
  \end{clm}

  \begin{clmproof}
    Suppose, for the sake of contradiction, that some component $C$ of $G-X$ has at least $t$  vertices in $Y$.     Choose a tree $T$ in $C$ with $|V(T)\cap Y|\geq t$  to lexicographically minimise $(|V(T)\cap\{\alpha_1,\ldots,\alpha_r\}|,|V(T)|)$. This is well-defined since any spanning tree $T$ of $C$ satisfies $|V(T)\cap Y|\geq t$.

    We claim that $V(T)\cap\{\alpha_1,\ldots,\alpha_r\}=\emptyset$. Otherwise $T$ contains $\alpha_j$ for some $j\in\{1,\ldots,r\}$. By \cref{no_surface_in_x}, $\alpha_j\not\in Y$. Now, $N_T(\alpha_j)$ is contained in $S_j\setminus\{\alpha_j\}$, which induces a connected subgraph of $G$ by \cref{sj_connected}. So $G[(V(T)\cup S_j)\setminus\{\alpha_j\}]$ is connected. Let $T'$ be a spanning tree of $G[(V(T)\cup S_j)\setminus\{\alpha_j\}]$. By \cref{no_surface_in_x}, $S_j\cap X=\emptyset$, so $T'$ is a subtree of $C$. By construction, $V(T)\cap Y \subseteq V(T')\cap Y$, so $|V(T')\cap Y|\geq t$ .    By \cref{s_disjoint}, $S_j\setminus\{\alpha_j\}$ avoids $\{\alpha_1,\dots,\alpha_r\}$. So $T'$ has fewer vertices in $\{\alpha_1,\dots,\alpha_r\}$ than $T$, contradicting the choice of $T$.  Hence $V(T)\cap\{\alpha_1,\ldots,\alpha_r\}=\emptyset$.

    If a leaf $x$ of $T$ is not in $Y$, then $T-x$ contradicts the minimality of $T$. So every leaf $x$ of $T$ is in $Y$. If $|V(T)\cap Y|>t$  and $x$ is any leaf of $T$, then $T-x$ again contradicts the minimality of $T$. So $|V(T)\cap Y|=t$.

    Let $T_1,\dots,T_t$ be a partition of $T$ into $t$ pairwise disjoint subtrees, each with exactly one vertex in $Y$.
    Now $|Y|\geq t(t-1)\geq t+s$ since $t\geq s+2$.
    Since $|V(T)\cap Y|=t$, there is a matching $\{a_1b_1,\dots,a_sb_s\}$ where $\{a_1,\dots,a_s\}=X$ and $b_1,\dots,b_s \in Y\setminus V(T)$, which avoids $\alpha_1,\dots,\alpha_r$ by \cref{no_surface_in_x}.
    Let $G'$ be the minor of $G$ obtained by contracting each $T_i$ into a vertex $z_i$, and contracting each edge $a_ib_i$ into a vertex $x_i$.  The subgraph $G'[\{z_1,\dots,z_t\}]$ is connected and thus contains a spanning tree.
    Since $a_ib_j$ is an edge of $G$ for distinct $i,j\in\{1,\dots,s\}$, in $G'$, the set $\{x_1,\dots,x_s\}$ is a clique, and it is complete to $\{z_1,\dots,z_t\}$. Hence, $G'[\{x_1,\dots,x_s,z_1,\dots,z_t\}]$ contains a graph in $\JJ_{s,t}$.
    Since $V(T)\cap\{\alpha_1,\ldots,\alpha_r\}=\emptyset$ and the matching $\{a_1b_1,\dots,a_sb_s\}$ avoids $\alpha_1,\dots,\alpha_r$, we have found a  graph in $\JJ_{s,t}$ as a minor of $G-\{\alpha_1,\dots,\alpha_r\}$. This contradiction shows that every component of $G-X$ has at most $t-1$ vertices in $Y$.
  \end{clmproof}

Note that the proofs of \cref{no_surface_in_x} and \cref{less_than_t_per} do not use \cref{p_size}, the upper bound on $|P(L)|$ implied by the fact that $L$ is $(s,p)$-good.  We use this fact later, in the proof of \cref{Heart2}, to apply both these claims to a graph $G^{+-}_{\uparrow\bullet}$ that contains a $K_{s,q}$-subgraph, (for a sufficiently large value of $q$) and satisfies all the requirements of \cref{lem:LW2} except for this upper bound on $|P(L)|$.

For the remainder of the proof, consider a copy $H$ of $K_{s,y}$ in $G$ that maximises $y$ (in particular, $y\ge t'$ and the two claims above apply). We again write the bipartition of $H$ as $\{X,Y\}$, with $|X|=s$ and $|Y|=y$.

  Let $C_1,\ldots,C_q$ be the components of $G-X$. By \cref{less_than_t_per}, each of $C_1,\ldots,C_q$ has at most $t-1$ vertices in $Y$. So $q\geq \frac{|Y|}{t-1}=\frac{y}{t-1} \geq t\geq s+2$. We now consider two cases.

  \paragraph{Case A. Some component $\boldsymbol{C_i}$ of $\boldsymbol{G-X}$ contains no vertex in $\boldsymbol{P}$:}

  Write $C:=C_i$ for brevity. Let $G':=G-V(C)$, and let $\mathcal{S}'$ be the restriction of $\mathcal{S}$ to $G'$.  Since $G'$ is a subgraph of $G$, it satisfies \cref{kh_minor_free,k2_t_free,k3_t_free}.  By \cref{no_surface_in_x} and \cref{sj_connected}, $G[S_j]=G'[S_j]$ or $G'[S_j]$ is empty for each $j\in\{1,\ldots,r\}$.  In the former case, $G'$ satisfies \cref{sj_connected}. In the latter case, $(\alpha_j,A_j)$ plays no role in $\mathcal{S}'$, so $G'$ also satisfies \cref{sj_connected}. Since $G'$ has fewer vertices than $G$, we may apply induction to $(G',L',\mathcal{S}')$, so $G'$ has an $L'$-colouring $\varphi'$ that satisfies \cref{p_mono_size}, \cref{apex_apex_proper}, and \cref{apex_proper}. We show how to extend $L'$ to a colouring of $G$.  To ensure that monochromatic components in $G'$ do not grow, our colouring of $C$ will properly extend $\varphi'$ to $G$.

  Let $G_C:=G[V(C)\cup X]$. Then $G_C$ satisfies \cref{kh_minor_free}--\cref{k3_t_free}   for the same reasons $G'$ does. Let $U_C:=\{v\in V(C):|N_G(v)\cap X|=s\} =Y \cap V(C)$ by the maximality of $y$. By \cref{less_than_t_per}, $|U_C|\le t-1$ and, by \cref{no_surface_in_x}, $U_C\cap S_j=\emptyset$ for each $j\in\{1,\ldots,r\}$.  Let $P_C:=X\cup U_C$. Observe that, for each $v\in V(C)$, $N_{G}(v)\cap P\subseteq X$. So  $|L(v)|\ge s+1-|P\cap X|$  by \cref{q_size}--\cref{other_size}.  This is helpful in verifying the claims about $|L_C(v)|$ for the list-assignment $L_C$ that we now define:
  \begin{compactitem}
    \item For each $x\in X$, let $L_C(x):=\{\varphi'(x)\}$.

    \item For each of the at most $t-1$ vertices $v\in U_C$, let $L_C(v)$ be a singleton subset of $L(v)\setminus\bigcup_{x\in X\cap N_{G_C}(v)} L_C(x)$.  This is possible since $v$ has exactly $s$ neighbours in $X$.

    \item Let $Q_C:=\{v\in N_{G_C}(P_C):|N_{G_C}(v)\cap P_C|\in\{1,\dots,s-1\}\}$. For each $v\in Q_C$, let $L_C(v):=L(v)\setminus\bigcup_{w\in N_{G_C}(v)\cap P_C} L_C(w)$.  The fact that $L$ satisfies \cref{q_size,other_size,q_proper} implies that $|L_C(v)|\ge s+1-|N_{G_C}(v)\cap P_C|$ for each $v\in Q_C$.

    \item For each $v\in N_{G_C}(P_C)\setminus Q_C$, let $L_C(v):=L(v)\setminus\bigcup_{x\in X\cap N_G(v)} L_C(x)$.  Observe that $v\not\in P_C$, so $v$ is adjacent to at most $s-1$ vertices in $X$. Since $L$ satisfies \cref{q_size}--\cref{q_proper} and $|L_C(x)|=1$ for each $x\in X$, this implies that $|L_C(v)|\ge 2$ for each $v\in  N_{G_C}(P_C)\setminus Q_C$.

    \item For each $v\in V(C)\setminus N_{G_C}[P_C]$, let $L_C(v):=L(v)$.  Since $L$ satisfies \cref{other_size}, $|L_C(v)|=|L(v)|\ge s+1$.
  \end{compactitem}
    The only vertices assigned singleton lists in $L_C$ are in $X\cup U_C=P_C$, so $P_C=P(L_C)$ and, by definition, $Q_C=Q(L_C)$. We now verify that, for any $p\ge p_0 \geq s+1$, $L_C$ is an $(s,p)$-good list-assignment for $G_C$.
    \begin{compactenum}[(g1)]
      \item[\cref{p_size}:] $|P_C| = |X|+|U_C| \le s+t-1 \le p_0 \le p$;
      \item[\cref{q_size}:] $|L_C(v)|\ge s+1-|N_G(v)\cap P_C|$ for each $v\in Q_C$, as discussed above;
      \item[\cref{u_size}:] $|L_C(v)|\ge 2$ for all $v\in N_{G_C}(P_C)\setminus Q_C$, as discussed above;
      \item[\cref{other_size}:] $|L_C(v)|\ge s+1$ for all $v\in V(G_C)\setminus N_{G_C}[P_C]$ since $L_C(v)=L(v)$; and
      \item[\cref{q_proper}:] $L_C(v)\cap L_C(u)=\emptyset$ for all $v\in Q_C$ and $u\in N_{G_C}(v)\cap P_C$.
   \end{compactenum}
   Each of the rules above explicitly ensures that $L_C(v)\cap L_C(x)=\emptyset$ for each $x\in X\cap N_{G_C}(v)$.  Therefore, any $L_C$-colouring of $G_C$ properly extends the precolouring of $X$ given by $\varphi'$.
   Let $\mathcal{S}_C$ be the restriction of $\mathcal{S}$ to $G_C$. By \cref{no_surface_in_x}, $(X\cup U_C)\cap S_j=\emptyset$ for each $j\in\{1,\dots,r\}$. So $S_j\cap P(L_C)\neq\emptyset$ if and only if $S_j\cap P(L)\neq\emptyset$ and $S_j\subseteq V(C)$. So $L_C$ trivially satisfies the requirements for being compatible with $\mathcal{S}_C$.

   Since $G-X$ has at least two components,  $|V(G_C)|<|V(G)|$. By induction applied to $(G_C,L_C,\mathcal{S}_C)$, there exists an $L_C$-colouring $\varphi_C$ of $G_C$ that satisfies \cref{p_mono_size}, \cref{apex_apex_proper}, and \cref{apex_proper}.  For each $x\in X$, $\varphi_C(x)=\varphi'(x)$, so $\varphi_C$ and $\varphi'$ can be combined to give a colouring $\varphi$ of $G$.  The colouring $\varphi$ inherits \cref{apex_apex_proper} and \cref{apex_proper} from $\varphi'$ and $\varphi_C$.  It also satisfies \cref{p_mono_size} since each vertex in $N_{G_C}(X)$ is assigned a colour distinct from the colours of its neighbours in $X$, so the monochromatic components that intersect $P$ are contained in $G'$.  Therefore, the colouring $\varphi$ of $G$ satisfies all requirements of the lemma.

  \paragraph{Case B. Each component $\boldsymbol{C_i}$ of $\boldsymbol{G-X}$ contains at least one vertex in $\boldsymbol{P}$:}

In this case, $X$ plays a  similar role to that of the separator $Z$ used in Case~D of the proof of \cref{lem:LW}. For each $i\in\{1,\ldots,q\}$, let $G_i:=G[V(C_i)\cup X]$ and let $\mathcal{S}_i$ be the restriction of $\mathcal{S}$ to $G_i$.  Since $G_i$ is a subgraph of $G$ it satisfies \cref{kh_minor_free,k2_t_free,k3_t_free} and, by \cref{no_surface_in_x,sj_connected}, $G_i$ and $L_i$ also satisfy \cref{sj_connected}.

  Let $L_1$ be the restriction of $L$ to $G_1$.  Then $L_1$ is $(s,p)$-good for $G_1$ (since $L$ is $(s,p)$-good for $G$) and $L_1$ is compatible with $\mathcal{S}_1$ (since $L$ is compatible with $\mathcal{S}$).  Let $P_1:=P(L_1)$.  Observe that $|P_1|\le |P|-q+1$ since $P_1$ avoids at least one vertex of $P$ in each component of $G-X$ other than $C_1$. Recall that $q\ge 2$, so $G_1$ has fewer vertices than $G$. Thus we may apply induction  to $(G_1,L_1,\mathcal{S}_1)$, so $G_1$ has a colouring $\varphi_1$ that satisfies \cref{apex_apex_proper,apex_proper} and such that the union of monochromatic components of $G_1$ that intersect $P_1$ has size at most $f(|P_1|)\le f(|P|-q+1)$.

  We now define a list-assignment $L_i$ for each $i\in\{2,\ldots,q\}$. This is similar to the list-assignment $L_C$ used in the proof of Case A except that we cannot prevent monochromatic components in $P$ from entering $C_i$ (since $P$ already contains vertices of $C_i$).  Instead, we use the fact that $P_i:=P(L_i)$ is much smaller than $P$ to control the growth of these components.

  For each $i\in\{2,\ldots,q\}$, let
  $P_i:= (P\cap V(G_i))\cup X$, let
  $Q_i:=\{v\in N_{G_i}(P_i):|N_{G_i}(v)\cap P_i|\in\{1,\dots,s-1\}\}$,
  and define the list-assignment $L_i$ as follows:
  \begin{compactitem}
    \item For each $x\in X$, let $L_i(x):=\{\varphi_1(x)\}$.

    \item For each $v\in Q_i$, let $L_i(v):=L(v)\setminus\bigcup_{w\in P_i\cap N_{G_i}(v)} L_i(w)$. The fact that $L$ satisfies \cref{q_size,other_size,q_proper} ensures that  $|L_i(v)|\ge s+1-|P_i\cap N_{G_i}(v)|$.

    \item For each $v\in V(G_i)\setminus (P_i\cup Q_i)$, let $L_i(v):=L(v)$.
  \end{compactitem}
  From these definitions, it follows immediately that $P_i:=P(L_i)$ and $Q_i:=Q(L_i)$.  We now verify that $L_i$ is an $(s,p)$-good list-assignment for $G_i$.
  \begin{compactenum}[(g1)]
    \item[\cref{p_size}:] $|P_i| \le |P|-q+|X|+1 = |P|+s-q+1 \le |P|\le p$ (since $q\geq s+1$);

    \item[\cref{q_size}:] $|L_i(v)|\ge s+1-|N_{G_i}(v)\cap P_i|$ for each $v\in Q_i$;

    \item[\cref{u_size}:] $|L_i(v)|\ge 2$ for each $v\in N_{G_i}(P_i)\setminus Q_i$, since $L_i(v)=L(v)$;

    \item[\cref{other_size}:] $|L_i(v)|\ge s+1$ for all $v\in V(G_i)\setminus N_{G_i}[P_i]$, since $L_i(v)=L(v)$; and

    \item[\cref{q_proper}:] $L(v)\cap L(u)=\emptyset$ for all $v\in Q_i$ and $u\in N_{G_i}(v)\cap P_i$, by definition.
  \end{compactenum}
  Let $\mathcal{S}_i$ be the restriction of $\mathcal{S}$ to $G_i$.  Since $P_i\setminus P\subseteq X$,  \cref{no_surface_in_x} implies that any $j$-trigger set intersected by $P_i$ is also intersected by $P$, so $L_i$ is compatible with $\mathcal{S}_i$ because $L$ is compatible with $\mathcal{S}$.

  We are almost done. By assumption, $P$ contains at least one vertex in each of $C_1,\ldots,C_q$, so $|P_i|\le |P|+s+1-q$.  Since $G_i$ has fewer vertices than $G$, we may apply induction to $(G_i,L_i,\mathcal{S}_i)$, so $G_i$ has an $L_i$-colouring $\varphi_i$ that satisfies \cref{apex_apex_proper} and \cref{apex_proper}.  Furthermore, the union of monochromatic components of $G_i$ that intersect $P_i$ has size at most $f(|P|+s+1-q)$.

  Since $\varphi_i(x)=\varphi_j(x)$ for each $i,j\in\{1,\ldots,q\}$ and each $x\in X$, the colourings $\varphi_1,\ldots,\varphi_q$ can be combined to give a colouring $\varphi$ of $G$ that satisfies \cref{apex_apex_proper} and \cref{apex_proper}.  In this colouring, each monochromatic component that intersects $P$ intersects $P_i$ for some $i\in\{1,\ldots,q\}$.  Therefore, the total number of vertices in all monochromatic components that intersect $P$ is at most
  \begin{align*}
    q\cdot f(|P|+s+1-q)
    & = q\cdot t^{\tfrac{|P|+p_0+s+1-q}{t-s-1}} \\
    &= t^{\tfrac{|P|+p_0+s+1-q}{t-s-1} + \tfrac{\log q}{\log t}} \\
    &= t^{\tfrac{|P|+p_0}{t-s-1} + \tfrac{s+1-q}{t-s-1}+\tfrac{\log q}{\log t}} \\
    &\le t^{\tfrac{|P|+p_0}{t-s-1}} & \text{(proved below)} \\
    &= f(|P|) \enspace .
  \end{align*}
  The above inequality is proved as follows. Since $q\geq t\geq s+2\geq 5$, we have $t(q-s-1) \ge q(t-s-1)$ and $\tfrac{q}{\log q} \ge \tfrac{t}{\log t}$, which together imply $\tfrac{q-s-1}{t-s-1}\ge \tfrac{q}{t} \ge \tfrac{\log q}{\log t}$ and $\tfrac{s+1-q}{t-s-1}+\tfrac{\log q}{\log t}\leq 0$. Thus $\varphi$ is an $L$-colouring of $G$ that satisfies \cref{p_mono_size}, \cref{apex_apex_proper}, and \cref{apex_proper}, as desired.
\end{proof}

\section{Proofs of the Main Results}

This section completes the proofs of \cref{Jst} and the $s=3$ case of  \cref{apex1}. A key ingredient is the next section on
$K_{2,t}$-subgraphs.

\subsection{\boldmath \texorpdfstring{$K_{2,t}$}{K2t}-Subgraphs in Embedded Graphs}
\label{surface_stuff}
\label{K2tSubgraphs}

This section introduces some old and some new tools for eliminating $K_{2,t}$-subgraphs in surface-embedded graphs, so that we can eventually apply \cref{lem:LW2}.

A \hdefn{$d$}{island} of a graph $G$ is a non-empty subset $I$ of $V(G)$ such that $|N_G(v)\setminus I|\le d$ for each $v\in I$. A $d$-island $I$ of $G$ is a \hdefn{$(d,c)$}{island} of $G$ if $|I|\le c$.  The concept of islands was introduced in the context of clustered colouring by \citet{EO16}, who made the following simple but important observation.

\begin{obs}[\citep{EO16}] \label{island_colouring}
For any $(d,c)$-island $I$ of a graph $G$, any $(d+1)$-colouring of $G-I$ with clustering at most $c$ can be properly extended to a $(d+1)$-colouring of $G$ with clustering at most $c$ by colouring each vertex $v\in I$ by a colour not used on $N_G(v)\setminus I$.
\end{obs}

Recall \cref{Genusg} by \citet{DN17}, which says that every  graph of bounded Euler genus has a $4$-colouring with bounded clustering. To prove this, \citet{DN17} showed that graphs of bounded Euler genus have $3$-islands of bounded size and applied \cref{island_colouring} inductively (see \citep[Section~3.5]{WoodSurvey} for a concise exposition). Their proof uses strongly sublinear separators, and we adapt it to show the following.

\begin{lem}\label{three_islands}
  There is a constant $c>0$ such that for every graph $G$ of Euler genus at most $g$ and for every non-empty set $Y\subseteq V(G)$ such that $|Y|\ge 4|N_G(Y)|+4g$, there exists a $(3,c(g+1))$-island $I$ of $G$ with $I\subseteq Y$.
\end{lem}

\begin{proof}
  The following lemma is folklore (see Lemma~15 in \citep{WoodSurvey}): Let $G$ be a graph such that for fixed $\alpha>0$ and $\beta\in(0,1)$, every subgraph $H$ of $G$ has a balanced separator of size at most $\alpha|V(H)|^{1-\beta}$. Then, for every $\epsilon\in(0,1]$, there exists $S \subseteq V(G)$ of size at most $\epsilon |V(G)|$ such that each component of $G - S$ has at most $f(\alpha,\epsilon,\beta):=\ceil{2 ( \frac{\alpha}{\epsilon (2^\beta-1)} )^{1/\beta}}$ vertices. Here, a \defn{balanced separator} in a graph $G$ is a set $S\subseteq V(G)$ such that every component of $G-S$ has at most $\frac12|V(G)|$ vertices.

  Let $G$ and $Y$ be as defined in the lemma. We may delete vertices not in $N_G[Y]$. Now assume that $V(G)=N_G[Y]$. The above lemma with $\beta=\frac12$ and $\alpha=\Theta(\sqrt{g+1})$ is applicable to $G$ and its subgraphs by a separator result of \citet{GHT-JAlg84}. We now apply the method of \citet{DN17}. By the above lemma applied to $G[Y]$ with $\epsilon := \frac{1}{16}$, there exists $X\subseteq Y$ of size at most $\tfrac{1}{16} |Y|$ such that if $K_1,\dots,K_r$ are the components of $G[Y\setminus X]$, then each $K_i$ has at most $f(\alpha,\frac{1}{16},\tfrac12)$ vertices, which is at most $c(g+1)$ for some constant $c$. Let $e(K_i)$ be the number of edges of $G$ with at least one endpoint in $K_i$. By Euler's formula,
  \begin{align*}
    \sum_{i=1}^r e(K_i)
    & \leq |E(G)| < 3(|V(G)|+g) = 3(|Y|+|N_G(Y)|+g) \\
    & \leq \tfrac{15}{4} |Y|
    \leq 4 \, |Y \setminus X |
    = 4 \sum_{i=1}^r |V(K_i)| .
  \end{align*}
  Hence, $e(K_i) < 4\,|V(K_i)|$ for some $i$. Repeatedly remove vertices from $K_i$ with at least four neighbours outside $K_i$. Doing so maintains the property that $e(K_i) < 4\,|V(K_i)|$. Thus, the final set $I$ is non-empty. By construction, $I\subseteq Y$, and $I$ is a $3$-island of $G$ of size at most $|V(K_i)| \leq c(g+1)$.
\end{proof}

Consider the following loose interpretation of \cref{three_islands} obtained by setting $S:=N_G(Y)$.  If $S$ is small compared to $V(G-S)$, then $V(G-S)$ contains a small $3$-island of $G$.  In our application, the set $S$ will consist of a few parts in a layered partition $(\PP,\LL)$ of $G$.  Unfortunately, a few parts in a layered partition $(\PP,\LL)$ of a bounded genus graph $G$ may contain the vast majority of the vertices of $G$. The following lemma shows that removing a constant number of (arbitrarily large) parts from a layered partition $(\PP,\LL)$ of a bounded genus graph $G$ leaves a graph that contains $(3,c)$-islands of $G$ or that is skinny with respect to $\LL$.

\begin{lem}\label{skinny_components}
  Let $\ell\geq 1$ and $g\geq 0$ be integers, and let $c$ be the constant from \cref{three_islands}. Let $G$ be a graph of Euler genus $g$, let $\LL$ be a layering of $G$, and let $S$ be a subset of $V(G)$ that is $\ell$-skinny with respect to $\LL$.  Then $V(G-S)$ is $(44\ell+4g)$-skinny with respect to $\LL$ or $V(G-S)$ contains a $(3,c(g+1))$-island of $G$.
\end{lem}

\begin{proof}
  Assume for the sake of contradiction that $V(G-S)$ has no $(3,c(g+1))$-island of $G$ and that $|V(G-S)\cap L|\geq 44\ell + 4g$ for some layer $L\in\LL$. For any set $R\subseteq\LL$, let $C_R$ be the set of all vertices in $G-S$ that are in layers in $R$. Let $R$ be a consecutive set of layers with
  \begin{equation}
      \label{CR-condition}
      | C_R| \geq  4\ell(|R|+10) + 4g
  \end{equation}
  and $|R|$ maximum. This is well-defined since $R = \{L\}$ satisfies this condition ($|C_{\{L\}}| \geq 4\ell(1+10) + 4g$).

  If $|C_R| \geq 4 |N_G(C_R)| + 4g$ then $C_R$ contains a 3-island of $G$ with size at most $c(g+1)$ by \cref{three_islands}. So
  $ |C_R| < 4 |N_G(C_R)| + 4g$.

  If $R=\LL$, then $C_R=V(G-S)$ and thus $N_G(C_R)\subseteq S$ and $|N_G(C_R)|\leq |S| \leq \ell \,|R|$. It follows that
  \[
    |C_R|  < 4|N_G(C_R)|+4g \leq 4\ell \,|R|+4g,
  \]
  which contradicts the choice of $R$. So $R\neq\LL$. Let $R'$ be the set of layers in $R$ plus the layers before and after $R$ (at least one of these exists). By the choice of $R$, $R'$ fails \cref{CR-condition}. Thus
  \begin{align*}
        |C_{R'}|
   < 4\ell (|R'|+10) + 4g
   \leq 4\ell (|R|+12) + 4g
   \leq |C_R| + 8\ell.
  \end{align*}
  Hence
  \begin{align*}
    |C_R| + 8\ell
     >   |C_{R'}|  & \geq |C_R| + |N_G(C_R)| - |N_G(C_R) \cap S|\\
    & > |C_R| + \tfrac14 |C_R|  - g - \ell (|R|+2)\\
    &  = \tfrac54  |C_R| - g - \ell (|R|+2).
  \end{align*}
Therefore
  $ \tfrac14 |C_R|  <  8\ell + g + \ell (|R|+2)$ and
  $  |C_R|  <  4g + 4\ell (|R|+10)$,
  which contradicts the choice of $R$.
\end{proof}

For any surface $\Sigma$ and any connected subset $\Delta\subseteq\Sigma$, let $\partial\Delta$ denote the boundary of $\Delta$ and let $\interior(\Delta):=\Delta\setminus{\partial\Delta}$ denote the interior of $\Delta$.  We say that two subsets $\Delta_1$ and $\Delta_2$ of $\Sigma$ \defn{overlap} if $\interior(\Delta_1)\cap\interior(\Delta_2)\neq\emptyset$, and $\Delta_1$ and $\Delta_2$ \defn{strictly overlap} if they overlap, but neither $\interior(\Delta_1)$ nor $\interior(\Delta_2)$ contains the other.

In the following, we treat the vertices and edges of a graph $G$ embedded in a surface $\Sigma$ as points and curves in $\Sigma$, respectively. Consequently, paths in $G$ are treated as simple curves in $\Sigma$ and  cycles
in $G$ are treated as simple closed curves in $\Sigma$.  For a subset $R\subseteq\Sigma$, let $G[R]:=G[V(G)\cap R]$ and $G-R:=G[V(G)\setminus R]$.\label{minus_def} Here we include an edge $uv$ in $G[R]$ if and only if both $u$ and $v$ are in $R$ (the edge $uv$ is not necessarily contained in $R$).
A \hdefn{$d$}{disc} $\Delta$ in $G$ (with respect to $\Sigma$) is a closed topological disc of $\Sigma$ whose boundary is a cycle of length $d$ in $G$.  Note that each contractible cycle of length $d$ in $G$ bounds at least one $d$-disc.

\begin{lem}\label{k2t_4disc}
  For any integer $g\ge 0$ and for any embedding of $K_{2,3g+2}$ in a surface $\Sigma$ of Euler genus $g$, there exists a $4$-disc in $K_{2,3g+2}$ with respect to $\Sigma$.
\end{lem}

\begin{proof}
  It is known that if $x$ and $y$ are two vertices in a graph $G$ that is cellularly embedded in $\Sigma$, and $P_1,\ldots,P_k$ are internally disjoint paths with ends $x$ and $y$, such that for any $1\le i<j\le k$, the cycle $P_i\cup P_j$ is non-contractible in $\Sigma$, then $k\le 3g+1$  (see \cite[Proposition 4.2.7]{MoharThom}).

  Now consider any (not necessarily cellular) embedding  of a copy of $K_{2,3g+2}$ in $\Sigma$. Let $\{\{x,y\},\{v_1,\ldots,v_{3g+2}\}\}$ be the bipartition of $K_{2,3g+2}$. Assume for the sake of contradiction that there is no pair $i<j$ such that the embedding contains a 4-disc bounded by $xv_iyv_j$. Triangulate all faces of the embedding that are not topological discs, so that each resulting triangular face bounds a topological disc (see for example \cite[Section 3.1]{MoharThom}).

  The resulting supergraph $G$ is cellularly embedded in $\Sigma$. The copy of $K_{2,3g+2}$ under consideration is a subgraph of $G$, with the property that the paths $P_i:=xv_i y$ of $K_{2,3g+2}$ are internally disjoint, and  $P_i\cup P_j$ is a non-contractible cycle
  for any $i<j$, which contradicts the above property.
\end{proof}

\begin{figure}
    \begin{center}
      \begin{tabular}{ccc}
    \includegraphics[page=1,width=.3\textwidth]{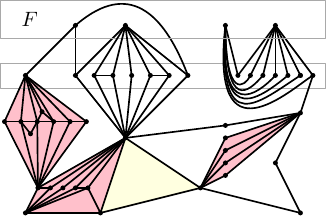} &
    \includegraphics[page=2,width=.3\textwidth]{ggpgpp.pdf} &
    \includegraphics[page=3,width=.3\textwidth]{ggpgpp.pdf} \\
    $G$ & $G_{\circ}$ & $G_{\bullet}^-$
      \end{tabular}
    \end{center}
    \caption{The graphs $G$, $G_{\circ}$ and $G_{\bullet}^-$ in \cref{cover_discs_i,cover_discs_ii}.}
    \label{ggpgpp}
\end{figure}

The following lemma, illustrated in \cref{ggpgpp}, will allow us to eliminate $K_{2,t}$-subgraphs that contain vertices in the bottom layer $L^x_{\uparrow}$ of the raised drape $G^{x}_{\uparrow}$ without disturbing the upper layers $L^x_1,\ldots,L^x_{5}$.  (This lemma will be applied with $F:=F^x:=L^x_{\le 4}$ so that $N_{G_{\uparrow}}[F]=L^x_{\le 5}$.)

Consider a graph $G$ embedded in a surface $\Sigma$ and a subset $F\subseteq V(G)$. We say that a $d$-disc $\Delta$ in $G$ with respect to $\Sigma$ is \hdefn{$F$}{avoiding} if $\Delta$ has no vertex of $F$ on its boundary and no vertex of $N_G[F]$ in its interior. An \hdefn{$F$}{simplification} $G_{\circ}$ of $G$ is obtained as follows: let $\mathcal{C}$ be a set of closed topological discs in $\Sigma$ such that:
\begin{compactenum}[(a)]
    \item each $\Delta\in \mathcal{C}$ is a $3$-disc or a $4$-disc in $G$ with respect to $\Sigma$;\label[paren]{cover_discs_first}
    \item no pair of discs in $\mathcal{C}$ overlap;
    \item every disc in $\mathcal{C}$ is $F$-avoiding;
    \item $|V(G)\cap (\bigcup_{\Delta\in \mathcal{C}} \interior(\Delta))|$ is maximised.\label[paren]{cover_discs_last}
\end{compactenum}
Then $G_{\circ}:=G- (\bigcup_{\Delta\in \mathcal{C}}\interior(\Delta))$ is an \hdefn{$F$}{simplification} of $G$.

\begin{lem}\label{cover_discs_i}
  Let $G$ be a graph embedded in a surface $\Sigma$ of Euler genus $g$, let $F\subseteq V(G)$ be such that $G[F]$ has at most $c$ connected components, and let $G_{\circ}$ be an $F$-simplification of $G$.
  Then there exists an integer $t:=t(g,c)$ such that $G_{\circ}-F$ is $K_{2,t}$-subgraph-free.
\end{lem}

\begin{proof}
  Let $\mathcal{C}$ be the set of discs used for the given $F$-simplification $G_\circ$.  A \defn{chord} in a $d$-disc $\Delta$ of $G$ having a boundary cycle $v_1,\ldots,v_d$ is an edge $v_iv_j$ of $G$ with endpoints in $\{v_1,\ldots,v_d\}$ such that the interior of $v_iv_j$ is contained in the interior of $\Delta$.  A $d$-disc is \defn{chord-free} if it has no chords.   We may assume that each $4$-disc in $\mathcal{C}$ is chord-free since each non-chord-free $4$-disc $\Delta$ can be replaced by two interior disjoint $3$-discs $\Delta_1$ and $\Delta_2$ whose union is $\Delta$ and whose intersection is a chord of $\Delta$. This operation preserves \cref{cover_discs_first}--\cref{cover_discs_last} and does not change $V(G)\cap (\bigcup_{\Delta\in \mathcal{C}}\interior(\Delta))$.

 We may assume that $\bigcup_{\Delta\in \mathcal{C}} \Delta$ contains every $F$-avoiding chord-free $4$-disc of $G$. Indeed, the only reason not to include such a disc $\Delta$ would be that it overlaps one or more other discs in $\mathcal{C}$. Since all discs in $\mathcal{C}$ are $3$-discs or chord-free $4$-discs, no disc in $\mathcal{C}$ strictly overlaps $\Delta$.\footnote{This uses the fact that if two chord-free $4$-discs $\Delta_1$ and $\Delta_2$ strictly overlap, then $\Delta_1\cup\Delta_2$ is also a chord-free $4$-disc whose interior contains strictly more vertices than $\interior(\Delta_1)\cup\interior(\Delta_2)$, so $\mathcal{C}$ would contain neither $\Delta_1$ nor $\Delta_2$.}  Therefore, all discs in $\mathcal{C}$ that overlap $\Delta$ are contained in $\Delta$, so we can remove them all and replace them with $\Delta$.  This does not change $V(G)\cap (\bigcup_{\Delta\in \mathcal{C}}\interior(\Delta))$ since this set is already of maximum size, and thus it does not change the simplification $G_\circ$.

  Suppose that $G_{\circ}-F$ contains a subgraph $Z$ isomorphic to $K_{2,t}$. Let $X:=\{v,w\}$ and $Y:=\{y_1,\ldots,y_t\}$ be the two parts in the bipartition of $V(Z)$.  Define a graph $J$ with vertex set $V(J):=Y$ that contains an edge $y_iy_j$ if and only if the cycle $vy_iwy_j$ is contractible.  By Ramsey's Theorem \cite{Ramsey30}, there exists a sufficiently large integer $t:=R(2c+6,3g+2)$ such that the graph $J$ contains a clique of size $2c+6$ or a stable set of size $3g+2$.  In the latter case, this stable set $Y'$ defines a subgraph of $G_{\circ}[X\cup Y']$ isomorphic to $K_{2,3g+2}$ such that no $4$-cycle in this subgraph bounds a $4$-disc. This contradicts \cref{k2t_4disc}. Therefore, we may assume that $J$ contains a clique of size $2c+6$ with vertex set $Y'$.

  Let $k:=2c+6$.
  Then $G_{\circ}[X\cup Y']$ contains a subgraph $L$ isomorphic to $K_{2,k}$ in which every $4$-cycle bounds a $4$-disc of $L$ with respect to $\Sigma$.  By \cref{lem:planardisc}, it follows that $L$ is embedded in a closed $4$-disc $\Delta$ of $L$ with respect to $\Sigma$. By renaming the vertices of $Y'$ to $y_1,\ldots,y_k$ if necessary, we may assume that the boundary of $\Delta$ is the cycle $v y_1 w y_k$ and, for each $i\in\{1,\ldots,k-1\}$, the 4-cycle $v y_i w y_{i+1}$ bounds a face of $L$.

  For each $i\in\{1,\ldots,\floor{(k-1)/2}\}$, consider the $4$-disc $\Delta_{2i}$ bounded by $v y_{2i-1} w y_{2i+1}$ that contains $y_{2i}$ in its interior.  Since $y_{2i}\in V(G_{\circ})$, $\Delta_{2i}\not\in \mathcal{C}$.  Therefore, $\Delta_{2i}$ contains a vertex of $N_G[F]$ in its interior (and is therefore not $F$-avoiding) or $\Delta_{2i}$ contains the edge $vw$ (and is therefore not chord-free).\footnote{The fact that $\Delta_{2i}$ cannot contain the edge $y_{2i-1}y_{2i+1}$ follows from the Jordan Curve Theorem.}

  \begin{compactitem}
    \item Since $G$ is a simple graph, there is at most one value of $i$ such that $\Delta_{2i}$ contains $vw$.

    \item If $\Delta_{2i}$ contains a vertex of $N_G[F]$ in its interior then, since $V(Z)\cap F=\emptyset$, the interior of $\Delta_{2i}$ contains an entire component of $F$. Therefore there are at most $c$ values of $i$ for which $\Delta_{2i}$ contains a vertex of $N_G[F]$.
  \end{compactitem}
  Therefore, $\floor{(k-1)/2}\le c+1$ which implies that $k< 2c+6$, a contradiction.  Therefore, $G_{\circ}-F$
  does not contain a subgraph isomorphic to $K_{2,t}$ for $t\ge R(2c+6,3g+2)$.
\end{proof}

For distinct vertices $v$ and $w$ in a graph $G$, we say that $v$ is \defn{dominated} by $w$ (in $G$) if $N_G(v)\subseteq N_G(w)$. (Note that by this definition a vertex cannot dominate its neighbour). If $v$ is dominated in $G$ by $w$, then $w$ \defn{dominates} $v$ in $G$.  We say that $v$ is \defn{dominated} in $G$ if there exists $w\in V(G)\setminus\{v\}$ that dominates $v$ in $G$.  Two properties of this definition are immediate:
\begin{inparaenum}[(i)]
  \item If $v$ is a dominated vertex of a connected graph $G$, then $G-v$ is also connected.
  \item If $v$ is dominated by $w$ in $G$ then $v$ is dominated by $w$ in any induced subgraph of $G$ that includes both $v$ and $w$.
\end{inparaenum}
For any $R\subseteq V(G)$ and any $v\in V(G)$ we say that $v$ is \hdefn{$R$}{dominated} in $G$ if there exists some $w\in R\setminus\{v\}$ that dominates $v$ in $G$.

\begin{obs}\label{no_3_dominated}
  Let $G$ be a graph embedded in the plane having an outer face bounded by a cycle $F$ in $G$.  Let $R:=V(G)\setminus V(F)$ be such that $N_G(x)\subseteq V(F)$ for each $x\in R$. Then any $R$-dominated vertex in $R$ has degree at most $2$.
\end{obs}

\begin{proof}
  If $x\in R$ has three or more neighbours in $F$, then the only face of $G[V(F)\cup\{x\}]$ with three neighbours of $x$ on its boundary is the outer face.  Thus, no vertex $y\in R$ can dominate $x$.
\end{proof}

\begin{obs}\label{non_dominated}
  Let $G$ be a graph embedded in the plane having an outer face bounded by a cycle $F$ in $G$ and such that $N_G(x)\subseteq V(F)$ for each $x\in V(G)\setminus V(F)$. Let $v$ and $w$ be distinct vertices of $F$, and let $R:=N_G(v)\cap N_G(w)\setminus V(F)$.  Then $|R|\le 2$ or some vertex of $R$ is $R$-dominated in $G$.
\end{obs}

\begin{proof}
  Suppose that $R$ contains three vertices $x$, $y$, and $z$.  By definition, each of these vertices is adjacent to both $v$ and $w$, and thus has at least two neighbours in $F$. If any of these vertices has exactly two neighbours in $F$, then it is dominated by the other two and we are done, so assume each of $x,y,z$ has a neighbour in $F$ that is distinct from $v$ and $w$. By symmetry, we can assume that $y$ lies in the interior of the 4-disc bounded by $vxwz$. Since this cycle only intersects $F$ in $\{v,w\}$, it follows from the Jordan Curve Theorem that $y$ has no neighbour in $F$ distinct from $v$ and $w$, which is a contradiction.
\end{proof}

For an embedded graph $G$ and $F\subseteq V(G)$, an \hdefn{$F$}{contraction} $G_{\bullet}$ of $G$ is obtained as follows:  Let $G_{\circ}$ be an $F$-simplification of $G$ obtained by a set $\mathcal{C}$ of $3$- and $4$-discs. In the original graph $G$, contract each component $C$ of $G-V(G_{\circ})$ into a single vertex $v_C$, and call the resulting graph $G_{\bullet}$. An \hdefn{$F$}{compression} $G_{\bullet}^-$ of $G$ is obtained from $G_{\bullet}$ by repeatedly deleting any vertex $v\in V(G_{\bullet})\setminus V(G_{\circ})$ that is dominated in the current graph by some vertex $w\in V(G_{\bullet})\setminus V(G_{\circ})$ that belongs to the same disc $\Delta\in\mathcal{C}$ as $v$, until no such dominated vertex remains.

To make the end of this section more concrete, we explain how the notions of $F$-contraction and $F$-compression will be used in \cref{JstMinorFree}. Once we have raised a curtain, our aim  is to clean its embedded part to avoid large copies of $K_{2,t}$. It is safe to remove components $C$ that are attached to at most two vertices of the embedded part of the raised curtain, because these can be coloured easily using our strategy of breaking skinny components using one additional colour. These components $C$ correspond to the
vertices of $G_\bullet-V(G_\bullet^-)$ (to which they have been contracted), and \cref{no_3_dominated,non_dominated} tell us that they are indeed attached to the embedded part by at most two vertices. We now prove that once these vertices are removed, the maximum size of a copy of $K_{2,t}$ disjoint from $F$ is bounded.

\begin{lem}\label{cover_discs_ii}
  Let $G$ be a graph embedded in a surface $\Sigma$ of Euler genus $g$, let $F\subseteq V(G)$ be such that $G[F]$ has at most $c$ connected components, and let $G_{\bullet}^-$ be an $F$-compression of $G$.
  Then there exists an integer $t:=t(g,c)$ such that $G_{\bullet}^--F$ is $K_{2,t}$-subgraph-free.
\end{lem}

\begin{proof}
  Let $G_{\circ}$ be the $F$-simplification of $G$ used to create $G_{\bullet}^-$, obtained from some set $\mathcal{C}$ of $3$- and $4$-discs, and let $R:=V(G_{\bullet}^{-})\setminus V(G_{\circ})$.
  The vertices in $R$ form an independent set in $G_{\bullet}^{-}$ and each vertex in $R$ is adjacent to at most four vertices on the boundary cycle of some $3$- or $4$-disc in $\mathcal{C}$. Therefore, each vertex in $R$ has degree at most $4$.  Let $Z$ be a $K_{2,t}$-subgraph in $G_{\bullet}^{-}$ with parts $X:=\{v,w\}$ and $Y:=\{y_1,\ldots,y_{t}\}$. Then each vertex in $X$ has degree at least $t$. Therefore, $\{v,w\}\subset V(G_{\circ})$ or $t\le 4$.  In the latter case, we are done, so now assume that $\{v,w\}\subset V(G_{\circ})$.

  Suppose that $V(Z)\cap F=\emptyset$.   By \cref{cover_discs_i}, $|Y\cap V(G_{\circ})|\le t'$ for some $t':=t'(g,c)$. Therefore, $|Y\cap R|\ge t-t'$.  By \cref{non_dominated}, for each disc $\Delta\in \mathcal{C}$ there are at most two vertices of $R$ that are adjacent to both $v$ and $w$ in the interior of $\Delta$.  Therefore, there are at least $(t-t')/2$ discs in $\mathcal{C}$ that contain both $v$ and $w$ on their boundary cycles.  There are at most two such discs in which $v$ and $w$ are consecutive in the boundary cycles.  Therefore, there are $d\ge (t-t')/2-2$ $4$-discs in $\mathcal{C}$ with $v$ and $w$ on, but not consecutive on, their boundary cycles.  Each of these $d$ discs defines two length-$2$ paths in $G_{\circ}$ from $v$ to $w$ and each of these length-$2$ paths is shared by at most two discs. Therefore, $G_{\circ}-F$ contains a $K_{2,d}$-subgraph, so $d< t'$. Thus, $(t-t')/2-2 \le d < t'$ so $t< 3t'+4$.  Therefore, $G_{\bullet}^--F$ does not contain a $K_{2,t}$-subgraph for $t:=3t'+4$.
\end{proof}

\subsection{\boldmath \texorpdfstring{\cref{Jst}: $\JJ_{s,t}$}{Jst}-Minor-Free Graphs}
\label{JstMinorFree}

This section finishes the proof of our main theorem, \cref{Jst}, which says that $\JJ_{s,t}$-minor-free graphs are $(s+1)$-colourable with bounded clustering. We use the same general strategy used in the proof of \cref{apex1} (for $s\ge 4$) given in \cref{apex1proof}, with some additional steps.  In \cref{JstCurtain} we show that every $\mathcal{J}_{s,t}$-minor-free graph $J$ is a lower-minor-closed tree of $\PP$-connected upward-connected $(s-3,\ell)$-curtains.  So that we can apply the top-down colouring strategy from \cref{Lungs-curtain} we show, in \cref{Heart2}, that any precolouring of at most $s$ vertices in the top of the root torso of a $(s-3,\ell)$-curtain $\torso{J}{B_x}$ of $J$ can be properly extended to an $(s+1)$-colouring of the lower torso $\ltorso{J}{B_x}$.  This is the most intricate part of the proof, which uses all of the machinery developed thus far:
\begin{enumerate}
    \item Like the proof of \cref{apex1} (for $s\ge 4)$, the first step works with the raised curtain $\ltorso{J}{B_x}_\uparrow$.  Although $\ltorso{J}{B_x}_\uparrow$ has bounded treewidth, we cannot immediately apply the techniques in \cref{BoundedTreewidth}.  In particular \cref{lem:LW2} has requirements that limit the largest value $k'$ such that $K_{2,k'}$ and $K_{3,k'}$ subgraphs appear in (certain subgraphs of) $\ltorso{J}{B_x}_\uparrow$.  The restrictions on $K_{3,k'}$ subgraphs only apply to the embedded parts of $\ltorso{J}{B_x}$, so these are satisfied by $\ltorso{J}{B_x}_\uparrow$ by Euler's Formula.
    \item In order to satisfy the restrictions on $K_{2,k'}$ subgraphs in $\ltorso{J}{B_x}$ we work with an $F$-compression $\ltorso{J}{B_x}^-_{\uparrow\bullet}$ of $\ltorso{J}{B_x}_\uparrow$.  Using \cref{cover_discs_ii} we show that $\ltorso{J}{B_x}^-_{\uparrow\bullet}$ does not contain $K_{2,k'}$ subgraphs for arbitrarily large values of $k'$.
    \item We add additional vertices $\alpha_1,\ldots,\alpha_r$ to  $\ltorso{J}{B_x}^-_{\uparrow\bullet}$ and colour the resulting graph using \cref{lem:LW2}.  We then extend this colouring in the natural way to the $F$-contraction $\ltorso{J}{B_x}_{\uparrow\bullet}$.  The addition of $\alpha_1,\ldots,\alpha_r$ ensures that each connected component of $\ltorso{J}{B_x}_{\uparrow\bullet}-B_x$ avoids all the colours that were assigned to (at most $s-3$) major apex vertices that dominate that component and leaves one colour available that can be used to break long skinny monochromatic components using \cref{component_breaking}.
    \item We then undo all of the contractions in $\ltorso{J}{B_x}_\uparrow$ used to create $\ltorso{J}{B_x}_{\uparrow\bullet}$.  Each vertex of $\ltorso{J}{B_x}_\uparrow$ inherits the colour of the vertex in $\ltorso{J}{B_x}_{\uparrow\bullet}$ into which it was contracted. As in the proof of \cref{apex1} (for $s\ge 4$) this leads to arbitrarily large monochromatic components.  Like the proof of \cref{apex1}, these components are long (they span many layers) but, unlike the proof of \cref{apex1}, these components are not necessarily skinny, so we cannot immediately apply \cref{component_breaking} to break long skinny monochromatic components.
    \item To resolve the issue with monochromatic components that are not skinny, we remove a sequence of bounded-size islands $I_1,\ldots,I_m$, each of which is contained in a single part of the layered partition $(\PP_y,\LL_y)$ of a single drape $\ltorso{\ltorso{J}{B_x}}{C_y}$ of the curtain $\ltorso{J}{B_x}$.  Each island $I_i$ is a $3$-island in the lower layers of the embedded part of $\ltorso{\ltorso{J}{B_x}}{C_y}$.  Since $\ltorso{\ltorso{J}{B_x}}{C_y}$ has at most $s-3$ major apex vertices, each island $I_i$ is an $s$-island in $\ltorso{\ltorso{J}{B_x}}{C_y}$.
    \item After removing $I_1,\ldots,I_m$ (letting $I:=\bigcup_{i=1}^m I_i$), we can then use \cref{skinny_components} to show that each monochromatic component of $\ltorso{J}{B_x}-I$ is skinny.  This allows us to apply \cref{component_breaking} to break long skinny monochromatic components of $\ltorso{J}{B_x}-I$ into components of bounded size.
    \item To complete colouring, we reintroduce the $3$-islands $I_m,\ldots,I_1$, colouring them so that each vertex in $I_i$ avoids the colours of its at most $3$ neighbours in $\ltorso{J}{B_x}-(\bigcup_{j=i+1}^m I_j)$ as well as the colours of the at most $s-3$ major apex vertices of the drape $\ltorso{\ltorso{J}{B_x}}{C_y}$ that contains $I_i$.
\end{enumerate}
With this proof overview complete, we now proceed.  The following lemma shows that we can restrict our attention to trees of $\PP$-connected upward-connected $(s-3,\ell)$-curtains.

\begin{lem}
\label{JstCurtain}
For integers $s\geq 3$ and $t\geq 2$, every $\JJ_{s,t}$-minor-free graph is a lower-minor-closed tree of $\PP$-connected upward-connected $(s-3,\ell)$-curtains for some $\ell=\ell(s,t)$.
\end{lem}

\begin{proof}
Let $K_s \oplus P_t$ be the complete join of $K_s$ and a $t$-vertex path, which is an element of $\JJ_{s,t}$. Observe that $K_2 \oplus P_t$ is planar, implying $K_s \oplus P_t$ is $(s-2)$-apex.  By definition, every $\JJ_{s,t}$-minor-free graph $J$ is $(K_s \oplus P_t)$-minor-free.
Therefore, by \cref{ApexMinorFreeStructure}, $J$ is a lower-minor-closed tree of $\PP$-connected upward-connected $(s-3,\ell)$-curtains for some $\ell=\ell(s,t)$.
\end{proof}

At this point the names of graphs under consideration are becoming quite long, so we will use overlines instead of subscript $0$ to denote the surface-embedded part of an almost-embedded graph. There should not be any confusion with the usual notation for the complement of a graph, a notion which is not used anywhere in the present paper. The contractions performed in $\torso{\torso{J}{B_x}}{C_y}$ to obtain $\torso{\torso{J}{B_x}}{C_y}_{\uparrow}$ involve edges whose endpoints are all in the lower layers of $\torso{\torso{J}{B_x}}{C_y}$.  Thus, these edges are all in the surface-embedded part $\overline{\torso{\torso{{J}}{B_x}}{C_y}}$ of the $(s-3,\ell)$-almost-embedded graph $\torso{\torso{J}{B_x}}{C_y}$. Every edge we contract while raising the drape $\torso{\torso{J}{B_x}}{C_y}$ is also present in the surface-embedded part $\overline{\ltorso{\ltorso{{J}}{B_x}}{C_y}}$ of the lower torso $\ltorso{\ltorso{J}{B_x}}{C_y}$. We will write $\overline{\torso{\torso{{J}}{B_x}}{C_y}}_{\uparrow}$ and $\overline{\ltorso{\ltorso{{J}}{B_x}}{C_y}}_\uparrow$ to denote the surface-embedded graphs obtained from $\overline{\torso{\torso{{J}}{B_x}}{C_y}}$ and $\overline{\ltorso{\ltorso{{J}}{B_x}}{C_y}}$, respectively, by performing the same contractions used to obtain  $\torso{\torso{J}{B_x}}{C_y}_{\uparrow}$ from $\torso{\torso{J}{B_x}}{C_y}$.

\begin{lem}\label{no_more_k2t}
  Let $J$ be a tree of $\PP$-connected upward-connected $(s-3,\ell)$-curtains described by the tree-decomposition $(B_x:x\in V(T))$ and, for each $x\in V(T)$, let $(C_y:y\in V(T_x))$ be a tree-decomposition that describes the $(s-3,\ell)$-curtain $\torso{J}{B_x}$ and for each $y\in V(T_x)$, let $(\PP_y,\LL_y)$ be the $\ell$-layered partition that describes the $(s-3,\ell)$-drape $\torso{\torso{J}{B_x}}{C_y}$ and let $(L_1,L_2,\ldots):=\LL_y$.

  Let $G:=\torso{\torso{J}{B_x}}{C_y}$ for some $x\in V(T)$ and some $y\in V(T_x)$ and let $A$ be the major apex set of $G$. Let $F:=V(\overline{G})\cap L_{\le 4}$ and let $\overline{G}^{\,\,-}_{\uparrow\bullet}$ be an $F$-compression of $\overline{G}_{\uparrow}$. Let $G_{\uparrow\bullet}^-$ be the graph obtained from $G_\uparrow$ by performing the same sequence of edge contractions and vertex deletions used to obtain $\overline{G}^{\,\,-}_{\uparrow\bullet}$ from $\overline{G}_{\uparrow}$.  Then $G_{\uparrow\bullet}^{-}- (A\cup L_{\le 4})$ is $K_{2,t'}$-subgraph-free for some $t':=t'(\ell)$.
\end{lem} 

\begin{proof}
  The lemma would immediately follow from \cref{cover_discs_ii} if  the number of components of $\overline{G}[F]$ were bounded by some $t'(\ell)$. This may not be the case, so an extra step is required.

  Recall that $\overline{G}[L_1]$ contains exactly the vertices of $\overline{G}$ that participate in the vortices $G_1,\ldots,G_\ell$ of the almost-embedded graph $\overline{G}$. For each $i\in\{1,\ldots,\ell\}$, $G_i$ is defined in terms of a $\overline{G}$-clean topological disc $D_i$ that has the vertices of $\overline{G}\cap G_i$ on its boundary. Define the embedded graph $Q_\uparrow$ by starting with $\overline{G}_\uparrow$ and adding a vertex $v_i$ in the interior of $D_i$ that is adjacent to each vertex of $G_i\cap \overline{G}$, for each $i\in\{1,\ldots,\ell\}$.  Define the graph $Q^{-}_{\uparrow\bullet}$ analogously, starting from the graph $\overline{G}^{\,\,-}_{\uparrow\bullet}$.  Let $L_0:=\{v_1,\ldots,v_\ell\}$.  Then $Q_{\uparrow}[L_0\cup L_1]$ has at most $\ell$ components and, since $\LL_y$ is upward-connected, \cref{UpwardConnected} implies that $Q_{\uparrow}[L_0\cup F]$ has at most $\ell$ components.  Furthermore, $Q^{-}_{\uparrow\bullet}$ is an $(L_0\cup F)$-compression of $Q_{\uparrow}$.  By \cref{cover_discs_ii}, $Q_{\uparrow\bullet}^{-}-L_{\le 4}$ is $K_{2,t'}$-subgraph-free for some $t':=t'(\ell)$.  Since $G^{-}_{\uparrow\bullet}-(A\cup L_{\le 4})$ is a subgraph of $Q^{-}_{\uparrow\bullet}-L_{\le 4}$, this completes the proof.
\end{proof}

Before the main event, we need a simple lemma about treewidth.

\begin{lem}\label{lem:apexify}
Let $G$ be a graph of treewidth at most $k$, and let $C_1,\ldots,C_r$ be pairwise disjoint connected subgraphs in $G$. Let $G^+$ be the graph obtained from $G$ by adding, for each $i\in\{1,\dots,r\}$, a new vertex $v_i$ whose neighbourhood is a subset of $V(C_i)$. Then $G^+$ has treewidth at most $2k+1$.
\end{lem}

\begin{proof}
Consider a tree-decomposition  $(B_x:x\in V(T))$  of $G$, such that $|B_x|\le k+1$ for each node $x\in V(T)$. Let $B_x':=B_x\cup\{v_i: V(C_i)\cap B_x\neq\emptyset\}$ for each $x\in V(T)$. Since $C_i$ is connected, the set $\{x\in V(T): B_x\cap C_i\ne \emptyset\}$ induces a subtree of $T$. Thus $(B_x' : x\in V(T))$ is a tree-decomposition of $G^+$. If a new vertex $v_i$ is in $B_x'$, then $B_x$ contains a vertex of $C_i$. Since $C_1,\dots,C_r$ are pairwise disjoint, $|B_x'|\le 2|B_x|\le 2k+2$ for each $x\in V(T)$. This shows that $G^+$ has treewidth at most $2k+1$, as desired.
\end{proof}

We are now ready to prove the following lemma, which is the last big step in the proof of \cref{Jst}.

\begin{lem}\label{Heart2}
  There is a function $f$ with the following properties. Let $s\geq 3$, $t\geq s+2$ and $\ell\ge 2$ be integers.  Let $J$ be a $\JJ_{s,t}$-minor-free lower-minor-closed tree of $\PP$-connected upward-connected $(s-3,\ell)$-curtains described by $\TT_0:=(B_x:x\in V(T_0))$.    Let $S$ be a set of at most $s$ vertices contained in the top of the root torso of $\torso{J}{B_x}$, for some $x \in V(T_0)$.  Then for any $(s+1)$-colouring of $S$, there is an $(s+1)$-colouring of the lower torso $\ltorso{J}{B_x}$ of $J$ that properly extends the given precolouring of $S$ and has clustering at most $f(s,t,\ell)$.
\end{lem}

\begin{proof}
  Let $\mathcal{T}:=(C_y:y\in V(T))$  be the tree-decomposition that describes the curtain $\mathdefn{G}:=\torso{J}{B_x}$.  For each $y\in V(T)$, let $\mathcal{E}_y:=(A_y,\hat{A}_y,\overline{G^y},G^y_1,\ldots,G^y_{\ell})$ be the $(s-3,\ell)$-almost-embedding of the $(s-3,\ell)$-drape $G^y:=\torso{\torso{J}{B_x}}{C_y}$ and let $(\PP_y,\LL_y:=(L^y_1,L^y_2,\ldots))$ be the $\ell$-layered partition that  describes $\torso{\torso{J}{B_x}}{C_y}$.

  Let $H:=\ltorso{J}{B_x}$ and, for each $y\in V(T)$, let $H^y:=\ltorso{\ltorso{J}{B_x}}{C_y}$.
  For each $y\in V(T)$, $H^y$ is a spanning subgraph of $G^y$. Since $H$ is a spanning subgraph of $G$, $\mathcal{T}$ is a tree-decomposition of $H$. For each $y\in V(T)$, $\mathcal{E}_y$ defines an $(s-3,\ell)$-almost-embedding $(A_y,\hat{A}_y,\overline{H^y},H^y_1,\ldots,H^y_\ell)$ of $H^y$. Raise the $(s-3,\ell)$-curtain $G$ to obtain the raised $(s-3,\ell)$-curtain $G_{\uparrow}$.  For each $y\in V(T)$, this defines the raised drape $G^y_{\uparrow}$ and the layering $\LL^y_{\uparrow}:=(L^y_1,\ldots,L^y_5,L^y_\uparrow)$ of $G^y-A_y$.
  Define the raised curtain $\mathdefn{H_\uparrow}:=\ltorso{J}{B_x}_\uparrow$ and each raised drape $\mathdefn{H^y_\uparrow}:=\ltorso{\ltorso{J}{B_x}}{C_y}_\uparrow$ in (what is by now)  the obvious way.  By \cref{raised_minor}, $H_\uparrow$ is a minor of $J$ and is therefore $\JJ_{s,t}$-minor-free.

  For each $y\in V(T)$, let $F_y:=V(\overline{H^y})\cap L^y_{\le 4}$,  let \defn{$\overline{H^y}_\uparrow$} be the surface-embedded graph obtained by applying the contractions used to obtain $H^y_\uparrow$ from $H^y$ to the graph $\overline{H^y}$, and let $\overline{H^y}_{\uparrow\bullet}$ and $\overline{H^{y}}^{\,\,-}_{\uparrow\bullet}$ be the $F_y$-contraction and $F_y$-compression, respectively, of $\overline{H^y}_\uparrow$.  Recall that $F_y$-contraction and $F_y$-compression are defined in terms of $F_y$-avoiding discs meaning that the boundary of each disc contains no vertex of $F_y$ and the interior of each disc contains no vertex of $N_{\overline{H^y}_\uparrow}(F_y)$.  In the current setting, this means that the interior of each disc contains only vertices from $L^y_\uparrow$.  For each $y\in V(T)$, applying the same operations used to obtain $\overline{H^y}_{\uparrow\bullet}$ and $\overline{H^{y}}^{\,\,-}_{\uparrow\bullet}$ to $H^y_\uparrow$ yields graphs $H^y_{\uparrow\bullet}$ and $H^{y-}_{\uparrow\bullet}$. We also define $G^{y-}_{\uparrow\bullet}$ as the graph obtained by applying the $F_y$-compression operations to $G^y_\uparrow$.

  Let $\mathdefn{H_{\uparrow\bullet}}:=\bigcup_{y\in V(T)} H^y_{\uparrow\bullet}$ and let $\mathdefn{H^{-}_{\uparrow\bullet}}:=\bigcup_{y\in V(T)} H^{y-}_{\uparrow\bullet}$. In words, both these graphs are obtained from the curtain $\ltorso{J}{B_x}$ by simultaneously raising all the drapes and then performing compressions and contractions. Each of the graphs $H^{-}_{\uparrow\bullet}$ and $H_{\uparrow\bullet}$ is a minor of $H_\uparrow$ and is therefore $\JJ_{s,t}$-minor-free. Let $A:=\bigcup_{y\in V(T)}A_y$ be the union of the major apex sets of the drapes of $H$. For each $i\in \mathbb{N}$, define $L_i:=\bigcup_{y\in V(T)} L^y_i$.  Note that the sets $L_1$ and $L_2$ are not necessarily disjoint, but $L_i$ and $L_j$ are disjoint for any $j\ge 3$ and $i\neq j$.  Define $\mathdefn{L_{\uparrow\bullet}}:=V(H_{\uparrow\bullet}-(A\cup L_{\le 5}))$ and $\mathdefn{L^-_{\uparrow\bullet}}:=V(H^-_{\uparrow\bullet}-(A\cup L_{\le 5}))=L_{\uparrow\bullet}\cap V(H^-_{\uparrow\bullet})$.  Note that each vertex in $L_{\uparrow\bullet}$ is obtained by contracting a connected subgraph of $H[L_{\ge 6}]$.

  By the definition of lower torso, $H[L_{\ge 3}]=G[L_{\ge 3}]$. Since $H^{y-}_{\uparrow\bullet}$ and $G^{y-}_{\uparrow\bullet}$ are obtained by performing identical operations on $H^y_\uparrow$ and $G^y_\uparrow$, this implies that $H^{y-}_{\uparrow\bullet}[\bigcup_{i=3}^5 L_i\cup L_{\uparrow\bullet}] = G^{y-}_{\uparrow\bullet}[\bigcup_{i=3}^5 L_i\cup L_{\uparrow\bullet}]$ for each $y\in V(T)$.  We will make use of this to apply \cref{no_more_k2t} (which is a statement about torsos) to lower torsos.

  Now we prepare to apply \cref{lem:LW2}, which requires the addition of some extra vertices.  Let $R_1,\ldots,R_r$ be the vertex sets of the connected components of $H_{\uparrow\bullet}[N_{H_{\uparrow\bullet}-A}[L_{\uparrow\bullet}]]$. For each $i\in\{1,\ldots,r\}$, let $y_i$ be the unique node of $T$ such that $R_i\subseteq V(H^{y_i}_{\uparrow\bullet})$.  Note that $y_1,\ldots,y_r$ are not necessarily distinct.
  Starting with $H_{\uparrow\bullet}$, define the graph $H_{\uparrow\bullet}^+$  by introducing vertices $\alpha_1,\ldots,\alpha_r$, where $\alpha_i$ is adjacent to every vertex in $R_i\cap L_{\uparrow\bullet}$, for each $i\in\{1,\ldots,r\}$.  Note that the sets $R_1,\ldots,R_r$ contain vertices in $L_5$ but $\alpha_1,\ldots,\alpha_r$ are only adjacent to vertices in $L_{\uparrow\bullet}$.

  \begin{figure}[htbp]
  \begin{center}
    \includegraphics[scale=1.2]{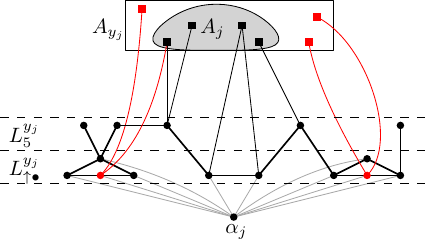}
  \end{center}
  \caption{The vertices of $S_j$ are depicted with circles. The vertices of $S_j^-$ are depicted with black circles, and the vertices of $S_j\setminus S_j^-\subseteq V(H_{\uparrow\bullet}^{+})\setminus V(H_{\uparrow\bullet}^{+-})$ are depicted with red circles.}
  \label{sjmoins}
\end{figure}

  Let $H_{\uparrow\bullet}^{+-}:= H_{\uparrow\bullet}^{+}[V(H^-_{\uparrow\bullet})\cup\{\alpha_1,\ldots,\alpha_r\}]$.  In words, $H_{\uparrow\bullet}^{+-}$ is obtained from $H_{\uparrow\bullet}^{+}$ by deleting the vertices that are created in some $F_y$-contraction, but are not present in the corresponding $F_y$-compression. For each $j\in\{1,\ldots,r\}$, let $S_j:= R_j\cup\{\alpha_j\}$ and let $S^-_j:= S_j\cap V(H^{+-}_{\uparrow\bullet})$.  We then define $A_j := A_{y_j} \cap N_{H^{+-}_{\uparrow\bullet}}(S^-_j)$. See Figure \ref{sjmoins} for an illustration. 
  Observe that $S^-_j\subseteq S_j\subseteq L^{y_j}_{5}\cup L^{y_j}_{\uparrow\bullet}\cup\{\alpha_j\}$, and that 
  $S^-_j=N_{H^{+-}_{\uparrow\bullet}-A_j}^2[\alpha_j]$, which agrees with  property \cref{no_a_in_s} of admissible set for the graph
  $H^{+-}_{\uparrow\bullet}$.

  The set $\mathcal{S}:=\{(\alpha_1,A_1),\ldots,(\alpha_r,A_r)\}$ is $(s-3)$-admissible with respect to $H_{\uparrow\bullet}^{+-}$ for the following reasons:
  \begin{compactenum}
    \item[\cref{a_size}:] $|A_j|\le s-3$, by the definition of ($s-3,\ell)$-curtain.

    \item[\cref{no_a_in_X}:] $\alpha_j$ is a vertex of $H_{\uparrow\bullet}^{+-}$ and $N_{H_{\uparrow\bullet}^{+-}}[\alpha_j]\cap A_j\subseteq (L_{\uparrow\bullet}\cup\{\alpha_j\})\cap A_j=\emptyset$, for each $j\in\{1,\ldots,r\}$.

    \item[\cref{no_a_in_s}:] We need to show that $\left(\bigcup_{i=1}^r A_i\right) \cap \left(\bigcup_{j=1}^r S^-_j\right)=\emptyset$. Suppose otherwise, so that $v\in A_i\cap S^-_{j}$ for some $i,j\in\{1,\ldots,r\}$.  Then $i\neq j$ and $y_i\neq y_{j}$ since $A_i\subseteq A_{y_i}$ does not intersect $L^{y_i}_5\cup L^{y_i}_{\uparrow\bullet}\cup\{\alpha_i\}\supseteq S^-_i$. Therefore, by the definition of curtain, $v\in A_{y_i}\cap (L^{y_{j}}_1\cup L^{y_{j}}_2 \cup A_{y_{j}})$, which contradicts the assumption that $v\in S^-_{j}\subseteq L^{y_{j}}_5 \cup L^{y_{j}}_{\uparrow\bullet}$.

    \item[\cref{a_in_n_s}:] $A_j \subseteq N_{H^{+-}_{\uparrow\bullet}}(S^-_j)$ by the definition of $A_j$.

    \item[\cref{s_disjoint}:] $S^{-}_1,\ldots,S^{-}_r$ are pairwise disjoint because $R_1,\ldots,R_r$ are pairwise disjoint.
  \end{compactenum}

  We want to apply \cref{lem:LW2} to $H^{+-}_{\uparrow\bullet}$ with $\mathcal{S}$ as our $(s-3)$-admissible set.  We first argue that $H^{+-}_{\uparrow\bullet}$ has bounded treewidth.  By \cref{raised_treewidth}, there is a function $c_1$ such that $H_\uparrow$ has treewidth at most $c_1(s,\ell)$. Let $c_1:=c_1(s,\ell)$.  Since $H_{\uparrow\bullet}$ is obtained from $H_\uparrow$ by contracting connected subgraphs, $H_{\uparrow\bullet}$ also has treewidth at most $c_1$.  Since $R_1,\ldots,R_r$ are pairwise disjoint and $H_{\uparrow\bullet}[R_j]$ is connected for each $j\in\{1,\ldots,r\}$,   \cref{lem:apexify} shows that the addition of $\alpha_1,\ldots,\alpha_r$ to $H_{\uparrow\bullet}$ creates a graph $H_{\uparrow\bullet}^+$ with treewidth at most $2c_1+1$. The removal of vertices in $S_j\setminus S^-_j$ to obtain $H^{+-}_{\uparrow\bullet}$  does not increase treewidth.

  In addition to having bounded treewidth, \cref{lem:LW2} requires that $H_{\uparrow\bullet}^{+-}$ satisfy several additional conditions, which we now verify.
  \begin{compactenum}
    \item[\cref{kh_minor_free}:] $H^{+-}_{\uparrow\bullet}-\{\alpha_1,\ldots,\alpha_r\}$ is $\JJ_{s,t}$-minor-free because $H^{+-}_{\uparrow\bullet}-\{\alpha_1,\ldots,\alpha_r\}=H^{-}_{\uparrow\bullet}$ (by definition) and $H^{-}_{\uparrow\bullet}$ is a minor of $H_\uparrow$ which, by \cref{raised_minor} is a minor of $J$, which is $\JJ_{s,t}$-minor-free.

    \item[\cref{sj_connected}:] We must show that $H^{+-}_{\uparrow\bullet}[S^-_j\setminus\{\alpha_j\}]$ is connected, for each $j\in\{1,\ldots,r\}$.  By definition, $H_{\uparrow\bullet}[R_j]=H^{+}_{\uparrow\bullet}[R_j]=H^{+}_{\uparrow\bullet}[S_j\setminus\{\alpha_j\}]$ is connected.  By the definition of $F_{y_j}$-compression, for each vertex $v\in S_j\setminus S^-_j$, there exists $w\in S^-_j\cap L_{\uparrow\bullet}\subseteq R_j$ that dominates $v$ in $H_{\uparrow\bullet}$.  Therefore, $H^{+-}_{\uparrow\bullet}[S^-_j]-\alpha_j$ is obtained from $H_{\uparrow\bullet}[R_j]$ by repeatedly removing a dominated vertex. It follows that $H^{+-}_{\uparrow\bullet}[S^-_j\setminus\{\alpha_j\}]$ is connected, since removing a dominated vertex from a connected graph does not disconnect the graph.
    \item[\cref{k2_t_free}:] We must show that $H^{+-}_{\uparrow\bullet}[S^-_j\setminus\{\alpha_j\}]$ is $K_{2,t'}$-subgraph-free, for each $j\in\{1,\ldots,r\}$ and some $t'=t'(\ell)$. By definition, $S^-_j\setminus\{\alpha_j\}\subseteq V(H^{y_j-}_{\uparrow\bullet})\cap (L^{y_j}_5\cup L_{\uparrow\bullet})\subseteq V(G^{y_j-}_{\uparrow\bullet}-(L^{y_j}_{\le 4}\cup A_{y_j}))$.  By \cref{no_more_k2t}, $G^{y_j-}_{\uparrow\bullet}-(L^{y_j}_{\le 4}\cup A_{y_j})$ is $K_{2,t'}$-subgraph-free (for some $t'=t'(\ell)$).
    Since $H^{+-}_{\uparrow\bullet}[S^-_j\setminus\{\alpha_j\}]$ is a subgraph of $G^{y_j-}_{\uparrow\bullet}-(L^{y_j}_{\le 4}\cup A_{y_j})$,  it is also $K_{2,t'}$-subgraph-free.
    \item[\cref{k3_t_free}:] For each $j\in\{1,\ldots,r\}$,
    $H^{+-}_{\uparrow\bullet}[N^4_{H^{+-}_{\uparrow\bullet}-A_j}[\alpha_j]\setminus\{\alpha_j\}] \subseteq H^{y_j-}_{\uparrow\bullet}[L_{\uparrow\bullet}\cup L^{y_j}_5\cup L^{y_j}_4\cup L^{y_j}_3]\subseteq \overline{H^{y_j}}^{\,\,-}_{\uparrow\bullet}$ has an embedding on a surface of Euler genus at most $\ell$.  Therefore $H^{+-}_{\uparrow\bullet}[N^4_{H^{+-}_{\uparrow\bullet}-A_j}[\alpha_j]\setminus\{\alpha_j\}]$ is $K_{3,k'}$-subgraph-free (for some $k':=k'(\ell)$) by Euler's Formula.
  \end{compactenum}

  Define the list-assignment $L_0:V(H^{+-}_{\uparrow\bullet})\to2^{\{1,\ldots,s+1\}}$ for $H^{+-}_{\uparrow\bullet}$ as follows:  For each $v\in S$, let $L_0(v)$ be the singleton set that contains the colour assigned to $v$ by the given precolouring of $S$.  For each $w\not\in S$, $L_0(w):=\{1,\ldots,s+1\}\setminus \bigcup\{L_0(v):v\in N_{H^{+-}_{\uparrow\bullet}}(w)\cap S \}$. Since $|S|\le s$, this ensures that each vertex has a non-empty list $L_0$ and that any $L_0$-colouring of $H^{+-}_{\uparrow\bullet}$ properly extends the precolouring of $S$.

  Let $P:=P(L_0)$ be the set of vertices precoloured by $L_0$ and let $Q:=Q(L_0)$ be the set of vertices with at least one and not more than $s-1$ precoloured neighbours (see the beginning of \cref{BoundedTreewidth} for the definitions of $P$, $Q$, and $(s,p)$-good). We now define a second list-assignment $L$. The need for $L_0$ and $L$ comes from the existence of vertices with fewer than $s$ neighbours in $S$ but at least $s$ neighbours in $P$; the number of such vertices is unbounded so we do not want to precolour them, but we require them to avoid colours used by their neighbours in $S$. For each $v\in Q$, let
  $L(v):=\{1,\ldots,s+1\}\setminus \bigcup\{L_0(w):w\in N_{H^{+-}_{\uparrow\bullet}}(v)\cap P\}$.  For each $v\not\in Q$, let $L(v):=L_0(v)$.  By the definition of $Q$, each vertex $v$ of $Q$ has at most $s-1$ neighbours in $P$, and thus satisfies $|L(v)|\ge 2$. It follows that $P=P(L)$ and therefore $Q=Q(L)$.  The list-assignment $L$ is a specialization of $L_0$, and since any $L_0$-colouring of $H^{+-}_{\uparrow\bullet}$ properly extends the precolouring of $S$, any $L$-colouring of $H^{+-}_{\uparrow\bullet}$ properly extends the precolouring of $S$.

  We claim that $L$ is trivially compatible with $\mathcal{S}$ because $P$ does not contain any vertex of any trigger set $S_j$.  Indeed, $P$ consists of $S$ along with vertices of $H^{+-}_{\uparrow\bullet}$ that have $s$ neighbours in $S$. Since $S$ is in the top of the root torso of $G$, $S\subseteq A\cup L_1$.  Since $S^-_j\subseteq L^{y_j}_5\cup L^{y_j}_{\uparrow\bullet}\cup\{\alpha_j\}$ and $\alpha_j\notin A\cup L_1$, we have $S^-_j\cap S=\emptyset$ for each $j\in\{1,\ldots,r\}$.  All that remains is to check that no vertex in $S^-_j$ has $s$ neighbours in $S$.  Let $v\in S^-_j$. If $v=\alpha_j$ then $N_{H^{+-}_{\uparrow\bullet}}(v)=N_{H^{+-}_{\uparrow\bullet}}(\alpha_j)\subseteq L^{y_j}_{\uparrow\bullet}$, so $|N_{H^{+-}_{\uparrow\bullet}}(v)\cap S|=0$ and $v\notin P$. If $v\in S^-_j\setminus\{\alpha_j\}$ then $v\in L^{y_j}_5\cup L^{y_j}_{\uparrow\bullet}$, so $N_{H^{+-}_{\uparrow\bullet}}(v)\subseteq L^{y_j}_4\cup L^{y_j}_5\cup L^{y_j}_{\uparrow\bullet}\cup A_j$ and therefore $|N_{H^{+-}_{\uparrow\bullet}}(v)\cap S|\le |A_j|\le s-3$, so $v\notin P$.  Thus $S^-_j\cap P=\emptyset$ for each $j\in\{1,\ldots,r\}$.

  Now we check that $L$ satisfies conditions \cref{q_size}--\cref{q_proper} for being $(s,p)$-good.
  \begin{compactenum}[(g1)]
    \item[\cref{q_size}:] $|L(v)|\ge s+1-|N_{H^{+-}_{\uparrow\bullet}}(v)\cap P|$, by definition, for all $v\in Q(L)$;
    \item[\cref{u_size}:] $|L(v)|\ge 2$ for all $v\in N_{H^{+-}_{\uparrow\bullet}}(P)\setminus Q$ since each such vertex has at most $s-1$ neighbours in $S$ and therefore has a list of size at least $(s+1)-(s-1)=2$;
    \item[\cref{other_size}:] $|L(v)|=s+1$ for all $v\in V(H^{+-}_{\uparrow\bullet})\setminus N_{H^{+-}_{\uparrow\bullet}}[P]$; and
    \item[\cref{q_proper}:] $L(v)\cap L(u)=\emptyset$ for all $v\in Q$ and $u\in N_{H^{+-}_{\uparrow\bullet}}(v)\cap P$ by definition.
  \end{compactenum}

  We did not establish \cref{p_size}---the upper bound on the size of $P$---because the size of $P$ cannot be upper bounded by any function of $s$, $t$, and $\ell$.  However, each vertex in $P\setminus S$ is adjacent to $s$ vertices in $S$. Thus, if we let $Y:=P\setminus S$ and $q:=|Y|$, then either $q=0$, or $|S|=s$ and $H^{+-}_{\uparrow\bullet}[P]$ contains a $K_{s,q}$-subgraph with parts $S$ and $Y$.  By \cref{no_surface_in_x} (whose proof relies only on \cref{k2_t_free,k3_t_free} and the admissibility of $\mathcal{S}$) and \cref{less_than_t_per} (whose proof relies only on \cref{kh_minor_free} and the admissibility of $\mathcal{S}$), this implies that $q < t'$ (where $t':=t'(s,t,\ell)$ is defined as in \cref{lem:LW2}) or that each component of $H^{+-}_{\uparrow\bullet}-S$ contains at most $t-1$ vertices of $Y$. In the former case, $|P| < t'+s$ and $L$ satisfies \cref{p_size} for $p\ge t'+s$, so we can apply \cref{lem:LW2} to all of $H^{+-}_{\uparrow\bullet}$ to obtain an $L$-colouring with clustering at most $c_2(c_1,s,t,\ell)$, for some function $c_2$. Let $c_2:= c_2(c_1,s,t,\ell)$. In the latter case, we can apply \cref{lem:LW2} independently on $H_C:=H^{+-}_{\uparrow\bullet}[S\cup V(C)]$ for each component $C$ of $H^{+-}_{\uparrow\bullet}-S$.  In each application, the size of the precoloured set is at most $s+t-1$ (so $L$ is $(s,s+t-1)$-good) and we obtain an $L$-colouring that properly extends the precolouring of $S$.\footnote{The colouring of each graph $H_C$ is similar to the colouring of the graph $G_C$ in Case~A of the proof of \cref{lem:LW2}.}  Therefore, each monochromatic component is either contained in $S$ or is contained in a single component of $H^{+-}_{\uparrow\bullet}-S$.  Therefore, all of these colourings can be combined to provide an $L$-colouring of $H^{+-}_{\uparrow\bullet}$ with clustering at most $c_3(c_1,s,t,\ell)$ for some function $c_3$. Let $c_3:=c_3(c_1,s,t,\ell)$.

  In either case, we obtain an $L$-colouring $\varphi$ of $H^{+-}_{\uparrow\bullet}$ with clustering at most $c_4:=\max\{c_2,c_3\}$, where $\varphi(a)\neq\varphi(\alpha_j)$ for each $j\in\{1,\ldots,r\}$ and each $a\in A_j$, and where $\varphi(u)\neq\varphi(a)$ for each $j\in\{1,\ldots,r\}$, each $u\in N_{H^{+-}_{\uparrow\bullet}}(\alpha_j)$ and each $a\in N_{H^{+-}_{\uparrow\bullet}}(u)\cap (A_j\cup\{\alpha_j\})$.  Since $H^{-}_{\uparrow\bullet}$ is a subgraph of $H^{+-}_{\uparrow\bullet}$, $\varphi$ is also a colouring of $H^{-}_{\uparrow\bullet}$.  By the definition of $L$, $\varphi$ is an $(s+1)$-colouring of $H^{-}_{\uparrow\bullet}$ that properly extends the precolouring of $S$. We now extend $\varphi$ to obtain an $(s+1)$-colouring of $H$ which properly extends the precolouring of $S$. We do this in several steps.

  For each $y\in V(T)$, each $v\in L^y_{\uparrow}$ is obtained by contracting a connected subgraph $X_v$ of $H^y[L^y_{\ge 6}]$. Each $v\in L^y_{\uparrow\bullet}$ is obtained by contracting a connected subgraph $X'_v$ of $H^y_\uparrow[L^y_{\uparrow}]$. Thus, each $v\in L^y_{\uparrow\bullet}$ corresponds to a connected subgraph $X_v = \bigcup_{u \in V(X'_v)} X_u$ of $H^y[L^y_{\ge 6}]$. For each $v\in L^y_{\uparrow}$, $V(X_v)$ is contained in a single part of the partition $\PP_y$, so $V(X_v)$ is $\ell$-skinny with respect to the layering $\LL_y$ of $H^y-A_y$. However, for $v\in L^y_{\uparrow\bullet}\setminus L^y_{\uparrow}$, it is not necessarily the case that $V(X_v)$ is $\ell$-skinny with respect to $\LL_y$. Our first step is to deal with this issue, by removing a sequence of $s$-islands from $H$, so that we are in a position to apply \cref{skinny_components}, followed by \cref{component_breaking}.   
  
  Let $c_5(\ell)$ be a function that is discussed further below.  We greedily remove a sequence of vertex sets $I_1,\ldots,I_m$ where $I_i$ is a $(3,c_5(\ell))$-island of $\overline{H^y}-(\bigcup_{j=1}^{i-1}I_j)$ and $I_i\subseteq V(X_v)$ for some $y\in V(T)$ and some $v \in L^y_{\uparrow\bullet}$. This process ends only when no such $I_i$ exists.  Let $I_1,\ldots,I_m$ be the sequence of islands removed in this way and let $I:=\bigcup_{i=1}^m I_i$.  For each $i\in\{1,\ldots,m\}$, $I_i\subseteq L^{y_j}_{\ge 6}$ for some $j\in\{1,\ldots,r\}$.  For each $v \in I_i$, $N_H(v)\subseteq N_{\overline{H^y}}(v)\cup A_{y}$.  Therefore $I_i$ is an $(s,c_5(\ell))$-island of $H-(\bigcup_{j=1}^{i-1}I_j)$.
  
  Recall that $A=\bigcup_{y\in V(T)}A_y$ is the union of the major apex sets of the drapes of $H$. We argue that we are now in a position to apply \cref{component_breaking} to $H-(A\cup I)$.  We now show that there exists $c_5=c_5(\ell)$ such that $V(X_v)\setminus I$ is $c_5$-skinny with respect to $\LL$, for each $v\in V(H_{\uparrow\bullet})\setminus V(H)$.  Let $v\in V(H_{\uparrow\bullet})\setminus V(H)$.  If $v\in V(H_\uparrow)$, then $V(X_v)\subseteq P_v$ for some $P_v\in\mathcal{P}_y$ and $y\in V(T)$, so $V(X_v)$ is $\ell$-skinny with respect to $\mathcal{L}_y$.  Otherwise, $v\in V(H_{\uparrow\bullet})\setminus V(H_{\uparrow})$.  In this case, $X_v$ is obtained by contracting a connected subgraph $X'_v$ of $\overline{H^y}_{\uparrow}$ while performing the $F_y$-contraction of $\overline{H^y}_{\uparrow}$.  By the definition of $F_y$-contraction, $N_{\overline{H^y}_{\uparrow}}(V(X'_v))$ is a set of $3$ or $4$ vertices in $\overline{H^y}_\uparrow$.  Therefore $N_{\overline{H^y}}(V(X_v))$ is contained in $3$ or $4$ parts of $\mathcal{P}_y$.

  Consider the graph $R:=\overline{H^y}[N_{\overline{H^y}-I}(V(X_v)\setminus I)]$ and the subset $U:=V(R)\setminus V(X_v)$ of $V(R)$.  Since $U\subseteq N_{\overline{H^y}}(V(X_v))$, $U$ is contained in $3$ or $4$ parts of $\PP_y$.  Therefore $U$ is $4\ell$-skinny with respect to the layering of $R$ induced by $\LL_y$.  By \cref{skinny_components}, $V(R-U)$ contains a $(3,c_5(\ell))$-island of $R$ or $V(R-U)$ is $c_5(\ell)$-skinny with respect to $\LL_y$, for $c_5(\ell):=\max\{c(\ell+1), 180\ell\}$, where $c$ is the constant from \cref{three_islands}.  Let $c_5:=c_5(\ell)$.  If $V(R-U)$ contains a $(3,c_5)$-island $D$ of $R$ then, for each $v\in D$, $N_{\overline{H^y}-I}(v)\subseteq V(R)$, so $N_{\overline{H^y}-I}(v)=N_{R}(v)$.  Thus $D$ is a $(3,c_5)$-island of $\overline{H^y}-I$. This is not possible because $D\subseteq V(R-U)= V(X_v)\setminus I$ and $V(X_v)\setminus I$ has no $(3,c_5)$-island of $\overline{H^y}-I$, by the definition of $I=\bigcup_{i=1}^m I_i$.  Thus $V(X_v)\setminus I = V(R-U)$ is $c_5$-skinny with respect to $\LL$.

  We now do a sequence of steps:
  First we extend $\varphi$ from $H^-_{\uparrow\bullet}$ to $H_{\uparrow\bullet}$.  Then we use the colouring $\varphi$ of $H_{\uparrow\bullet}$ to obtain a colouring $\varphi$ of $H-I$ such that each $\varphi$-monochromatic component of $H-I$ is skinny with respect to the global layering $\LL:=(L_1\cup L_2,L_3,L_4,\ldots)$ of $H-(A\cup I)$.  Then we apply \cref{component_breaking} to obtain a colouring $\varphi'$ of $H-I$ with bounded clustering.  After recolouring a specific subset of vertices, we finally properly extend $\varphi'$ to $H$ by colouring the islands $I_m,\ldots,I_1$ in the opposite order from which they were removed.
    
  The extension from $H^-_{\uparrow\bullet}$ to $H_{\uparrow\bullet}$ is easy:  Let $Z:= V(H_{\uparrow\bullet})\setminus V(H^{-}_{\uparrow\bullet})$.  By the definition of $F$-compression, the vertices in $Z$ form an independent set in $H_{\uparrow\bullet}$.  Each $v\in Z$ belongs to $N_{H^{+}_{\uparrow\bullet}}(\alpha_j)\subset L^{y_j}_{\uparrow\bullet}$ for exactly one $j\in\{1,\ldots,r\}$.  By \cref{no_3_dominated}, $v$ has degree at most 2 in $\overline{H^{y_j}}_{\uparrow\bullet}$. Since $|A_{y_j}|\le s-3$, there exists a colour $\phi\in\{1,\ldots,s+1\}\setminus \{\varphi(x):x\in A_{y_j}\cup \{\alpha_j\}\cup N_{\overline{H^{y_j}}_{\uparrow\bullet}}(v)\}$.  We set $\varphi(v):=\phi$.  This properly extends the colouring $\varphi$ to a colouring of $H_{\uparrow\bullet}$, that still has clustering at most $c_4$, since each new monochromatic component has only a single vertex in $Z$. 

  This colouring also has the following critical property: For each $v\in L_{\uparrow\bullet}\setminus Z$, $\varphi(v)\notin\{\varphi(a):a\in (N_{H_{\uparrow\bullet}}(v)\cap A_j)\cup\{\alpha_j\}\}$ by Condition \cref{apex_proper2} of \cref{lem:LW2}. For such a vertex $v\in L_{\uparrow\bullet}\setminus Z$, $N_{H_{\uparrow\bullet}}(v)\cap A_{j}=N_{H_{\uparrow\bullet}}(v)\cap A_{y_j}$ by the definition of $A_j$, and thus $\varphi(v)\notin\{\varphi(a):a\in (N_{H_{\uparrow\bullet}}(v)\cap A_{y_j})\cup\{\alpha_j\}\}$. For each $v\in L_{\uparrow\bullet}\cap  Z$, it follows from the paragraph above that also $\varphi(v)\notin\{\varphi(a):a\in (N_{H_{\uparrow\bullet}}(v)\cap A_j)\cup\{\alpha_j\}\}$. 
  Thus, for each $v\in L_{\uparrow\bullet}$, $\varphi(v)\neq \varphi(a)$ for each $a\in N_H(v)\cap A$.

  To extend $\varphi$ from $H_{\uparrow\bullet}$ to $H-I$, we set $\varphi(w):=\varphi(v)$ for each $v\in L_{\uparrow\bullet}$ and each $w\in V(X_v)\setminus I$.  Since the $\varphi$-monochromatic components of $H_{\uparrow\bullet}$ have size at most $c_4$ and $V(X_v)\setminus I$ is $c_5$-skinny with respect to $\LL$, it follows that each $\varphi$-monochromatic component of $H-I$ is $c_4c_5$-skinny with respect to $\LL$.  Furthermore, the colouring $\varphi$ of $H-I$ inherits the following property from the colouring $H_{\uparrow\bullet}$:  For each $w\in L_{\ge 6}$, $\varphi(w)\neq\varphi(a)$ for each $a\in N_H(w)\cap A$.

  In order to apply \cref{component_breaking}, which requires a layering, it is helpful to leave $A$ aside for now and focus instead on $H-(A\cup I)$. We have an $(s+1)$-colouring of $H-(A\cup I)$ and a layering $\LL$ of $H-(A\cup I)$ such that each monochromatic component is $c_4c_5$-skinny with respect to the layering $\LL$ and, for each component $C$ of $H[L_{\ge 6}]$ there is a colour $\phi_C$ that is not used on any vertex of $C$.  These are exactly the requirements needed to apply  \cref{component_breaking} to $H-(A\cup I)$. Before we do so, we mention a subtle issue that we will have to deal with after its application, when we will reintroduce the vertices of $A$ coloured with their original colour from $\varphi$.
  

 We have applied \cref{lem:LW2} to the graph $H^{+-}_{\uparrow\bullet}$, which does not contain the vertices in $Z$.  This ensures that for each $j\in\{1,\ldots,r\}$, each $v\in (R_j\setminus Z)\cap L_{\uparrow\bullet}$ is assigned a colour in $\{1,\ldots,s+1\}\setminus \{\varphi(a):a\in N_H(v)\cap A_j\cup\{\alpha_j\}\}$ and $\alpha_j$ is assigned a colour distinct from the colours assigned to vertices in $A_j$.  However, the colour of $\alpha_j$ may appear at some vertex in $A_{y_j}\supseteq A_j$.\footnote{We cannot define $A_j=A_{y_j}$ because this would violate condition \cref{a_in_n_s} of admissible sets, which is required for the proof of \cref{lem:LW2}.}  For a vertex $v'\in Z\cap R_j$, it is not necessarily the case that $N_H(v')\cap A_{y_j}\subseteq A_j$.  In particular $v'$ may be adjacent to a vertex $a\in A_{y_j}\setminus A_j$ with $\varphi(a)=\varphi(\alpha_j)$.  If we apply \cref{component_breaking} to $H-(A\cup I)$, then the colour $\varphi(\alpha_j)$ will be used to break large monochromatic components of $H-(A\cup I)$. In this process, an unbounded number of vertices in $X_{v'}$ may receive the colour $\varphi(\alpha_j)$.  These vertices may all be adjacent to $a$, which would lead to a monochromatic component of colour $\varphi(\alpha_j)$ of unbounded size once we reintroduce the vertices of $A$ with their original colours from $\varphi$. We will thus need to recolour the vertices from $X_{v'}$ with a colour that is different from the colour breaker and the colours of the neighbours in the apex set. Fortunately, at least one such colour will always be available, as we will see shortly.

  So, we finally apply \cref{component_breaking} to the graph $H-(A\cup I)$ with the layering $\mathcal{L}$ and the colouring $\varphi$.  Let $\varphi'$ be the resulting colouring of $H-(A\cup I)$.  Then $\varphi'$ is an $(s+1)$-colouring of $H-( A\cup I)$ with clustering at most $(2s+6)c_4c_5$ and that satisfies conditions \ref{breaker_untouched_top}, \ref{breaker_untouched_even}, \ref{breaker_component_colours}, \ref{breaker_only_recolour_2vf} of \cref{component_breaking} for $c=s+1$ and $t=5$.

 We now reintroduce $A$ in order to extend $\varphi'$ from $H-( A\cup I)$ to $H-I$. 
 For each $a\in A$, let $\varphi'(a):=\varphi(a)$.  Recall that, for each $j\in\{1,\ldots,r\}$, \cref{lem:LW2} ensures that $\varphi(\alpha_j)\neq \varphi(a)$ for each $a\in A_j$.   For each $w\in L_{\ge 6}$, there exists $v\in L_{\uparrow\bullet}$ such that $w\in V(X_v)$.  Assume first that $v\in R_j\cap L^{y_j-}_{\uparrow\bullet}$, for some $j\in\{1,\ldots,r\}$. Then 
 $N_{H^-_{\uparrow\bullet}}(v)\cap A\subseteq A_j$.  For each such vertex $v$, each $w\in V(X_v)\setminus I$, has $\varphi'(w)\in\{\varphi(w),\varphi(\alpha_j)\}$, and thus $\varphi'(w)\neq\varphi'(a)$ for each $a\in N_H(w)\cap A$. Consider now the case where $w\in V(X_v)$ with $v\in R_j\cap (L^{y_j}_{\uparrow\bullet}\setminus L^{y_j-}_{\uparrow\bullet})\subseteq Z$, for some $j\in\{1,\ldots,r\}$. As explained above, we have $\varphi'(w)=\varphi(a)$ for some $a\in N_H(w)$ precisely if $a\in A_{y_j}\setminus A_j$ and $\varphi'(w)=\varphi(\alpha_j)\ne \varphi(v)$, that is $w$ was recoloured with a colour breaker during the application of  \cref{component_breaking}. For every such vertex $v\in R_j \cap Z$, we choose an arbitrary colour $\phi(v)$ in the set $\{1,\ldots,s+1\}\setminus\{\varphi(a):a\in N_{H_{\uparrow\bullet}}[v]\}$.  Such a colour exists because $N_{H_{\uparrow\bullet}}[v]$ has size at most $|A_{y_j}|+1 + |N_{\overline{H^{y_j}_{\uparrow\bullet}}}(v)|\le (s-3)+1+2=s$. We then recolour with colour $\phi(v)$ every $w\in X_v$ such that $\varphi'(w)=\varphi(\alpha_j)\ne \varphi(v)$. Let $\varphi''$ denote the resulting colouring of $H-I$.

 We now argue that $\varphi''$ still properly extends the precolouring of $S$ and has clustering at most $(2s+6)c_4c_5$. The former follows from the fact that $\varphi$ was properly extending the colouring of $S$ and newly coloured vertices have been assigned colours distinct from the  colours of their neighbours in $A$. To prove the latter, consider a monochromatic component $M$ in $\varphi''$. If it contains an element of $A$ then it must be also a monochromatic component of $\varphi$ in $H_{\uparrow}$, and thus it has size at most $c_4$. Otherwise $M$ is disjoint from $A$ and then it must be a monochromatic component in $H-(A\cup I)$. If $M$ does not intersect any vertex $w$ such that $\varphi''(w)\ne \varphi'(w)$ (that is, such that $w$ was recoloured between $\varphi'$ and $\varphi''$), then $M$ is a monochromatic component of   $H-(A\cup I)$ with respect to $\varphi'$ and thus has size at most $(2s+6)c_4c_5$, by definition of $\varphi'$. So we can assume that $M$ contains some vertex $w$ such that $w\in X_v$ and $v\in R_j\cap Z$, for some $j\in \{1,\ldots,r\}$, and moreover $\varphi''(w)=\phi(v)$. By definition, $\phi(v)$ is distinct from the colours appearing in the closed neighbourhood of $v$ in the colouring $\varphi$ of $H_{\uparrow\bullet}$. Since $Z$ is an independent set, for each neighbour $u$ of $v$ in $H_{\uparrow\bullet}-A$, the vertices of $X_{u}$ are not recoloured between $\varphi'$ and $\varphi''$, and therefore $M$ is contained in $X_v$. Since every neighbour of $w$ that is coloured $\phi(v)$ is in the same layer as $w$ (by condition \ref{breaker_only_recolour_2vf} of \cref{component_breaking}), $M$ lies in the intersection of $V(X_v)\setminus I$ with a single layer of the layering $\mathcal{L}$. Recall that we have proved that $V(X_v)\setminus I$ is $c_5$-skinny with respect to $\mathcal{L}$. It follows that $M$ has size at most $c_5$, as desired.
 
Finally, we properly extend $\varphi''$ from $H-I$ to $H$ by assigning colours to the vertices of $I_m,I_{m-1},\ldots,I_1$. 
Assign each vertex $v\in I_i$ a colour not assigned to its at most three already-coloured neighbours in $V(H-(A\cup \bigcup_{j=1}^i I_j))$ and its at most $s-3$ neighbours in $A$. This is possible since there are $s+1$ colours. For each $i\in\{m,\ldots,1\}$, the colouring of $H-(\bigcup_{j=1}^{i-1} I_j)$ properly extends the colouring of $H-(\bigcup_{j=1}^i I_j)$.  Therefore, this does not increase the size of any existing monochromatic component and does not create any new monochromatic component of size greater than $\max_{i=1}^m|I_i|\le c_5\le (2s+6)c_4c_5$.

  Finally, we argue that $\varphi''$ is a colouring of $H$ that properly extends the given precolouring of $S$. First observe that $\varphi''(v)=\varphi(v)$ for each $v\in L_{\le 5}$ and that $L_{\le 5}\subseteq V(H^-_{\uparrow\bullet})$.  Thus $\varphi''(v)=\varphi(v)\notin\{\varphi(a):a\in N_H(v)\cap S\}$ for each $v\in L_{\le 5}$, since \cref{lem:LW2} guarantees this property for the colouring $\varphi$ of $H^{+-}_{\uparrow\bullet}$. For a vertex $v\in L_{\ge 6}$, $N_H(v)\cap S\subseteq N_H(v)\cap A$, since $S\subseteq A\cup L_1\cup L_2$ is contained in the top of $G$.  As established above, $\varphi''(v)\neq\varphi(a)$ for each $a\in N_H(v)\cap A$. Thus, $\varphi''$ properly extends the precolouring of $S$. This completes the proof.
\end{proof}

We now prove the main theorem.

\begin{proof}[Proof of \cref{Jst}]
  First consider the case $s\ge 3$. We may assume that $t\geq s+2$. Let $J$ be a $\JJ_{s,t}$-minor-free graph.  By \cref{JstCurtain}, $J$ is a lower-minor-closed tree of $\mathcal{P}$-connected upward-connected $(s-3,\ell)$-curtains for some $\ell=\ell(s,t)$. \Cref{Heart2} then implies there is a function $c$ such that for each curtain $\torso{J}{B_x}$ of $J$, any $(s+1)$-colouring of at most $s$ vertices in the top of the root torso of $\torso{J}{B_x}$ can be properly extended to an $(s+1)$-colouring of the lower torso $\ltorso{J}{B_x}$ with clustering at most $c(s,t)$. Since $J[B_x]$ is a subgraph of $\ltorso{J}{B_x}$, \cref{Lungs-curtain} then implies that $J$ is $(s+1)$-colourable with clustering $c(s,t)$.

  Now consider the $s=2$ case. Let $G$ be a $\JJ_{2,t}$-minor-free graph. Since the complete join $K_2\oplus P_t$ is a planar graph in $\JJ_{2,t}$, by the Grid Minor Theorem~\citep{RS-V}, $G$ has treewidth at most some $\ell_t$. We apply the inductive argument of \cref{lem:LW2} with $s=2$ and $r=0$ (taking $\mathcal{S}=\emptyset$ and $P=\emptyset$). Although \cref{lem:LW2} formally states $s\geq 3$ to accommodate $(s-3)$-admissible sets, when $r=0$ there are no apex pairs $(\alpha_j,A_j)$ or trigger sets $S_j$, rendering the admissibility hypothesis vacuous. Specifically, if $G$ contains no $K_{2,y}$-subgraph for a sufficiently large $y$, \cref{lem:LW} applies directly. Otherwise, we find a $K_{2,y}$-subgraph with parts $X=\{x_1,x_2\}$ and $Y$. \Cref{no_surface_in_x} is vacuously true, and \cref{less_than_t_per} holds since each component of $G-X$ contains at most $t-1$ vertices of $Y$ (otherwise $G$ would contain a $K_2\oplus P_t$ minor). Cases~A and~B proceed identically, as the only role of $\mathcal{S}$ in the induction is ensuring the precoloured sets avoid $\bigcup_{j=1}^r S_j$, which is trivially empty. Thus, $G$ is $3$-colourable with clustering bounded by a function of $t$.

  For the case $s=1$, we may colour each connected component of $G$ independently, so assume $G$ is connected. Select an arbitrary vertex $v\in V(G)$ and let $L_i:=\{x\in V(G):\dist_G(v,x)=i\}$ for each $i\ge 0$. Colour each vertex in $L_i$ with $i\bmod 2$. Any monochromatic component $C$ is contained in a single layer $L_i$. In a BFS tree rooted at $v$, contracting the paths from $v$ to each vertex of $C$ yields a minor isomorphic to $K_1 \oplus T$ for some spanning tree $T$ of $C$. Since $K_1 \oplus T \in \JJ_{1,|V(C)|}$ and $G$ is $\JJ_{1,t}$-minor-free, we have $|V(C)| < t$.
\end{proof}

 \subsection{\boldmath \texorpdfstring{\cref{apex1}}{Theorem~7}: The \texorpdfstring{$s=3$}{s=3} Case}
 \label{KstFree3}
 \label{apex1proof3}

We now complete the proof of \cref{apex1}.  \Cref{apex1proof} already establishes \cref{apex1} for $s\ge 4$.  Observe that, in this proof, the restriction $s\ge 4$ appears because the graph $H_{\uparrow}^+$ obtained from $H_{\uparrow}$ by adding a vertex $\alpha$ may contain arbitrarily large $K_{3,t}$-subgraphs. This is because $G[L_5\cup L_{\uparrow}]$ may contain arbitrarily large $K_{2,t}$-subgraphs with the $t$ `right-part' vertices in $L_{\uparrow}$.  The additional vertex $\alpha$ then completes this into a $K_{3,t}$-subgraph, so \cref{lem:LW0} cannot be applied with $s=3$. We encountered a similar issue in the proof of \cref{Heart2} and used \cref{no_more_k2t} to eliminate $K_{2,t}$-subgraphs from $G[L_5\cup L_{\uparrow}]$. The following paragraph explains how the same techniques can be incorporated into the proof of \cref{apex1} in order to handle the case $s=3$.

\begin{proof}[Proof sketch of \cref{apex1} for $s=3$]
  By \cref{ApexMinorFreeStructure}, the $X$-minor-free $K_{3,t}$-sub\-graph-free graph $J$ that we want to colour is a lower-minor-closed tree of $\mathcal{P}$-connected upward-connected $(0,\ell)$-curtains described by $\TT_0:=(B_x:x\in V(T_0))$, for some $\ell:=\ell(X)$.  By \cref{Lungs-curtain}, it suffices to show that for each curtain $\torso{J}{B_x}$ any precolouring of any set $S$ of at most three vertices in the top of the root torso of $\torso{J}{B_x}$ can be properly extended to a $4$-colouring of $\ltorso{J}{B_x}$ with clustering at most $f(X)$.  The remainder of this discussion shows how to do this using (a simpler version of) the steps used in the proof of \cref{Heart2}.

  Let $G:=\torso{J}{B_x}$ and let $H:=\ltorso{J}{B_x}$ for some $x\in V(T_0)$. Let $\mathcal{T}:=(C_y:y\in V(T))$ be the tree-decomposition that describes the $(0,\ell)$-curtain $G$ (and therefore also describes $H$). For each $y\in V(T)$, let $G^y:=\ttorso{J}{B_x}{C_y}$ and let $H^y:=\ltorso{\ltorso{J}{B_x}}{C_y}$.  For each $y\in V(T)$, let $(\PP_y,\LL_y:=(L^y_1,L^y_2,\ldots))$ be the neat $\PP$-connected upward-connected $(w,\ell)$-layered partition that describes the drape $G^y$.  For each $i\in\NN$, let $L_i:=\bigcup_{y\in V(T)} L^y_i$ and let $\LL:=(L_1\cup L_2, L_3, L_4,\ldots)$, so that $\LL$ is a layering of $G$ and of $H$.

  We raise the $(0,\ell)$-curtain $G$ as in the proof of \cref{Heart2} and perform the same contractions in $H$ to obtain $H_{\uparrow}$.  Next, define the graphs $H_{\uparrow\bullet}$, and $H^-_{\uparrow\bullet}$ exactly as in the proof of \cref{Heart2}.
  The graph $H^{+-}_{\uparrow\bullet}$ is then obtained from $H^-_{\uparrow\bullet}$ by adding a single vertex $\alpha$ that dominates $L_{\uparrow\bullet}$ and setting the colour of $\alpha$ to be $m:=4$.
  By \cref{no_more_k2t}, $H^-_{\uparrow\bullet}[L_5\cup L_{\uparrow\bullet}]$ is $K_{2,t'}$-subgraph-free for some $t':=t'(\ell)$. Therefore $H_{\uparrow\bullet}^{+-}$ is $K_{3,\max\{t',t\}}$-subgraph-free.  The treewidth of $H_{\uparrow\bullet}^{+-}$ is not more than the treewidth of $H_{\uparrow}$ plus $1$ (due to the addition of $\alpha$). Thus, by \cref{raised_treewidth}, $H_{\uparrow\bullet}^{+-}$ has treewidth at most $k:=k(\ell)$.

  Since $H_{\uparrow\bullet}^{+-}$ has bounded treewidth and is $K_{3,\max\{t',t\}}$-subgraph-free, it can be $4$-coloured using \cref{lem:LW0} with the precoloured set $S\cup\{\alpha\}$.  Since no vertex in $H_{\uparrow\bullet}^{+-}$ has more than $3$ neighbours in $S$, the resulting colouring properly extends the precolouring of $S\cup\{\alpha\}$ and has clustering at most some $f(\ell)$. We extend this to a $4$-colouring of $H$ in several steps. First, exactly as in the proof of \cref{Heart2}, we extend the colouring to $H_{\uparrow\bullet}$ by assigning to each deleted vertex $u \in V(H_{\uparrow\bullet}) \setminus V(H^-_{\uparrow\bullet})$ a colour different from its at most two neighbours in $H^-_{\uparrow\bullet}$ and from $\alpha$. Then, for each connected subgraph $X_v$ of $H$ that is contracted to obtain a vertex $v$ of $H_{\uparrow\bullet}$, set the colour of each vertex in $X_v$ to the colour of $v$.

  At this point, $H$ has $\varphi$-monochromatic components that are neither of bounded size nor $f(\ell)$-skinny with respect to $\LL$.  As in the proof of \cref{Heart2}, we deal with this by removing islands, applying \cref{skinny_components} to establish that all $\varphi$-monochromatic components are skinny, and then applying \cref{component_breaking} (with $t=5$) to break large $\varphi$-monochromatic components using the colour ($4$) of $\alpha$.  Briefly, we remove a sequence of $(3,f(\ell))$-islands $I_1,\ldots,I_m$ from $H$, each of which is contained in $V(X_v)$ for some $v\in V(H_{\uparrow\bullet})$ until no such $(3,f(\ell))$-islands remain (letting $I:=\bigcup_{i=1}^m I_i$).  By \cref{skinny_components}, each monochromatic component in the resulting colouring of $H-I$ is $f(\ell)$-skinny with respect to the layering $\LL$ of $H$.  Each vertex in $L_{\ge 6}$ avoids the colour $\varphi(\alpha)=4$, so  \cref{component_breaking} can be used to ensure that all $\varphi$-monochromatic components of $H-I$ have bounded size.  Finally, we colour the islands $I_m,\ldots,I_1$ in this order so that the colouring of 
  $H-\bigcup_{j=i-1}^m I_j$ properly extends the colouring of 
  $H-\bigcup_{j=i}^m I_j$ for each $i\in\{2,\ldots,m+1\}$. This does not increase the size of any $\varphi$-monochromatic component and does not create any new $\varphi$-monochromatic component of size greater than $\max_{i=1}^m |I_i|\le f(\ell)$, and the resulting colouring still properly extends the original precolouring of the set $S$.
\end{proof}

\section{Algorithms}
\label{algorithms}

All the results in this paper are algorithmic and give polynomial-time algorithms for finding the colourings promised by \cref{Kh,Kst,Jst,apex1}. More specifically, there exists a constant $c$ such that for $n$-vertex graphs:
\begin{compactitem}
  \item the colouring in \cref{Kh} can be found in $f(h)\cdot n^c$ time for some $f \colon \NN\to\NN$;
  \item the colourings in \cref{Kst,Jst} can be found in $f(s,t)\cdot n^c$ time for some $f \colon \NN^2\to\NN$; and
  \item the colouring in \cref{apex1} can be found in $f(|X|,s,t)\cdot n^c$ time for some $f \colon \NN^3\to\NN$.
\end{compactitem}
We now elaborate on the steps required to compute these colourings:
\begin{enumerate}
    \item The tree-decomposition $\TT_0:=(B_x:x\in V(T_0))$ and almost-embeddings of \cref{DvoThoOriginal} can be computed in $f(|X|)\cdot n^c$ time \citep{DvoTho}.

    \item From this decomposition, we easily derive the tree of curtains decomposition $\TT:=(B_x:x\in V(T))$ of $J$. The decomposition itself is read off from the tree-decomposition $\TT_0$ as are the tree-decompositions $\TT_x:=(C_y:y\in V(T_x))$ that describe the curtain $\torso{J}{B_x}$ for each $x\in V(T)$. The layered partition of each drape is obtained using existing algorithms for product structure of almost-embedded graphs without major apex vertices \citep{DJMMUW20,BMO22}.

    \item Raising each curtain $\torso{J}{B_x}$ to obtain $G_{\uparrow}:=\torso{J}{B_x}_\uparrow$ and $H_{\uparrow}:=\ltorso{J}{B_x}_\uparrow$ is easily done in time linear in the size of the curtain.

    \item In each drape $G^y:=\ttorso{J}{B_x}{C_y}$ of each raised curtain $\torso{J}{B_x}$, we first identify a maximal set of disjoint $L^y_{\le 4}$-avoiding $3$-discs and $4$-discs in the surface-embedded graph $\overline{H}^y:=\overline{\ltorso{\ltorso{J}{B_x}}{C_y}}$. These can easily be identified in $O(n^8)$ time by enumerating all of the $3$- and $4$-cycles (each of which defines at most $2$ discs) and making a directed acyclic graph with these discs as vertices in which the direction of each edge is determined by the containment relationship. The sources in this graph give the desired set of discs.  This allows us to compute the $L^{y}_{\le 4}$-compression $\overline{H^{y}}_{\uparrow\bullet}$ as well as the related graph $H^y_{\uparrow\bullet}$ and $H^{+-}_{\uparrow\bullet}$.  
    
    \item The bounded treewidth results in \cref{lem:LW,lem:LW2} are applied to the graph $H^{+-}_{\uparrow\bullet}$. The proofs of these two results are inductive and give easy recursive polynomial-time algorithms.
    
    \item Identifying the islands $I_1,\ldots,I_m$ that must be removed from connected subgraphs $X_v$, for $v\in L_{\uparrow\bullet}$ is done as follows (letting $I:=\bigcup_{i=1}^m I_i$):  If a subgraph $X_v-I$ is not $c_5$-skinny then, by \cref{three_islands}, it contains a $(3,c_5)$-island $I_{m+1}$ in the appropriate surface-embedded graph $\overline{H^y}-I$.  The island can be found and eliminated using the separator method of \citet[Theorem~8]{DN17} outlined in \cref{three_islands}. This increases the total size of all islands $I_1,\ldots, I_m$ and so we can repeat the process a total of at most $n$ times.

    \item The other steps---lowering the curtains, extending the colouring, and breaking long skinny components using \cref{component_breaking}---are easy to implement in linear time. 
\end{enumerate}

\subsection*{Acknowledgements}

This research was initiated at the Barbados Graph Theory Workshop held at the Bellairs Research Institute in April 2019. Thanks to the other workshop participants for creating a productive working atmosphere. We would also like to thank Francis Lazarus for the helpful discussion about the proof of \cref{lem:planardisc}. Finally, thanks to Zden\v{e}k Dvo\v{r}\'{a}k for pointing out~\citep[Section~9]{RS-XVII} and for helpful conversations about the algorithmic aspects of \cref{DvoThoOriginal}.

{\fontsize{10pt}{11pt}
\selectfont
\def\soft#1{\leavevmode\setbox0=\hbox{h}\dimen7=\ht0\advance \dimen7
  by-1ex\relax\if t#1\relax\rlap{\raise.6\dimen7
  \hbox{\kern.3ex\char'47}}#1\relax\else\if T#1\relax
  \rlap{\raise.5\dimen7\hbox{\kern1.3ex\char'47}}#1\relax \else\if
  d#1\relax\rlap{\raise.5\dimen7\hbox{\kern.9ex \char'47}}#1\relax\else\if
  D#1\relax\rlap{\raise.5\dimen7 \hbox{\kern1.4ex\char'47}}#1\relax\else\if
  l#1\relax \rlap{\raise.5\dimen7\hbox{\kern.4ex\char'47}}#1\relax \else\if
  L#1\relax\rlap{\raise.5\dimen7\hbox{\kern.7ex
  \char'47}}#1\relax\else\message{accent \string\soft \space #1 not
  defined!}#1\relax\fi\fi\fi\fi\fi\fi}

}

\appendix
\section{Proof of \texorpdfstring{\cref{DvoThoCorollaryWithMinors}}{Lemma~12}}
\label{minors_in_proof}

In this appendix, we prove \cref{DvoThoCorollaryWithMinors}. Similar results are known in the literature~\citep[Section~9]{RS-XVII}, but we include our proof for completeness and precision. We begin with a version that requires only moderate changes to the tree-decomposition of \cref{DvoThoOriginal}.

\begin{lem}\label{DvoThoCorollary}
  For every integer $k\ge 1$ and every $k$-apex graph $X$ there exists an integer $\ell$ such that every $X$-minor-free graph $G$ has a rooted tree-decomposition $(B_x\colon x\in V(T))$ such that:
  \begin{enumerate}[(1)]
    \item for each $x\in V(T)$, the torso $\torso{G}{B_x}$ is a $(k-1,\ell)$-almost-embedded graph; and
    \item for each edge $xy$ of $T$ where $y$ is the parent of $x$;
    \begin{enumerate}[(a)]
      \item $B_x\cap B_y$ is contained in the top of $\torso{G}{B_x}$,
      \item $B_x\cap B_y$ is contained in the near-top of $\torso{G}{B_y}$ or
   $|B_x\cap B_y|\le k+2$; and
      \item $B_x\cap B_y$ contains at most three vertices not in the top of $\torso{G}{B_y}$.
    \end{enumerate}
  \end{enumerate}
\end{lem}

\begin{proof}
  First, apply \cref{DvoThoOriginal} to obtain a tree-decomposition $(B_x\colon x\in V(T))$ of $G$ in which each torso $\torso{G}{B_x}$ satisfies (i) and (ii).  For each $x\in V(T)$, denote by $(A_x,\hat{A}_x,G^x_0,G^x_1,\ldots,G^x_r)$  the $(k-1,\ell_0)$-almost-embedding of $\torso{G}{B_x}$. For convenience, assume the number of vortices equals $\ell_0$. We keep the same tree-decomposition $(B_x\colon x\in V(T))$.  In order to establish (2) we will, for each node $x\in V(T)$, create at most five additional vortices in $\torso{G}{B_x}$ that are vertex-disjoint from each other and from existing vortices. After these modifications, $\torso{G}{B_x}$ is $(k-1,\ell_0+5)$-almost-embeddable, so (1) is satisfied with $\ell=\ell_0+5$.

  Let $x$ be a non-root node of $T$ and let $y$ be the parent of $x$. Suppose that $B_x\cap B_y$ is not contained in the top of $\torso{G}{B_x}$. Let $C$ be the set of non-top vertices of $B_x\cap B_y$. By definition, $C$ induces a complete subgraph $K$ in $\torso{G}{B_x}$ and, since $C$ contains no vertices in the top of $\torso{G}{B_x}$, the vertices and edges of $K$ are all contained in the embedded part $G_0^x$ of $\torso{G}{B_x}$.  By (ii), every 3-cycle in $G^x_0$ bounds a 2-cell face. Thomassen’s 3-Path Property~\cite[Proposition~3.5]{Thom90} states that if any three internally disjoint paths $P_1, P_2, P_3$ from a vertex $u$ to a vertex $v$ in an embedding ($G^x_0$ in our case) are such that two of the three cycles $C_{i,j} := P_i \cup P_j$ ($1 \leq i < j \leq 3$) are contractible, then all three cycles $C_{i,j}$ are contractible (or equivalently, bound an open topological disc on the surface). Assume for the sake of contradiction that $K$ contains a cycle, $J$, that does not bound a topological disc in the embedding of $G^x_0$ and let $J$ be a shortest such cycle. By (ii),  $|J|\geq 4$. Since $J$ is a subgraph of the clique $K$, $J$ has a chord in $K$. The two cycles defined by the chord and the two half-cycles of $J$ are contractible, by the minimality of $J$. Thus by Thomassen’s 3-Path Property, $J$ is also contractible, which is a contradiction. It then follows from \cref{lem:planardisc} that the embedding of $K$ in $G^x_0$ is in a topological disc. Since no complete graph on five or more vertices has an embedding in a topological disc, $|V(K)|\leq 4$.  For each vertex $v\in C$, define a trivial vortex that contains only $v$.  Doing this for each $x\in V(T)$ ensures that (2a) is satisfied and introduces at most four new vortices in each torso $\torso{G}{B_x}$.

  To establish (2b), suppose that $B_x\cap B_y$ contains a vertex $v$ not in the near-top of $\torso{G}{B_y}$.  Then $C:=(B_x\cap B_y)\setminus A_y\subseteq N_{G^y_0}(v)$.  If $|C|\le 3$, then $|B_x\cap B_y|\le |A_y|+|C|\le k+2$ and there is nothing more to prove.  Suppose, therefore, that $|C|\ge 4$. Now $C$ induces a complete subgraph $K$ in $G^y_0$ so, by the argument in the previous paragraph, $|C|=4$.  Since each of the four triangles of $K$ is the boundary of a 2-cell face, $G^y_0$ is embedded in the sphere, and $V(G^y_0)=V(K)=C$. In this case, we create a trivial vortex in $\torso{G}{B_y}$ that contains one vertex of $G^y_0$.  With this new vortex, every vertex of $B_y$ is contained in the near-top of $\torso{G}{B_y}$, so (2b) is again trivially satisfied.

  To establish (2c), suppose that $B_x\cap B_y$ contains a set $C$ of four vertices not in the top of $\torso{G}{B_y}$.  Since no vertex of $C$ is in the top of $\torso{G}{B_y}$, each edge of the complete graph $\torso{G}{B_y}[C]$ is an edge of $G^y_0$.
 By the argument above, this implies that $G^y_0=K_4$ and is embedded on the sphere.  As in the previous paragraph we handle this by creating a trivial vortex that contains a single vertex of $C$.  (Note that establishing both (2b) and (2c) requires the addition of at most one vertex per torso.)
\end{proof}

For an edge $xy$ in a tree $T$, let $T_{x:y}$ be the component of $T-xy$ that contains $x$. For a tree-decomposition $\TT:=(B_x:x\in V(T))$ of a graph $G$, let $G_{x:y}^+:=G[\bigcup_{z\in V(T_{x:y})} B_z]$ and let $G_{x:y}:=G^+_{x:y}-(B_x\cap B_y)$.  The next lemma shows that subgraphs of almost-embedded graphs are also almost-embedded.

\begin{lem}\label{almost_subgraph}
   Let $G$ be a graph, let $G'$ be a subgraph of $G$, and let $(A,\hat{A},G_0,G_1,\ldots,G_r)$ be an $(a,\ell)$-almost-embedding of $G$.  Then $G'$ has an $(a,2\ell+1)$-almost-embedding $(A',\hat{A}', G_0',G_1',\ldots,G_r')$, on the same surface $\Sigma$ as $G$, where $A'=A\cap V(G')$, $\hat{A}'=\hat{A}\cap V(G')$, and $G_i'=G_i[V(G')]$ for each $i\in\{1,\ldots,r\}$.
\end{lem}

\begin{proof}
    First, consider the case when $G'$ is an induced subgraph of $G$.  The statement of the lemma already defines $A'$, $\hat{A}'$, and the graphs $G_1',\ldots,G_r'$ that form the vortices of $G'$.  There are two related issues:
    \begin{inparaenum}[(i)]
      \item the statement of the lemma does not define the graph $G_0'$, and
      \item each vortex $G_i$ has a path-decomposition whose bags are indexed by $V(G_i)\cap V(G_0)$, so removing vertices in $G_i$ makes (the path-decompositions of) $G_i'$ undefined.
    \end{inparaenum}
    We use the freedom provided by (i) to deal with the problem raised by (ii). The only thing we change is that for some vortices $G_i$ some interior vertices in $G_i-V(G_0)$ will be moved to $G_0'$ in order to ensure that the vortices $G_i'$ have an appropriate path-decomposition of width at most $2\ell+1$.

    Let $X:=V(G)\setminus V(G')$ be the set of vertices that are removed from $G$ to obtain $G'$. For each $i\in\{1,\ldots,r\}$, let $X_i:=X\cap V(G_i)$ be the vertices removed from vortex $G_i$ and let $Y_i:=X_i\cap V(G_0)$ be the set of vertices removed from the boundary of $G_i$ .

    Fix some $i\in\{1,\ldots,r\}$, let $B_i:=\{v_1,\ldots,v_p\}:=V(G_i)\cap V(G_0)$ be the vertices on the boundary of $G_i$, and let $(C_1,\ldots,C_p)$ be the path-decomposition of $G_i$, where $C_j$ is the bag associated with $v_j$ for each $j\in\{1,\ldots,p\}$.  First remove the elements of $X_i\setminus Y_i$ from $G_i$ and from each $C_1,\ldots,C_p$.  Now $C_1,\ldots,C_p$ is a path-decomposition of $G_i-(X_i\setminus Y_i)$ whose width has not increased.  It remains to remove the vertices in $Y_i$.

    First consider the easy case:  If $v_j\in Y_i$ and $C_j$ contains a vertex $w\not\in B_i$, then redefine $w$ to be a vertex of $G_0'$ and embed it at the same location as $v_j$. Leave the contents of the bag $C_j$ unchanged (but $C_j$ is now associated with $w$). Now $w$ is still a vertex of $G_i'$, but it is a boundary vertex of $G_i'$ so it is added to $B_i$. Note that each of the edges incident to $w$ belongs to $G_i'$ or has an endpoint in $\hat{A}$, so $w$ is an isolated vertex of $G_0'$ placed at the same point as $v_j$, which is not a vertex of $G_0'$.  This modification does not increase the width of the path-decomposition $(C_1,\ldots,C_p)$.

    After handling all the easy cases described in the previous paragraph, we are in the situation where $C_j\subseteq B_i$ for each $v_j\in Y_i$.  We handle these deletions in batches. Consider a maximal interval $\{a,a+1,\ldots,b\}$ such that $\{v_a,v_{a+1},\ldots,v_b\}\subseteq Y_i$.  If $a=1$ and $b=p$ then $G_i'=G_i-Y_i$ is the empty graph and we just remove it.

    We can therefore assume without loss of generality that $b < p$. If $a=1$ then, for any $j \le b$, $C_j\setminus\{v_1,\ldots,v_b\}\subseteq C_{b+1}$
    and we can simply remove $v_1,\ldots,v_b$ from all bags and remove the bags $C_1,\ldots,C_b$ from the path-decomposition. What remains is still a path-decomposition since, for any edge $vw$ of $G_i'$ with $v,w\in C_j$ and $j\le b$, $v$ and $w$ are also contained in $C_{b+1}$.

    Thus, we may assume $1 < a \le b < p$.  Then, for any edge $vw$ of $G_i'$ with $v,w\in C_j$ and $a\le j\le b$, $\{v,w\}\subseteq C_{a-1}\cup C_{b+1}$.  To handle this case, we remove $v_a,\ldots,v_b$ from all bags, we remove the bags $C_a,\ldots,C_b$ from the path-decomposition, and we set $C_{a-1}:=C_{a-1}\cup C_{b+1}$. This modification ensures that if some bag $C_i$ with $a\le i \le b$ contains two vertices $v_c,v_d$,  with $c<a$ and $d>b$, then the new set $C_{a-1}$ contains both $v_c$ and $v_d$, thus after the removal of the vertices $v_a,\ldots, v_b$ and  bags $C_a,\ldots,C_b$, the resulting sequence of bags is still a path-decomposition.
    The bag $C_{a-1}$, which immediately precedes the interval $\{a,\ldots,b\}$, now contains the union of two of the original bags of the decomposition so its size is at most $2(\ell+1)$.

    The step described in the preceding paragraph removes the interval $\{a,\ldots,b\}$ and modifies $C_{a-1}$ so that its size is at most $2(\ell+1)$. Repeat this, handling the intervals $\{\{a_j,\ldots,b_j\}\}_{j=1}^q$ by increasing order of left endpoint $a_j$.  After removing $\{a_j,\ldots,b_j\}$, $C_{a_j-1}$ has size at most $2(\ell+1)$ and $C_{a_{j'}}$ has size at most $\ell+1$ for each $j'\ge j+1$.  When this process completes, no bag has size greater than $2(\ell+1)$, so it  gives a $(a,2\ell+1)$-almost-embedding of $G'$ that satisfies the conditions of the lemma.

    Finally, suppose that $G'$ is not an induced subgraph of $G$. First, suppose $G$ and $G'$ differ by the deletion of exactly one edge.  Then the given $(a,\ell)$-almost-embedding of $G$ is also an $(a,\ell)$-almost-embedding of $G'$ that satisfies the conditions of the lemma. Since the parameters $a$ and $\ell$ of the almost-embedding do not change, we can use this fact repeatedly for any $G'$ in which $V(G')=V(G)$.  To handle the (general) case where $V(G')\neq V(G)$ and $G'$ is not an induced subgraph of $G$, we apply the first argument above
    on the induced graph $G[V(G')]$ to obtain an $(a,2\ell+1)$-almost-embedding of $G[V(G')]$ that satisfies the requirements of the lemma.  Then $G'$ is a subgraph of $G[V(G')]$ with $V(G')=V(G[V(G')])$ so the given embedding of $G[V(G')]$ is also an embedding of $G'$ that satisfies the requirements of the lemma.
\end{proof}

The first step in proving \cref{DvoThoCorollaryWithMinors} is to ensure that some of the subgraphs $G_{x:y}$ are connected (see \citep{FN06,GJNW23} for related results that are proven using similar techniques).

\begin{lem}\label{no_stupid_clique_sums_i}
    Let $\TT:=(B_x:x\in V(T))$ be a rooted tree-decomposition of a graph $G$. Then $G$ has a rooted tree-decomposition $\TT':=(B'_x:x\in V(T'))$ such that:
    \begin{compactenum}
        \item for any $xy\in E(T')$ with $y$ the parent of $x$,
        \begin{compactenum}
            \item $G_{x:y}$ is connected and non-empty; and
            \item $N_G(V(G_{x:y}))=B'_{x}\cap B'_{y}$;
        \end{compactenum}
        \item there is a mapping $\rho:V(T')\to V(T)$ such that:
        \begin{compactenum}
            \item for $x\in V(T')$, $B_x'\subseteq B_{\rho(x)}$; and
            \item for each $xy\in E(T')$, $\rho(x)\rho(y)\in E(T)$.
        \end{compactenum}
    \end{compactenum}
\end{lem}

\begin{proof}
    If each edge $xy$ of $T$ satisfies (1a) and (1b) then the tree-decomposition $\TT$ and the identity mapping $\rho(x):=x$ already satisfies the requirements of the lemma, and we are done.  Otherwise, let $y$ be a minimum-depth node in $T$ that has a child $x$ such that $G_{x:y}$ is not connected or $N_G(V(G_{x:y}))\neq B_x\cap B_y$.  Let $C_1,\ldots,C_k$ be the components of $G_{x:y}$.  (Note that, if $B_x\subseteq B_y$, then $G_{x:y}$ has no vertices, so $k=0$, and the subtree $T_{x:y}$ will not contribute anything to $T'$.)

    For each $i\in\{1,\ldots,k\}$, let $T_i$ be a copy of $T_{x:y}$, let $x_i$ be the copy of $x$ in $T_i$,
    and let $\TT_i:=(B_z\cap N_G[V(C_i)]:z\in V(T_i))$.\footnote{Since each $T_i$ is a copy of $T_{x:y}\subseteq T$, we abuse notation slightly here and use vertices of $T_i$ as indices for the bags of $\TT$.}  Create a new tree $T'$ by joining $T_{y:x}$ to each of $T_1,\ldots,T_k$ using the edge $x_iy$.  Then $\TT':=(B_z:z\in V(T'))$ is a tree-decomposition of $G$ and the edges $x_1y,\ldots,x_ky$ satisfy (1a) and (1b). Note that each node $x'$ in $T'$ is either a node of $T$ or a copy of some node $x$ of $T$. In the former case we set $\rho(x'):=x'$ and in the latter case we set $\rho(x'):=x$.  Clearly the mapping  $\rho$ satisfies (2a) and (2b). Repeat this step as long as some edge does not satisfy (1a) or (1b).

    This process eventually eliminates all edges $yx$ incident to $y$ that do not satisfy (1a) or (1b) since the number of such edges is reduced by 1 at each iteration.  The process eventually eliminates all edges whose upper endpoint has the same depth as $y$ because there are a finite number of nodes with the same depth as $y$ and this process does not introduce any new nodes at this depth. After this any subsequent nodes processed have depth greater than $y$.  This process eventually terminates because every tree-decomposition it produces uses a tree whose height is no more than that of $T$.
\end{proof}

We also use the following result. For a graph $G$ and distinct vertices $u,v,w\in V(G)$, we say $G$ has  a \hdefn{$\{u,v,w\}$}{rooted $K_3$ minor} if there are pairwise disjoint pairwise adjacent connected subgraphs $A,B,C$ in $G$ with $u\in V(A)$, $v\in V(B)$ and $w\in V(C)$.

\begin{lem}[\citep{LinusWood}]
\label{no_rooted_k3}
    Let $G$ be a graph with three distinguished vertices $u$, $v$, and $w$.  If $G$ does not contain a $\{u,v,w\}$-rooted $K_3$ minor then $G$ contains a vertex $q$ such that each component of $G-q$ contains at most one of $u$, $v$, or $w$.
\end{lem}

We now prove the main result of this appendix, which we restate for convenience.

\DvoThoCorollaryWithMinors*

\begin{proof}[Proof of \cref{DvoThoCorollaryWithMinors}]
    Let $\TT$ be the rooted tree-decomposition of $G$ guaranteed by applying \cref{no_stupid_clique_sums_i} to the tree-decomposition given by \cref{DvoThoCorollary}.
    Suppose that $\TT$  does not satisfy the requirements of the lemma.  Therefore, there exists $y\in V(T)$ of minimum depth such that there is no faithful $\ltorso{G}{B_y}$-model in $G$ in which each vertex in the top of $\torso{G}{B_y}$ is in a branch set of size one.  We will repeatedly choose an edge $uv$ in $\ltorso{G}{B_y}$ and show the existence of a subgraph in $G[G_{x:y}\cup\{u,v\}]$ that \defn{represents} the edge $uv$ (and possibly some other edges).  The subgraphs assigned this way are pairwise disjoint, except for the vertices in $B_y$, so that the union of all such graphs contains a model of $\ltorso{G}{B_y}$.  This process may fail for some child $x$ of $y$, in which case we will make adjustments to a subtree rooted at one of the children $x$ of $y$ that will lift a vertex $q\in V(G_{x:y})$ into $B_y$.  This will introduce new edges (incident to $q$) in $\ltorso{G}{B_y}$.  We will be able to represent these new edges and this change will eliminate the edges of the clique $B_x\cap B_y$ that we are unable to represent.

    Consider any edge $uv$ of $\ltorso{G}{B_y}$ that is not yet represented.  If $uv$ is an edge of $G[B_y]$ then we say that $uv$ is \defn{represented} (by $uv$).  Note that this applies in particular to any edge $uv$ of $\ltorso{G}{B_y}$ for which $u$ or $v$ is in the top of $\torso{G}{B_y}$.  This will ensure that the branch set $G_v$ for any vertex $v$ in the top of $\torso{G}{B_y}$ consists of a single vertex.

    Otherwise $uv$ is not an edge in $G[B_y]$, and thus the vertices $u$ and $v$ are part of at least one adhesion set $B_x\cap B_y$ for some edge $xy$ of $T$ where $y$ is the parent of $x$. For each torso $\torso{G}{B_x}$ in a tree-decomposition, let $\ttop{G}{B_x}$ denote the top of $\torso{G}{B_x}$.
 Since neither $u$ nor $v$ is in the top of $\torso{G}{B_y}$, \cref{DvoThoCorollary}(2c) implies that $B_x\cap B_y\setminus \ttop{G}{B_y}$ contains $u$, $v$ and at most one other vertex.

    First, suppose that there exists an edge $xy\in E(T)$ such that $B_x\cap B_y \setminus \ttop{G}{B_y}=\{u,v\}$.
    Since $G_{x:y}$ is connected and non-empty and $\{u,v\}\subseteq N_G(V(G_{x:y}))$, there is a path in $G$ from $u$ to $v$ whose internal vertices are contained in $G_{x:y}$. In this case we say that $uv$ is \defn{represented} (by $G_{x:y}$).

    Otherwise, there exists a $3$-cycle $uvw$ in $\torso{G}{B_y}\setminus\ttop{G}{B_y}$ and an edge $xy\in E(T)$ such that $B_x\cap B_{y}\setminus \ttop{G}{B_y} =\{u,v,w\}$.  Suppose that at least one of $vw$ or $uw$ (say $vw$) is already represented.
    The edge $vw$ may be represented by itself or by a subgraph $G_{x':y}$ with $B_{x'}\cap B_{y}\setminus \ttop{G}{B_y}\supseteq \{v,w\}$.  Importantly, $vw$ is not represented by $G_{x:y}$.  Since $G_{x:y}$ is connected, non-empty, and adjacent to each of $u,v,w$ we can contract $G_{x:y}$ into a single vertex $q$ that is adjacent to each of $u,v,w$. Then by contracting the edge $qu$ we obtain a graph that contains $uv$ and $uw$.  In this case we say that the path $vuw$ is \defn{represented} (by $G_{x:y}$) and that the edges $uv$ and $uw$ are \defn{represented} (by $G_{x:y}$).

    If neither $uw$ nor $vw$  are already represented, but there are distinct edges $xy$ and $x'y$ of $T$ such that $\{u,v,w\}\subseteq B_x\cap B_{y}$ and $\{u,v,w\}\subseteq B_{x'} \cap B_{y}$ then the path $uvw$ can be represented by $G_{x:y}$ and the path $vuw$ be represented by $G_{x':y}$.  In this case we say that the $3$-cycle $uvw$ and the edges $uv$, $vw$, and $uw$ are \defn{represented} (by $G_{x:y}$ and $G_{x':y}$).

    We are left with the case where none of $uv$, $vw$, or $uw$ are represented and there is exactly one edge $xy\in E(T)$ such that $B_{x}\cap B_{y}\setminus \ttop{G}{B_y}=\{u,v,w\}$.  Consider the graph $H:=G[V(G_{x:y})\cup\{u,v,w\}]$.  Since $G_{x:y}$ is connected and $N_G(V(G_{x:y}))=B_x\cap B_y$, $H$ is connected.  Since none of the edges of the cycle $uvw$ are represented yet, none of these edges is present in $H$, so $H$ contains at least four vertices. If $H$ contains a $\{u,v,w\}$-rooted $K_{3}$-minor, then we can again contract edges of $H$ to obtain the edges of $uvw$.  In this case, we say that the $3$-cycle $uvw$, and the edges $uv$, $vw$, and $uw$ are \defn{represented} (by $G_{x:y}$).

    Otherwise, by \cref{no_rooted_k3}, $H$ contains a vertex $q$ such that each component of $H-q$ contains at most one of $u$, $v$, or $w$.  Since $G_{x:y}$ is connected and adjacent to each of $u,v,w$ this implies that $q\not\in\{u,v,w\}$.  In this case, we will lift the vertex $q$ into $B_y$ by modifying our tree-decomposition $\TT$.

    We now describe these modifications, first showing how to modify $\TT$ to obtain a new tree-decomposition $\TT':=(B'_\tau:\tau\in V(T'))$ with $q\in B'_y$ in a way that preserves the properties of \cref{no_stupid_clique_sums_i}. We then show how to obtain almost-embeddings of each torso of $\TT'$ so that these almost-embeddings collectively satisfy the properties of \cref{DvoThoCorollary} (with a slightly larger value of $\ell$).  After this modification, $q$ will be in $B'_y$ and $\ltorso{G}{B'_y}$ will include the edges $uq$, $vq$, and $wq$ but will not necessarily include the edges $uv$, $vw$, or $uw$.  We show that each of the edges $uq$, $vq$, and $wq$ can be represented by one of the graphs $G_{c_i:y}$ where $c_i$ is the root of the one of the trees in $T'-V(T)$ that was introduced to replaces $x$.

    \paragraph{Modified Tree-Decomposition:}

    Let $z$ be the minimum-depth node of $T$ with $q\in B_z$ and let $z_0,\ldots,z_p$ be the path in $T$ from $z_0:=z$ to $z_p:=y$. We define a sequence of tree-decompositions $\TT_0,\TT_1,\ldots,\TT_p$ where $\TT_0=\TT$, $T^0:=T$ and $\TT_{i}:=(B^i_x:x\in V(T^i))$.  Refer to \cref{toot}. The trees $T^i$ and $T^{i-1}$ are almost identical except that the node $z_{i-1}$ that appears in $T^{i-1}$ will be replaced with several nodes $\{c_1,\ldots,c_k\}=V(T_i)\setminus V(T_{i-1})$ in $T^i$. Each of the children $a_1,\ldots,a_d$ of $z_{i-1}$ in $T^{i-1}$ will become the child of a node in $\{c_1,\ldots,c_k\}$.  For each node $\tau\in V(T^i)\cap V(T^{i-1})\setminus\{z_i\}$, we define $B^i_{\tau}:=B^{i-1}_\tau$. For $z_i$ we define $B^i_{z_i}=B^{i-1}_{z_i}\cup\{q\}$.

    \begin{figure}
      \begin{center}
        \begin{tabular}{cc}
          \includegraphics[page=1]{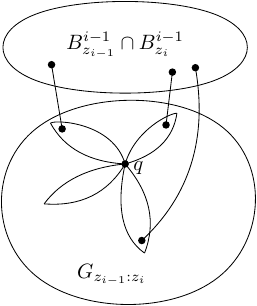} &
          \includegraphics[page=2]{figs/toot}
        \end{tabular}
      \end{center}
      \caption{Adding $q$ to $B_{z_i}$ may disconnect the graph $G_{z_{i-1}:z_i}$.}
      \label{toot}
    \end{figure}

    Let $C_1,\ldots,C_k$ be the connected components of $G_{z_{i-1}:z_i}-\{q\}$ and, for each $j\in\{1,\ldots,k\}$, let $G_j:=G[N_G[V(C_j)]]$.  In $T^i$ we replace $z_{i-1}$ with $k$ new nodes $c_1,\ldots,c_k$ and set $B^i_{c_j}:=B^{i-1}_{z_{i-1}}\cap V(G_{j})$ for each $j\in\{1,\ldots,k\}$. Let $a_1,\ldots,a_d$ be the children of $z_{i-1}$ in $T^{i-1}$. For each $j\in\{1,\ldots,d\}$, the graph $G_{a_j:z_{i-1}}$ is connected and its vertex set is a subset of $G_{z_{i-1}:z_i}-q$.  Therefore, $G_{a_j:z_{i-1}}$ is an induced subgraph of $C_s$ for some $s\in\{1,\ldots,k\}$.  In $T^i$, we make $a_j$ a child of $c_s$.  This completes the definition of $T^i$ and $\TT_i:=(B^i_\tau:\tau\in V(T^i))$.

    \paragraph{\boldmath Connectivity Conditions for $\TT_i$:}

    We claim that $\TT_i$ satisfies properties (2a) and (2b) of \cref{no_stupid_clique_sums_i}.
    That is, for each node $a\in V(T^i)$ with parent $b\in V(T^i)$, $G_{a:b}$ is connected and $N_G(V(G_{a:b}))=B_a\cap B_b$.  It suffices to check this condition for $a\in\{a_1,\ldots,a_d\}\cup\{c_1,\ldots,c_k\}$; for any other choice of $a$, the edge $ab$ appears in $T^{i-1}$ and $G_{a:b}$ is identical independent of whether it is defined in terms of $\TT_{i-1}$ or $\TT_i$.
    \begin{compactitem}
      \item  For each $j\in\{1,\ldots,k\}$, $C_j=G_{c_j:z_i}$ is connected and $N_G(C_j)=V(G_j)\setminus V(C_j)=B_{c_j}\cap B_{z_i}$.

      \item Now consider $G_{a_j:c_s}$ for some $j\in\{1,\ldots,d\}$ and (the relevant) $s\in\{1,\ldots,k\}$. We will show that $B_{a_j}\cap B_{z_{i-1}}\subseteq B_{a_j}\cap B_{c_s}\subseteq B_{a_j}\cap B_{z_{i-1}}$, which implies that $B_{a_j}\cap B_{c_s}=B_{a_j}\cap B_{z_{i-1}}$.  Since $T_{a_j:c_s}=T_{a_j:z_{i-1}}$ this shows that $G_{a_j:c_s}$ is connected and that $N_G(V(G_{a_j:c_s}))=B_{a_j}\cap B_{c_s}$.

      By definition, $B_{c_s}\subseteq B_{z_{i-1}}$, so $B_{a_j}\cap B_{c_s}\subseteq B_{a_j}\cap B_{z_{i-1}}$.  On the other hand,
      \[    V(G_s) = N_G[V(C_s)]
          \supseteq N_G[V(G_{a_j:z_{i-1}})]
          \supseteq N_G(V(G_{a_j:z_{i-1}})).
      \]
      Therefore $B_{c_s}=V(G_s)\cap B_{z_{i-1}} \supseteq N_G(V(G_{a_j:z_{i-1}}))\cap B_{z_{i-1}} = B_{a_j}\cap B_{z_{i-1}}$.  Therefore $B_{a_j}\cap B_{c_s}\supseteq B_{a_j}\cap B_{z_{i-1}}$.
    \end{compactitem}
    This process completes with a tree-decomposition $\TT_p=(B^p_\tau:\tau\in V(T^p))$ that satisfies conditions (2a) and (2b) of \cref{no_stupid_clique_sums_i} and in which $q\in B_{y}$.

    \paragraph{Almost-Embedding the New Torsos:}

    Now we explain how to almost-embed each torso of $\TT_i$ so that the torsos collectively satisfy the conditions of \cref{DvoThoCorollary}.  First, consider the torso $\torso{G}{B^i_z}=\torso{G}{B^i_{z_0}}$. This torso inherits an embedding from $\torso{G}{B^{i-1}_z}$ that includes $q$.  If $q$ is not in the top of $\torso{G}{B^{i-1}_z}$ then we make $q$ a trivial vortex in $\torso{G}{B^i_z}$, so that $q$ is in the top of $\torso{G}{B_{z}}$.

    Next, consider some torso $\torso{G}{B^i_{z_i}}$ for some $i\in\{1,\ldots,p-1\}$.  This torso inherits an embedding from $\torso{G}{B^{i-1}_{z_i}}$ that does not include $q$.  In the embedding of $\torso{G}{B^i_{z_i}}$, we make $q$ a non-major apex vertex and we make each vertex of $B^{i-1}_{z_{i-1}}\cap B^{i-1}_{z_i}$ that is not in the top of $\torso{G}{B^{i-1}_{z_i}}$ a top vertex of $\torso{G}{B^i_{z_i}}$ by creating a trivial vortex. By \cref{DvoThoCorollary}(2c), this results in the creation of at most three new vortices.  Note: It is important here that each vertex in $B^{i-1}_{z_{i-1}}\cap B^{i}_{z_{i}}$ (including $q$) is in the top of $\torso{G}{B^i_{z_i}}$.

    Next, consider one of the torsos $\torso{G}{B^i_{c_j}}$ where $c_j$ is one of the vertices in $V(T^i)\setminus V(T^{i-1})$ that replaces $z_{i-1}$.  Since $B^i_{c_j}\subseteq B^{i-1}_{z_{i-1}}$, this torso inherits an embedding from $\torso{G}{B^{i-1}_{z_{i-1}}}$ that includes $q$ in its top. We use this embedding as is.

    Finally, consider the torso $\torso{G}{B^p_y}=\torso{G}{B^p_{z_p}}$.  Recall that $uvw$ is the boundary of a $2$-cell face $D$ in the embedded part of $\torso{G}{B_y}$.  The torso $\torso{G}{B^p_y}$ inherits an embedding from $\torso{G}{B^p_y}$ that does not contain $q$.  To add $q$, we make a $q$ a vertex in the embedded part that is embedded in the interior of $D$.  We also embed each edge $qu$, $qv$, and $qw$ in $D$.  This completes the description of the embedding of $\torso{G}{B^i_\tau}$ for each $\tau\in V(T^i)$.

    \paragraph{\boldmath The New Torso $\torso{G}{B_y}$:}

    Now $\TT':=\TT_{p}$ is a tree-decomposition in which each torso is equipped with an almost-embedding.  Let $T':=T^p$, and let $B'_\tau:=B^p_\tau$ for each $\tau\in V(T')$.  Consider the torso $\torso{G}{B'_y}$.  In $T'$, the child $x\in V(T)$ of $y$ has been replaced by nodes $c_1,\ldots,c_k$ associated with the components $C_1,\ldots,C_k$ of $G_{x:y}-\{q\}$.  Each component $C_i$ is adjacent to $q$ and to at most one of $u$, $v$, or $w$. Since $N_G(B_{x})=B_{x}\cap B_{y}$, either $uq\in E(G)$ or some component $C_{i}$ is adjacent to both $q$ and $u$.  In the former case we say that the edge $uq\in E(\torso{G}{B_y})$ is \defn{represented} by itself. In the latter case, $uq$ is \defn{represented} by $C_i:=G_{c_{i}:y}$.  Similar comments hold for the edges $vq$ and $wq$ of $\torso{G}{B_{y}}$.

    Since each $C_i$ contains at most one of $u$, $v$, or $w$, none of the edges $yc_i$ causes an edge of the cycle $uvw$ to be in $\ltorso{G}{B_y}$.  Some of these edges may be present in $\torso{G}{B'_y}$, but this is caused by other children of $y$.  Any such edge will be dealt with when we consider unrepresented edges of $\ltorso{G}{B'_y}$ with the help of some child $x'$ of $y$ that is not one of $c_1,\ldots,c_k$.

    \paragraph{\boldmath Bounding the Increase in $\ell$:}

    Let $Z_{y,x}:=V(T')\setminus V(T)$.  (Note that this does not include the node $y=z_p$.)  Each vertex in $Z_{y,x}$ is obtained from $z_{i-1}$ for some $i\in\{1,\ldots,p\}$ by splitting $B^{i-1}_{z_{i-1}}$ into $B^i_{c_1},\ldots,B^i_{c_k}$, so that $z_{i-1}\in V(T^{i-1})$ is replaced by $c_1,\ldots,c_k\in V(T^i)$.  For $i\in\{1,\ldots,p-1\}$, each vertex of $B^i_{c_j}\cap B^i_{z_i}$ is contained in the top of $\torso{G}{B^i_{z_i}}$.  This implies that, for each edge $\tau\tau'\in T'[Z_{y,z}]$, $B'_{\tau}\cap B'_{\tau'}$ is in the top of $\torso{G}{B'_{\tau}}$ and in the top of  $\torso{G}{B'_{\tau'}}$.

    By construction, for each $\tau\in Z_{y,x}$, there exists a $\rho(\tau)\in V(T)$ such that $B_\tau\subseteq B_{\rho(\tau)}\cup\{q\}$.  The almost-embedding of $\torso{G}{B_\tau}$ inherits an almost-embedding from $B_{\rho(\tau)}$ and then does some modifications that increase the number of vortices by at most three and the number of non-major apex vertices by at most one.  By \cref{almost_subgraph}, the inherited embedding has at most $k-1$ major apex vertices, at most $\ell$ vortices, each of width at most $2\ell+1$, and at most $\ell$ apex vertices.  Since $\ell+3\le 2\ell+1$ for $\ell\ge 2$, the final embedding of $\torso{G}{B_\tau}$ is a $(k,2\ell+1)$-almost embedding, for each $\tau\in Z_{y,x}$.

    The preceding operation modifies the tree-decomposition $\TT$ in a way that only affects the subtree $T_{x:y}$ of $T$, and this only occurs at the unique child $x$ of $y$ such that $B_x\cap B_y\setminus \ttop{G}{B_y}=\{u,v,w\}$.  Therefore, we can perform this operation on each of the children $x_1,\ldots,x_c$ of $y$ for which it is required, and each such operation will only affect the subtree $T_{x_i:y}$.  Since $T_{x_1:y},\ldots,T_{x_c:y}$ are pairwise vertex-disjoint, $Z_{y,x_1},\ldots,Z_{y,x_c}$ are pairwise disjoint. For convenience, we still call the resulting tree $T'$ and the resulting tree-decomposition $\TT':=(B'_x:x\in V(T'))$. After doing this, every edge of the lower torso $\ltorso{G}{B_y}$ is represented.\footnote{The graphs $\torso{G}{B_y}$ and $\ltorso{G}{B_y}$ are defined here with respect to the new tree-decomposition $\TT'$.}  This ensures that $\ltorso{G}{B_y}$ is a faithful minor of $G$.

    Let $Z_y :=\bigcup_{i=1}^c Z_{y,x_i}$.
    We ran this process at $y$ because it was a minimum-depth node for which $\ltorso{G}{B_y}$ was not a faithful minor of $G$.  If there exists another node $y'$ for which $\ltorso{G}{B_{y'}}$ is not a faithful minor of $G$, then run this process again (for the minimum-depth such $y'$).  Because $y$ was a minimum-depth node, $y'$ is not an ancestor of $y$.  If  $y'\not\in Z_y$ then, again, $Z_y$ and $Z_{y'}$ will be disjoint because $T'[Z_y\cup\{y\}]$ and $T'[Z_{y'}\cup\{y'\}]$ are each connected and neither contains the root of the other.  If $y'\in Z_y$ then running the preceding process on $y'$ affects subtrees $T_{x_1':y'},\ldots,T_{x_c':y'}$ for one or more children $x_1',\ldots,x_c'$ of $y'$.  For each $i\in\{1,\ldots,c\}$, $B_{y'}\cap B_{x_i'}$ contains exactly three vertices not in the top of $\torso{G}{B_{y'}}$---namely three vertices $u',v',w'$ that form a $K_3$ in $\ltorso{G}{B_{y'}}$ but for which $G_{x_i':y'}$ has no $\{u',v',w'\}$-rooted $K_3$ minor. Suppose that $x_i'=\tau'$ for some $i\in\{1,\ldots,c\}$ and $\tau'\in Z_y$. Since $y'\neq y$, the parent $\tau$ of $\tau'$ is also in $Z_y$.  This implies that $x_i'\not\in Z_y$ since $B_{\tau}\cap B_{\tau'}$ is contained in the top of $\torso{G}{B_{\tau}}$.  Therefore $Z_y$ and $Z_{y'}$ are disjoint.

    Thus we can find a sequence of tree-decompositions $\TT_0=\TT,\TT_1,\ldots,\TT_d$ where $\TT_i:=(B^i_y:y\in V(T_i))$ and a sequence $y_1,\ldots,y_d$ of nodes, where $y_i\in \bigcap_{j=i-1}^d V(T_{j})$, $\ltorso{G}{B^{i-1}_{y_i}}$ is not a faithful minor of $G$, but $\ltorso{G}{B^{j}_{y_i}}$ is a faithful minor of $G$ for each $j\in\{i,\ldots,d\}$. By construction, for each $i\in\{1\ldots,d\}$, $Z_{y_i}\subseteq \bigcap_{j=i}^d V(T_{j})$ and the sets $Z_{y_1},\ldots,Z_{y_d}$ are pairwise disjoint.  The end result is a tree-decomposition $\TT_d:=(B^d_{z}:z\in V(T_d))$ of $G$.  For each $z\in V(T_d)$, the torso $\torso{G}{B_z}$ can be obtained by taking a subgraph of a torso of $\TT$ and then adding at most one new non-major apex vertex (referred to as $q$, above) and at most three new vortices (one for each vertex of $B^{i-1}_{z_{i-1}}\cap B^{i-1}_{z_i}$ not in the top, above).  As discussed above, this implies that, for $\ell\ge 2$, the resulting graph is $(k-1,2\ell+1)$-almost-embedded.  For each $x\in V(T)$, the faithful $\ltorso{G}{B_x}$-model $\{G_v:v\in B_x\}$ satisfies the ``furthermore'' clause of the lemma because the only vertices that are assigned a branch-set $G_v$ with more than one vertex (the vertices $u$, $v$, and $w$, above) are not in the top of the $\torso{G}{B_x}$.  Therefore, the tree-decomposition $\TT'$ satisfies the requirements of the lemma.
\end{proof}

The following theorem and its consequences are not used directly in the current paper, but have applications to other problems \cite{DHHJLMMRW25}.

\begin{thm}\label{apex_special_case}
    For every apex graph $X$, there exists a positive integer $\ell$ such that every $X$-minor-free graph $G$ has a rooted tree-decomposition $(B_x:x\in V(T))$ of adhesion at most $3$ such that for every $x\in V(T)$:
    \begin{enumerate}[(1)]
        \item $\ltorso{G}{B_x}$ is a minor of $G$,
        \item $\ltorso{G}{B_x}$ has a connected $\ell$-layered partition $(\PP_x,\LL_x)$, and
        \item if $y$ is the parent of $x$ then
        \begin{compactenum}
            \item every vertex in $B_x\cap B_y$ is contained in the first layer of $\LL_x$,
            \item no vertex in $B_x\cap B_y$ is contained in the first layer of $\LL_y$, and
            \item each vertex in $B_x\cap B_y$ is in a singleton part of $\PP_x$.
        \end{compactenum}
    \end{enumerate}
\end{thm}

\begin{proof}
    By \cref{ApexMinorFreeStructure} with $k=1$, $G$ is a lower-minor-closed tree of $(0,\ell_0)$-curtains described by $\TT:=(B_x:x\in V(T))$ where each torso $\torso{G}{B_x}$ is equipped with the connected $\ell_0$-layered partition $(\PP^0_x,\LL_x)$ for some $\ell_0:=\ell_0(X)$. This already satisfies (1), (2), and (3a). Although not stated in \cref{ApexMinorFreeStructure}, the construction used in the proof of \cref{ApexMinorFreeStructure} only makes $x$ a child of $y$ if $B_x\cap B_y$ contains a vertex $v$ not in the near-top of $\torso{G}{B_y}$ (this is the definition of a $k$-light edge). By property (1) of \cref{DvoThoCorollaryWithMinors}, the almost-embedding $(A^y,\hat{A}^y,G^y_0,G^y_1,\ldots,G^y_r)$ of $\torso{G}{B_y}$ has no major apex vertices (so $A^y=\emptyset$).  Since $B_x\cap B_y$ is a clique in $\torso{G}{B_y}$, $B_x\cap B_y\subseteq N_{\torso{G}{B_y}}[v]$.  Since $v$ is not in the top of $\torso{G}{B_y}$ and $\torso{G}{B_y}$ has no major apex vertices, $N_{\torso{G}{B_y}}[v]=N_{G^y_0}(v)$.  Since $v$ is not the near-top of $\torso{G}{B_y}$, $N_{G^y_0}[v]\cap (\hat{A}^y\cup V(G^y_1\cup\cdots\cup G^y_r))=\emptyset$.  Thus $B_x\cap B_y$ contains no vertices in the top of $\torso{G}{B_y}$, which is the first layer of $\LL_y$.  Thus, the construction so far already satisfies (3b) as well. 

    We now show that a slight modification of $\PP_x$ satisfies (3c) without greatly increasing the treewidth of $H_x:=\ltorso{G}{B_x}/\PP_x$.
    For each $x$ in $V(T)$ with parent $y$ in $V(T)$, each $v\in B_x\cap B_y$ is contained in the top of $\torso{G}{B_x}$ and is therefore contained in the first layer of $\LL_x$.  We now modify $\PP_x:=(L^x_1,L^x_2,\ldots)$ by breaking parts that intersect $B_x\cap B_y$ so that each $v\in B_x\cap B_y$ is in a singleton part, and each part is connected. Let $P$ be a part in $\PP_x$ that contains a vertex of $B_x\cap B_y$. By definition, $|P\cap L^x_2|\le \ell$. Since $G[P]$ is connected, each of the components $C_1,\ldots,C_r$ of $G[P]-(P\cap L^x_1)$ contains at least one vertex in $P\cap L^x_2$.  Therefore $G[P]-(P\cap L^x_1)$ has $r\le \ell$ connected components.  In $\PP_x$, replace $P$ with the at most $2\ell$ parts in $R_P:=\{V(C_1),\ldots,V(C_r)\}\cup\{\{v\}:v\in P\cap L^x_1\}$. In any tree-decomposition of $\ltorso{G}{B_x}/\PP_x$, the vertex that represents $P$ can be replaced by the vertices that represent each set in $R_P$. By definition, before this modification of $\PP_x$, $\ltorso{G}{B_x}/\PP_x$ had a tree-decomposition of width at most $\ell_0$.  Therefore, after this modification, $\ltorso{G}{B_x}/\PP_x$ has a tree-decomposition of width at most $\ell:= \ell_0(\ell_0+1)-1$.
\end{proof}

\Cref{apex_special_case} has several implications:
\begin{compactitem}
    \item For each $x\in V(T)$, $H_x:= \ltorso{G}{B_x}/\PP_x$ is a minor of $G$ (since $\ltorso{G}{B_x}/\PP_x$ is obtained by contracting connected subgraphs in a minor of $G$).

    \item For each $x\in V(T)$, $G[B_x]/\PP_x$ is a minor of $G$ (since $G[B_x]$ is a subgraph of $\ltorso{G}{B_x}$).

    \item $G[B_x]$ is isomorphic to a subgraph of $H_x\boxtimes P$, for some path $P$ and some minor $H_x$ of $G$ whose treewidth is at most $\ell$.
\end{compactitem}

\printindex

\end{document}